\documentclass[review]{elsarticle}
\usepackage{lineno,hyperref}
\usepackage{indentfirst,latexsym,bm}
\usepackage{mathrsfs}
\usepackage{CJK,fancybox,fancyhdr,color,graphicx,amsmath,amsthm,amssymb,enumerate,amsfonts}
\usepackage{anysize}
\usepackage[misc]{ifsym}
\usepackage[T1]{fontenc}
\usepackage{subfigure}
\usepackage[noend]{algpseudocode}
\usepackage{algorithmicx,algorithm}
\usepackage{xcolor}
\usepackage{colortbl,booktabs}

\newtheorem{thm}{Theorem}[section]
\newtheorem{cor}[thm]{Corollary}
\newtheorem{lem}[thm]{Lemma}

\newtheorem{prop}[thm]{Proposition}

\newtheorem{rem}[thm]{Remark}
\numberwithin{algorithm}{section}
\newcommand{\upcite}[1]{\normalfont\textsuperscript{\cite{#1}}}

\modulolinenumbers[5]

\journal{Journal of \LaTeX\ Templates}

\begin{document}

\begin{frontmatter}

\title{Convergence rates of the Kaczmarz-Tanabe method for linear systems}

\author[mymainaddress]{Chuan-gang Kang}\corref{mycorrespondingauthor}

\address[mymainaddress]{School of Mathematical Sciences, Tiangong University, Tianjin 300387,  Peoples R China}
\cortext[mycorrespondingauthor]{Chuan-gang Kang}
\ead{ckangtj@tiangong.edu.cn}

%
%



\begin{abstract}
  In this paper, we investigate the Kaczmarz-Tanabe method for exact and inexact linear systems. The Kaczmarz-Tanabe method is derived from the Kaczmarz method, but is more stable than that. We analyze the convergence and the convergence rate of the Kaczmarz-Tanabe method based on the singular value decomposition theory, and discover two important factors, i.e., the second maximum singular value of $Q$ and the minimum non-zero singular value of $A$, that influence the convergence speed and the amplitude of fluctuation of the Kaczmarz-Tanabe method (even for the Kaczmarz method). Numerical tests verify the theoretical results of the Kaczmarz-Tanabe method.\\

\noindent {\bf Keywords:} Kaczmarz-Tanabe method; Convergence rates; Singular value decompositon\\

\noindent {\bf Mathematics Subject Classification(2010)} 65F10 \sep 65F08 \sep 65N22 \sep 65J20
\end{abstract}
\end{frontmatter}

\section{Introduction}

\noindent The Kaczmarz method is one of the most popular iterative methods for image reconstruction in computerized tomography. It was proposed by the Polish mathematician Stefan Kaczmarz in \cite{Kaczmarz1937,Kaczmarz1993en}.
For the linear system with equations
\begin{equation}\label{linear.system}
  Ax=b,
\end{equation}
where $A=(a_{ij})\in \mathbb{R}^{m\times n}$ and $b\in \mathbb{R}^m$, and $x$ is an unknown vector. We denote the true solution with $x^*$ when \eqref{linear.system} is consistent. However, the linear problem \eqref{linear.system} may have no solution or multiple solutions, hence we are often asked to solve the minimum norm least-squares solution (i.e., Moore-Penrose generalized solution \cite{Engl96}) $x^\dagger$.

Let $A=(a_1,a_2,\ldots,a_m)^T$ and $b=(b_1,\ldots,b_m)^T$, then the classical form of Kaczmarz's algorithm \cite{Kaczmarz1937,Aboud2020} is described as
\begin{align}\label{kaczmarz.iteration}
  x_k=x_{k-1}+\frac{b_i-\langle a_i, x_{k-1}\rangle}{\|a_i\|_2^2}a_i,k=1,2,\ldots,
\end{align}
and the matrix-vector form is as follows,
\begin{align}\label{kaczmarz.iteration.2}
  x_k=(I-\frac{a_ia_i^T}{\|a_i\|_2^2})x_{k-1}+\frac{b_i}{\|a_i\|_2^2}a_i,k=1,2,\ldots,
\end{align}
where $i=(k\bmod m)+1$, $\langle x,y\rangle=x^Ty$ and $\|x\|_2=\sqrt{\langle x,x\rangle}$ denote the inner product and the square norm in $\mathbb{R}^n$, respectively. The iterative scheme of Kaczmarz's algorithm \eqref{kaczmarz.iteration} (or \eqref{kaczmarz.iteration.2}) sweeps through the equations of $Ax=b$ in a cyclic manner. In the first epoch, the processes of projecting the iterate $x_{k-1}$ orthogonally onto the solution hyperplane $\langle a_k, x\rangle=b_k$ and getting the new iterate $x_k$ are executed from $k=1$ until $k=m$. When $k\ge m$, then take $x_0=x_m$ and repeat the above process.

Kaczmarz's algorithm was unknown for more than 10 years after it was proposed \cite{Popa2018,Cegielski2015,Sznajder2015}. Until 1970, it was rediscovered as an algebraic reconstruction technique (ART) in computed tomography in \cite{Gordon1970}. Thereafter, K. Tanabe considered the Kaczmarz method in \cite{Tanabe1971} and investigated the convergence theory. He proved that the sequence of vectors generated by Kaczmarz's algorithm converges to the superposition of the Moore-Penrose solution and the orthogonal projection of the initial vector $x_0$ onto the null space $N(A)$.

T. Strohmer and R. Vershynin considered the randomized Kaczmarz method for the consistent and inconsistent linear systems and established the results of the exponential convergence rate in \cite{Strohmer2009}. In 2014, D. Needell  and J. A. Tropp considered block Kaczmarz's algorithm \cite{Needell2014paved} that used a random manner to pick up the projective subspace $\{x|A_\tau x=b_\tau\}$ at each step, where $\tau\in T=\{\tau_1,\ldots,\tau_r\}$ is a partition of the row indices of $A$.

Y. Jiao, B. Jin and X. Lu considered the preasymptotic convergence behavior of the randomized Kaczmarz method in \cite{Jiao2017preasymptotic} and illustrated its fast empirical convergence by analyzing the properties of the high- and low- frequency iterative errors.

K. Wei used the Kaczmarz method to solve systems of phaseless equation in \cite{Weike2016}, i.e., the generalized phase retrieval problem. He extended the Kaczmarz method for solving systems of linear equations by integrating a phase selection heuristic in each iteration. The preliminary convergence analysis has been presented for the randomized Kaczmarz methods.

C. Kang and H. Zhou considered the convergence of the Kaczmarz method \cite{Kang2021} and presented the convergence rate of the method for solving the exact and inexact linear systems based on the convergence theory in \cite{Tanabe1971}.

For further description, the following symbols will be used in this paper. The null and range spaces of $A$ will be denoted by $N(A)$ and $R(A)$, respectively. The rank of $A$ will be denoted by $\text{rank}(A)$. $\mathscr{S}^\bot$ will denote the orthogonal complement of a linear subspace $\mathscr{S}$. The Moore-Penrose generalized inverse \cite{Ben-Israel2003} of $A$ will be denoted by $A^\dagger$. The symbol $I$ will denote identity matrix of whatever size appropriate to the context. The transposition of a matrix or vector $X$ will be denoted by $X^T$. $\|A\|_2$ denotes the spectral norm of a matrix $A$ and is defined by
\begin{align*}
  \|A\|_2=\sup\limits_{x\in R^n\backslash\{0\}}\frac{\|Ax\|_2}{\|x\|_2}.
\end{align*}
Denote
\begin{align}\label{definition_AS_and_M}
  A_\mathcal{S}=(Q_1a_1,Q_2a_2,\ldots,Q_ma_m)^T, \quad M=\text{diag}(1/\|a_1\|_2^2,1/\|a_2\|_2^2,\ldots,1/\|a_m\|_2^2),
\end{align}
where,
\begin{align}\label{definition._pi_and_Qj}
  Q_j=P_mP_{m-1}\ldots P_{j+1} \quad(j=1,2,\ldots,m-1), \quad Q_m=I, \quad P_i=I-\frac{a_ia_i^T}{\|a_i\|_2^2} \quad(i=1,2,\ldots,m).
\end{align}
Some of these mathematical symbols were introduced by Tanabe in \cite{Tanabe1971}, so we try to quote these symbols in order to maintain their consistency, but there are still some difference, such as $Q_i(i=1,2,\ldots,m)$ defined in \eqref{definition._pi_and_Qj}, and $Q$ defined as follows
\begin{align}\label{definition_of_Q}
  Q=P_mP_{m-1}\ldots P_1.
\end{align}

K. Tanabe introduced the following results in \cite{Tanabe1971} and they are also valid for matrix $Q$ defined in \eqref{definition_of_Q}.
\begin{lem}\label{lemma1.1}  $Qx=x$ iff $x\in N(A)$.\end{lem}
\begin{lem}\label{lemma.2}\normalfont $\|Q\|_2\le 1$. If $\text{rank}(A)<n$ then $\|Q\|_2=1$.\end{lem}

\begin{thm}\label{thm.tilde.Q}\normalfont $\|\tilde{Q}\|_2=\sup\limits_{x\in R(A^T),\|x\|_2=1} \|Qx\|_2<1$, where $\tilde{Q}=QP_{R(A^T)}$ and $P_{R(A^T)}$ denotes the orthogonal projection onto the range space $R(A^T)$.\end{thm}

\begin{prop}\label{equality.result.2}\normalfont $I-A_\mathcal{S}^TMA=Q$.
\end{prop}
\proof \quad According to the definition of symbols \eqref{definition_AS_and_M} and \eqref{definition._pi_and_Qj}, we have
\begin{align*}
  Q&=P_m\ldots P_1=Q_1-Q_1\frac{a_1a_1^T}{\|a_1\|_2^2}=\ldots=Q_m-Q_m\frac{a_ma_m^T}{\|a_m\|_2^2}-\ldots-Q_1\frac{a_1a_1^T}{\|a_1\|_2^2}\\
   &=I-(Q_1a_1,Q_2a_2,\ldots,Q_ma_m)\text{diag}(1/\|a_1\|_2^2,1/\|a_2\|_2^2,\ldots,1/\|a_m\|_2^2)(a_1,a_2,\dots,a_m)^T\\
   &=I-A_\mathcal{S}^TMA.
\end{align*}

Proposition \ref{equality.result.2} was first introduced by K. Tanabe in 1971, and we redescribe the result because the definition of $Q_i$ here is somewhat different. The following theorem is derived from Theorem 8.1 in \cite{Tanabe1971}.

\begin{thm}\label{equality.result.1}\normalfont $(I-\tilde{Q})^{-1}A_{\mathcal{S}}^TM=A^\dagger$, where $\tilde{Q}$ is defined in Theorem \ref{thm.tilde.Q}.
\end{thm}

C. Kang and H. Zhou introduced the set-property (Lemma \ref{Range_result}) and the convergence rate (Theorem \ref{theorem.refer}) for the sequence of vectors $\{x_k\}$ generated by Kaczmarz's iteration \eqref{kaczmarz.iteration} in \cite{Kang2021}.
\begin{lem}\label{Range_result}\normalfont  For any $x_0\in R^n$, let $D_r=P_{N(A)}x_0+N(A)^{\bot}$, then for the sequence of vectors $\{x_k\}$ generated by Kaczmarz's iteration \eqref{kaczmarz.iteration} (or \eqref{kaczmarz.iteration.2}), there hold $x_k\in D_r, k=1,2,\ldots.$
\end{lem}

\begin{thm}\label{theorem.refer}\normalfont
Assume that \eqref{linear.system} is consistent and the sequence of vectors $\{x_k\}_{k=1}^m$ is generated by Kaczmarz's iteration \eqref{kaczmarz.iteration}. Then for any initial vector $x_0\in R^n$, there holds
\begin{align*}
   \|x_{k+1}-P_{N(A)}x_0-x^\dagger\|_2^2\le (1-\frac{1}{\|a_{k+1}\|_2^2\|(a_{k+1}^TP_{k+1})^\dagger\|_2^2})\|x_k-P_{N(A)}x_0-x^\dagger\|_2^2,
\end{align*}
where $\|(a_{k+1}^TP_{k+1})^\dagger\|_2\ge \|A^\dagger\|_2$.
\end{thm}

Theorem \ref{theorem.refer} was obtained based on Kaczmarz's iteration \eqref{kaczmarz.iteration}, and the coefficient factors on the right-hand side of the inequality are not easy to quantify. The iteration number of the Kaczmarz method is usually a integer multiple of the number of equations in order to make each equation work in the iterative algorithm. This feature allows us to extract a subsequence $\{y_k\}=\{x_{km}\}$ from the sequence of vectors $\{x_k\}$ and make it as a new iterative sequence.

Compared with the original sequence, the new iteration can be generated by the multiplication of matrix and vector and was named the Kaczmarz-Tanabe iteration in \cite{Popa2018}. The iterative scheme is described as follows,
\begin{align}\label{the_priodic_Kaczmarz_method}
  y_{k+1}=(I-A_{\mathcal{S}}^TMA)y_k+A_{\mathcal{S}}^TMb,\qquad k=0,1,2,\ldots.
\end{align}

In this paper, we consider the new iterative formula \eqref{the_priodic_Kaczmarz_method} of the Kaczmarz method to improve these convergence rate results in \cite{Kang2021}. The whole linear system are used in each iteration, therefore we can use the characteristic information about $A$ such as the condition number and the singular values, which can effectively avoid the difficulty to quantity the characteristic information of a standalone equation.

Our work is organized as follows. In Section \ref{section2}, we present the explicit form of the Kaczmarz-Tanabe method and consider its convergence and convergence rate for the exact linear system. In Section \ref{section3}, we consider the convergence rate of the Kaczmarz-Tanabe method for the linear system with perturbed right-hand side. In Section \ref{section4}, we present an sub-optimal algorithm, which can save much computational cost in forming $A_\mathcal{S}$ and $Q$. In Section \ref{section5}, we present some numerical tests to verify these theoretical results about convergence and convergence rate. Section \ref{section6} is the conclusion of this paper.

\section{The convergence rate of the Kaczmarz-Tanabe method for an exact linear system}\label{section2}

\noindent From Proposition \ref{equality.result.2}, the iterative formula \eqref{the_priodic_Kaczmarz_method} can also be described as
\begin{align}\label{iteration.scheme.2}
  y_{k+1}=Qy_k+A_{\mathcal{S}}^TMb,\qquad k=0,1,2,\ldots.
\end{align}
Obviously, there hold the following equalities,
\begin{align}\label{middle.formula}
  (I-A_\mathcal{S}^TMA)P_{N(A)}y_0^\delta=P_{N(A)}y_0^\delta, \quad (I-A_\mathcal{S}^TMA)x^\dagger=x^\dagger-A_\mathcal{S}^TMb.
\end{align}
Let $e_k=y_k-P_{N(A)}y_0-x^\dagger$ and $r_{k+1}=b-Ay_{k+1}$, then for any $k\ge 0$ there hold
\begin{align}\label{iterative.error.scheme}
  e_{k+1}=(I-A_\mathcal{S}^TMA)e_k
\end{align}
and
\begin{align}\label{residue.1}
  r_{k+1}=(I-AA_\mathcal{S}^TM)r_k.
\end{align}

\begin{lem}\label{range.result.2}\normalfont
  if $x\in N(A)$, then $Q_ix=x$ and $Q_i^Tx=x$, $i=1,2,\ldots,m$.
\end{lem}
\proof \quad From the condition $x\in N(A)$, we have $a_i^Tx=0, i=1,2,\ldots,m$. Therefore,
\begin{align*}
  Q_ix=Q_{i+1}x-Q_{i+1}\frac{a_{i+1}a_{i+1}^T}{\|a_{i+1}\|_2^2}x=Q_{i+1}x=\ldots=Q_mx=x.
\end{align*}
From $Q_i^T=P_{i+1}P_{i+2}\ldots P_m$ and $P_ix=x$, we also have the following equality.
\begin{align*}
  Q_i^Tx=P_{i+1}P_{i+2}\ldots P_mx=x.
\end{align*}

\begin{cor}\label{Kernal.result.3}\normalfont
  If $x\in N(A)$, then $A_\mathcal{S}x=0$ and $Q^Tx=x$.
\end{cor}

\begin{lem}\label{Range_result.3}\normalfont Assume that the sequence of vectors $\{y_k\}$ is generated by \eqref{iteration.scheme.2}, then $y_k\in D_r, k=1,2,\ldots.$
\end{lem}
\proof \quad For any vector $\tilde{x}\in N(A)$, there holds from Corollary \ref{Kernal.result.3} that
\begin{align*}
  \langle y_{k+1}-P_{N(A)}y_0,\tilde{x} \rangle
  &= \langle Qy_k+A_\mathcal{S}^TMb-P_{N(A)}y_0,\tilde{x} \rangle = \langle Qy_k,\tilde{x} \rangle +\langle A_\mathcal{S}^TMb,\tilde{x}\rangle-\langle P_{N(A)}y_0,\tilde{x}\rangle\\
  &=\langle y_k,Q^T\tilde{x}\rangle+\langle Mb,A_\mathcal{S}\tilde{x}\rangle-\langle P_{N(A)}y_0,\tilde{x}\rangle=\langle y_k-P_{N(A)}y_0,\tilde{x}\rangle.
\end{align*}
From the above recursive formula, it is obvious that
\begin{align*}
   \langle y_{k+1}-P_{N(A)}y_0,\tilde{x}\rangle =\langle y_0-P_{N(A)}y_0,\tilde{x}\rangle=0,
\end{align*}
which proves $y_{k+1}-P_{N(A)}y_0\in N(A)^\bot$, that is, $y_{k+1}\in D_r$ for any $k$.\qed
\begin{cor}\label{range.result4}\normalfont
  Under the condition of Lemma \ref{Range_result.3}, there hold $e_k\in N(A)^\bot, k=1,2,\ldots.$
\end{cor}

\begin{lem}\normalfont Let $W\in \mathbf{C}^{n\times n}$, and $W^H$ be its conjugate transpose matrix, then for any $x\in \mathbf{C}^n$ there hold
\begin{align*}
   \big\langle\frac{W+W^H}{2}x,x\big\rangle=\mathbf{Real}\langle Wx,x\rangle \quad \text{and} \quad \big\langle\frac{W-W^H}{2}x,x\big\rangle=\mathbf{Im}\big\langle Wx,x\big\rangle,
\end{align*}
where $\frac{W+W^H}{2}$ and $\frac{W-W^H}{2}$ are the symmetric and anti-symmetric parts of $W$, respectively.
\end{lem}

\begin{cor}\label{antisymmetric.cor}\normalfont Assume $\mathcal{W}\in R^{n\times n}$ is an anti-symmetric matrix, then there holds for any $x\in R^n$ that
\begin{align*}
 \langle\mathcal{W}x, x\rangle=0.
\end{align*}
\end{cor}
\begin{lem}\label{important.equality}\normalfont For the matrices $A_\mathcal{S}, M$ and $Q$ defined in \eqref{definition_AS_and_M} and \eqref{definition_of_Q}, there holds
\begin{align*}
   (A_\mathcal{S}^TMA)^T(2I-A_\mathcal{S}^TMA)=(I-Q^TQ)+2\cdot\frac{Q-Q^T}{2},
\end{align*}
where $(Q-Q^T)/2$ is the anti-symmetric part of the iterative matrix $Q$.
\end{lem}
\proof \quad From Proposition \ref{equality.result.2}, obviously there holds
\begin{align*}
   (A_\mathcal{S}^TMA)^T(2I-A_\mathcal{S}^TMA)=(I-Q)^T(I+Q)=(I-Q^TQ)+2\cdot\frac{Q-Q^T}{2}.
\end{align*}
\qed

The following theorem gives the relationship among the null space of $A, I-Q^TQ$ and $(I-Q^TQ)^\frac{1}{2}$.

\begin{thm}\label{kernal.space}\normalfont  $N(A)=N(I-Q^TQ)=N((I-Q^TQ)^\frac{1}{2})$.
\end{thm}

\proof \quad  If $\tilde{x}\in N(A)$, then from Corollary \ref{Kernal.result.3} there holds $\tilde{x}\subset N(I-Q^TQ)$. Conversely, when $\tilde{x}\in N(I-Q^TQ)$, if $\tilde{x}\notin N(A)$, then from the proof of Lemma 2 in \cite{Tanabe1971} we have $\|Q\|_2<1$, consequently $\|Q^TQ\tilde{x}\|_2\le \|Q\|_2^2\|\tilde{x}\|_2< \|\tilde{x}\|_2$, which is contrary to the condition of $\tilde{x}\in N(I-Q^TQ)$. This proves $N(I-Q^TQ)\subset N(A)$. The latter equality is obvious from $N((I-Q^TQ)^{1/2})\subset N(I-Q^TQ)$ and $\dim N((I-Q^TQ)^{1/2})=\dim N(I-Q^TQ)$.\qed

%
%
%

The following theorem gives the monotone result of the residues for the Kaczmarz-Tanabe method.
\begin{thm}\label{thm.monotone.residue.1}\normalfont  Let $\{r_k, k\ge 0\}$ be the sequence of vectors generated by \eqref{residue.1} and $\mathcal{L}=I-AA_\mathcal{S}^TM$, then there holds
\begin{align*}
  \|r_{k+1}\|_2^2=\|r_k\|_2^2-\langle(I-\mathcal{L}^T\mathcal{L})r_k,r_k\rangle.
\end{align*}
Moreover, if $\|\mathcal{L}\|_2\le 1$ then there also holds
\begin{align*}
  \|r_{k+1}\|_2^2\le\|r_k\|_2^2.
\end{align*}
\end{thm}
\proof \quad From \eqref{residue.1}, we have
\begin{align*}
  \|r_{k+1}\|_2^2=\|r_k\|_2^2-\langle MA_\mathcal{S}A^T(2I-AA_\mathcal{S}^TM)r_k,r_k\rangle,
\end{align*}
consequently,
\begin{align*}
  \|r_{k+1}\|_2^2&=\|r_k\|_2^2-\langle(I-\mathcal{L})^T(I+\mathcal{L})r_k,r_k\rangle\\
               &=\|r_k\|_2^2-\langle(I-\mathcal{L}^T\mathcal{L})r_k,r_k\rangle-2\langle\frac{\mathcal{L}-\mathcal{L}^T}{2}r_k,r_k\rangle.
\end{align*}
From Corollary \ref{antisymmetric.cor}, then it follows
\begin{align*}
  \|r_{k+1}\|_2^2=\|r_k\|_2^2-\langle(I-\mathcal{L}^T\mathcal{L})r_k,r_k\rangle.
\end{align*}
Especially, $I-\mathcal{L}^T\mathcal{L}$ is positive semi-definite when $\|\mathcal{L}\|_2\le 1$, thus $ \|r_{k+1}\|_2^2\le\|r_k\|_2^2$.\qed

%
\begin{thm}\label{convergence.error.free}\normalfont For any matrix $A$ with nonzero rows and any $m$ dimensional column vector $b$, let $\{y_k, k\ge 0\}$ be the sequence of vectors generated by \eqref{the_priodic_Kaczmarz_method}, then there hold
\begin{align}\label{important.result.1}
\|e_{k+1}\|_2^2\le\mathcal{K}\|e_k\|_2^2
\end{align}
and
\begin{align}\label{important.result.2}
  \|e_{k+1}\|_2^2\le\mathcal{K}^{k+1}\|e_0\|_2^2,
\end{align}
where $\mathcal{K}=1-\min\limits_{i=1,\ldots,p}\{1-\sigma_i^2,1\}$ and $\sigma_i$ is the singular value of $Q$.
\end{thm}

\proof \quad From \eqref{iterative.error.scheme}, we have
\begin{align*}
  \|e_{k+1}\|_2^2&=\langle(I-A_\mathcal{S}^TMA)e_k,(I-A_\mathcal{S}^TMA)e_k\rangle=\|e_k\|_2^2-\langle (A_\mathcal{S}^TMA)^T(2I-A_\mathcal{S}^TMA)e_k,e_k\rangle.
\end{align*}
From Lemma \ref{important.equality}, consequently,
\begin{align}\label{important.equality.1}
  \|e_{k+1}\|_2^2=\|e_k\|_2^2-\langle(I-Q^T)(I+Q)e_k,e_k\rangle=\|e_k\|_2^2-\langle(I-Q^TQ)e_k,e_k\rangle-2\langle\frac{Q-Q^T}{2}e_k,e_k\rangle.
\end{align}
From Lemma \ref{antisymmetric.cor} and \eqref{important.equality.1}, it follows that
\begin{align}\label{error.formula.1}
  \|e_{k+1}\|_2^2=\|e_k\|_2^2-\langle(I-Q^TQ)e_k,e_k\rangle.
\end{align}
We assume the singular value decomposition \cite{Ben-Israel2003} of $Q$ as follows,
\begin{align*}
  Q=U\Sigma V^T,\quad U=(U_1,U_2), \quad V=(V_1,V_2),
\end{align*}
where $U$ and $V$ are orthogonal matrices of order $m$ and $n$, respectively. $U_1$ and $U_2$ are $m\times p$ and $m\times (n-p)$ matrices, $V_1$ and $V_2$ are $n\times p$ and $n\times (n-p)$ matrices, respectively, and
\begin{align*}
  \Sigma=\text{diag}(\sigma_1,\sigma_2,\cdots,\sigma_p,0,\cdots,0)\in R^{m\times n},
\end{align*}\label{kernal.equality.2}
where $\sigma_1\ge\sigma_2\ge\cdots\ge\sigma_p>0$ and $\text{rank}(Q)=p$. From Lemma \ref{lemma.2} and $\sqrt{\lambda_{\max}(Q^TQ)}=\|Q\|_2\le 1$, therefore $I-Q^TQ$ is positive semi-definite, and it follows
\begin{align}\label{important.formular}
  \|e_{k+1}\|_2^2=\|e_k\|_2^2-\langle(I-Q^TQ)e_k,e_k\rangle=\|e_k\|_2^2-\|(I-Q^TQ)^{1/2}e_k\|_2^2.
\end{align}
Notice that
\begin{align*}
Q^TQ=V(\Sigma^T\Sigma)V^T=V\left [
  \begin{array}{cc}
    \bar{\Sigma} &0\\
    0&0
  \end{array}
  \right ]V^T,\quad
  I-Q^TQ=V(I-\Sigma^T\Sigma)V^T=V\left [
  \begin{array}{cc}
    I-\bar{\Sigma}^2 &0\\
    0&I
  \end{array}
  \right ]V^T,
\end{align*}
where $\bar{\Sigma}=\text{diag}(\sigma_1,\sigma_2,\cdots,\sigma_p)$, and the diagonal elements of $\bar{\Sigma}$ fall into the interval $(0,1]$, which means $I-\bar{\Sigma}^2$ is positive semi-definite, therefore,
\begin{align*}
  (I-Q^TQ)^{1/2}=V\left [
  \begin{array}{cc}
    (I-\bar{\Sigma}^2)^{\frac{1}{2}} &0\\
    0&I
  \end{array}
  \right ]V^T.
\end{align*}
Consequently,
\begin{align*}
  \big((I-Q^TQ)^{1/2}\big)^\dagger=V\left [
  \begin{array}{cc}
    \big((I-\bar{\Sigma}^2)^{\frac{1}{2}}\big)^\dagger &0\\
    0&I
  \end{array}
  \right ]V^T,
\end{align*}
and
\begin{align*}
  \big((I-\bar{\Sigma}^2)^{\frac{1}{2}}\big)^\dagger=\text{diag}(\lambda_1,\lambda_2,\ldots,\lambda_p),
\end{align*}
where
\begin{align*}
  \lambda_i=\left \{
    \begin{array}{ll}
       (1-\sigma_i^2)^{-\frac{1}{2}}, &\quad \text{if}\quad|\sigma_i|<1,\\
       0, &\quad \text{if}\quad |\sigma_i|=1.
    \end{array}
  \right.
\end{align*}
So
\begin{align}\label{importan.norm}
  \|\big((I-Q^TQ)^{1/2}\big)^\dagger\|_2=\max\limits_{i=1,2,\ldots,p}\{(1-\sigma_i^2)^{-\frac{1}{2}},1\}.
\end{align}
From Corollary \ref{range.result4} and Theorem \ref{kernal.space}, there holds $e_k\in N((I-Q^TQ)^{1/2})^\bot$. And from \eqref{importan.norm}, it follows that
\begin{align}\label{important.iequality}
   \|(I-Q^TQ)^{1/2}e_k\|_2\ge \frac{1}{\|\big((I-Q^TQ)^{1/2}\big)^\dagger\|_2}\|e_k\|_2=\min\limits_{i=1,2,\ldots,p}\{(1-\sigma_i^2)^{\frac{1}{2}},1\}\|e_k\|_2.
\end{align}
Then there holds from \eqref{important.formular} and \eqref{important.iequality} that
\begin{align*}
  \|e_{k+1}\|_2^2\le\big(1-\min\limits_{i=1,2,\ldots,p}\{(1-\sigma_i^2),1\}\big)\|e_k\|_2^2,
\end{align*}
and
\begin{align*}
  \|e_{k+1}\|_2^2\le\big(1-\min\limits_{i=1,2,\ldots,p}\{(1-\sigma_i^2),1\}\big)^{k+1}\|e_0\|_2^2.
\end{align*}
From the definition of symbol $\mathcal{K}$, \eqref{important.result.1} and \eqref{important.result.2} are proved. \hfill{$\square$}

Theorem \ref{convergence.error.free} presents the results of convergence rate for the Kaczmarz-Tanabe method. Furthermore, we also have the following result from Lemmas \ref{lemma1.1}, \ref{lemma.2} and Corollary \ref{range.result4}.
\begin{cor}\label{corollary_for_convergence_rate}\normalfont Under the conditions of Theorem \ref{convergence.error.free}, there also hold
\begin{align*}
\|e_{k+1}\|_2^2\le\bar{\mathcal{K}}\|e_k\|_2^2 \quad \text{and}\quad\|e_{k+1}\|_2^2\le\bar{\mathcal{K}}^{k+1}\|e_0\|_2^2
\end{align*}
where $\bar{\mathcal{K}}=1-\min\limits_{i=1,\ldots,p}\{1-\sigma_i^2,1\}\in [0,1)$, and $\sigma_i$ is the singular value that less than $1$ of $Q$.
\end{cor}
From Corollary \ref{corollary_for_convergence_rate}, we can obtain the convergence of the Kaczmarz-Tanabe method. Meanwhile, Corollary \ref{corollary_for_convergence_rate} is also equivalent to Corollary 9 in \cite{Tanabe1971}. Moreover, the convergence speed of the Kaczmarz-Tanabe method is closely related to the second maximum singular value of Q.

\section{The convergence rate of the Kaczmarz-Tanabe method for an inexact linear system}\label{section3}

\noindent If there is perturbation on the right-hand side of \eqref{linear.system}, i.e.,
\begin{align}\label{linear.system.perturbation}
  Ax=b^\delta,
\end{align}
where $b^\delta=b+\delta b$, then Kaczmarz-Tanabe's iteration \eqref{the_priodic_Kaczmarz_method} (or \eqref{iteration.scheme.2}) can be described as
\begin{align}\label{iteration.noise.1}
  y_{k+1}^\delta=(I-A_{\mathcal{S}}^TMA)y_k^\delta+A_{\mathcal{S}}^TMb^\delta=Qy_k^\delta+A_{\mathcal{S}}^TMb^\delta,\qquad k=0,1,2,\ldots.
\end{align}
Moreover it follows from \eqref{middle.formula} and \eqref{iteration.noise.1} that
\begin{align*}
  y_{k+1}^\delta-P_{N(A)}y_0^\delta-x^\dagger&=(I-A_\mathcal{S}^TMA)(y_k^\delta-P_{N(A)}y_0^\delta-x^\dagger)+A_\mathcal{S}^TM(b^\delta-b).
\end{align*}
Let $e_k^\delta=y_k^\delta-P_{N(A)}y_0^\delta-x^\dagger$ and $r_k^\delta=b^\delta-Ay_k^\delta$, consequently there also hold
\begin{align}\label{error.formula.2}
  e_{k+1}^\delta=(I-A_\mathcal{S}^TMA)e_k^\delta+A_\mathcal{S}^TM(b^\delta-b)=Qe_k^\delta+A_\mathcal{S}^TM(b^\delta-b)
\end{align}
and
\begin{align}\label{general_residue}
  r_{k+1}^\delta=(I-AA_\mathcal{S}^TM)r_k^\delta.
\end{align}

Theorem \ref{thm.monotone.residue.2} gives the monotonicity of the residues for the Kaczmarz-Tanabe method to solve \eqref{linear.system.perturbation}.
\begin{thm}\label{thm.monotone.residue.2}\normalfont Let $\{r_k^\delta, k\ge 0\}$ be the sequence of vectors generated by \eqref{general_residue} and $\mathcal{L}=I-AA_\mathcal{S}^TM$, then there holds
\begin{align*}
  \|r_{k+1}^\delta\|_2^2=\|r_k^\delta\|_2^2-\langle(I-\mathcal{L}^T\mathcal{L})r_k^\delta,r_k^\delta\rangle.
\end{align*}
Moreover, if $\|\mathcal{L}\|_2\le 1$ there also holds
\begin{align*}
  \|r_{k+1}^\delta\|_2^2\le\|r_k^\delta\|_2^2.
\end{align*}
\end{thm}
\begin{rem}\normalfont The results can be derived from the procedure of the proof to Theorem \ref{thm.monotone.residue.1}.\end{rem}

The following theorem gives the results of the convergence rate for the Kaczmarz-Tanabe method to solve the inexact linear system \eqref{linear.system.perturbation}.


\begin{thm}\label{convergence.rate.perturbation}\normalfont For any matrix $A$ with nonzero rows and any $m$-dimensional column vector $b^\delta$, let $\{y_k^\delta, k\ge 0\}$ be the sequence of vectors generated by \eqref{iteration.noise.1}, then there hold
\begin{align}\label{formular.estimation.1}
  \|e_{k+1}^\delta\|_2\le \mathcal{K}^{1/2}\|e_k^\delta\|_2+2\|A^\dagger\|_2\|b^\delta-b\|_2
\end{align}
and
\begin{align}\label{formular.estimation.2}
  \|e_{k+1}^\delta\|_2\le \mathcal{K}^{(k+1)/2}\|e_0^\delta\|_2+4(1-\mathcal{K}^{(k+1)/2})\max\limits_{i=1,\ldots,p}\{(1-\sigma_i^2)^{-1},1\}\|A^\dagger\|_2\|b^\delta-b\|_2,
\end{align}
where $\mathcal{K}$ and $\sigma_i$ are defined in Theorem \ref{convergence.error.free}.

\end{thm}
\proof \quad From \eqref{error.formula.2}, we have
\begin{align*}
  \|e_{k+1}^\delta\|_2^2
  &=\langle (I-A_\mathcal{S}^TMA)e_k^\delta+A_\mathcal{S}^TM(b^\delta-b),(I-A_\mathcal{S}^TMA)e_k^\delta+A_\mathcal{S}^TM(b^\delta-b)\rangle\\
  &=\|(I-A_\mathcal{S}^TMA)e_k^\delta\|_2^2+2\langle A_\mathcal{S}^TM(b^\delta-b),(I-A_\mathcal{S}^TMA)e_k^\delta\rangle+\langle A_\mathcal{S}^TM(b^\delta-b),A_\mathcal{S}^TM(b^\delta-b)\rangle\\
  &\le \big(\|(I-A_\mathcal{S}^TMA)e_k^\delta\|_2+\|A_\mathcal{S}^TM(b^\delta-b)\|_2\big)^2.
\end{align*}
Therefore,
\begin{align}\label{recursive.formula}
  \|e_{k+1}^\delta\|_2\le\|(I-A_\mathcal{S}^TMA)e_k^\delta\|_2+\|A_\mathcal{S}^TM(b^\delta-b)\|_2.
\end{align}
From \eqref{error.formula.1} and Lemma \ref{lemma.2}, then it follows from \eqref{recursive.formula} that
\begin{align}\label{important.inequality.perturbation}
  \|e_{k+1}^\delta\|_2\le \Big(1-\frac{1}{\|\big((I-Q^TQ)^{1/2}\big)^\dagger\|_2^2}\Big)^{1/2}\|e_k^\delta\|_2+\|A_\mathcal{S}^TM(b^\delta-b)\|_2.
\end{align}
Then from Theorems \ref{equality.result.1} and \ref{thm.tilde.Q}, we have $A_{\mathcal{S}}^TM=(I-\tilde{Q})A^\dagger$ and $\|\tilde{Q}\|_2\le 1$, consequently,
\begin{align}\label{important.formula.1}
  \|A_{\mathcal{S}}^TM\|_2=\|(I-\tilde{Q})A^\dagger\|_2\le\|(I-\tilde{Q})\|_2\|A^\dagger\|_2\le(1+\|\tilde{Q}\|_2)\|A^\dagger\|_2\le 2\|A^\dagger\|_2.
\end{align}
Hence from \eqref{importan.norm}, \eqref{important.inequality.perturbation} and \eqref{important.formula.1}, there holds
\begin{align*}
  \|e_{k+1}^\delta\|_2\le \big(1-\min\limits_{i=1,2,\ldots,p}\{1-\sigma_i^2,1\}\big)^{1/2}\|e_k^\delta\|_2+2\|A^\dagger\|_2\|b^\delta-b\|_2,
\end{align*}
which proves \eqref{formular.estimation.1}. Then from the recursion of \eqref{formular.estimation.1} there holds
\begin{align*}
  \|e_{k+1}^\delta\|_2&\le\mathcal{K}^{(k+1)/2}\|e_0^\delta\|_2+2\big(1-\mathcal{K}^{(k+1)/2}\big)(2-\mathcal{K})\cdot\max\limits_{i=1,\ldots,p}\{(1-\sigma_i^2)^{-1},1\}\|A^\dagger\|_2\|b^\delta-b\|_2\\
                      &\le\mathcal{K}^{(k+1)/2}\|e_0^\delta\|_2+4\big{(}1-\mathcal{K}^{(k+1)/2}\big{)}\cdot\max\limits_{i=1,\ldots,p}\{(1-\sigma_i^2)^{-1},1\}\|A^\dagger\|_2\|b^\delta-b\|_2.
\end{align*}
Thus \eqref{formular.estimation.2} is proved.\qed

Theorem \ref{convergence.rate.perturbation} presents the convergence rates of the Kaczmarz-Tanabe method for a perturbed linear system. From \eqref{formular.estimation.2}, the error of each iteration consists of two parts, i.e., the iterative error and the perturbed error. The iterative error is decreasing and the perturbed error is increasing about iteration number $k$, they are determined by the factor $\big{(}1-\min\limits_{i=1,\ldots,p}\{1-\sigma_i^2,1\}\big{)}^{(k+1)/2}$. However, the amplitude of the fluctuation is nevertheless related with the quantity $\max\limits_{i=1,\ldots,p}\{(1-\sigma_i^2)^{-1},1\}\|A^\dagger\|_2$. Moreover, the iterative error dominates the iterations to converge at the previous stage, and the perturbed error dominates the iterations to diverge at the latter stage. The more general conclusion is that the faster the iterations converge at the former stage and the faster the iterations diverge at the latter stage.

\section{The algorithm of the Kaczmarz-Tanabe method}\label{section4}

\noindent In this section, we present the algorithm of the Kaczmarz-Tanabe method in Algorithm \ref{The.periodic.Kaczmarz.solver} which can avoid a lot of repeated calculation. Algorithm \ref{The.periodic.Kaczmarz.solver} is close to optimal without considering parallelism. The calculation cost of the Kaczmarz-Tanabe method mainly comes from the generation of $A_{\mathcal{S}}$ and $Q$ (that is in the Process I of Algorithm \ref{The.periodic.Kaczmarz.solver}), and the cost of Process II is very small. Especially, once $A_\mathcal{S}$ and $Q$ are generated, then they can be used repeatedly for more scenes. However, the Kaczmarz method doesn't have this advantage.

\begin{algorithm}[t]
\caption{The Kaczmarz-Tanabe method solver}\label{The.periodic.Kaczmarz.solver}
\begin{algorithmic}[0]
\State Step 1. Given $A,b^\delta$, an initial guess $y_0^\delta$ and the maximum iteration number $K_{\max}$;
\State Step 2. generate the diagonal matrix $M$;
\State Step 3. $[m,n]=size(A)$;
\end{algorithmic}
\hspace*{0.1in}{\bf Process I. Compute $A_\mathcal{S}$ and $Q$}
\begin{algorithmic}[0]
\State Step 4. $i=1; Q=I_n$;
\State Step 5. $A_\mathcal{S}(:,i)=Q*A(i,:)^T$;
\State Step 6. $Q=Q-Q*A(i,:)^T*A(i,:)/\|A(i,:)\|_2^2$;
\State Step 7. $i=i+1$, if $i\le m$, go to Step 5; otherwise, go to step 8.
\end{algorithmic}
\hspace*{0.1in}{\bf process II. Perform Kaczmarz-Tanabe's iteration}
\begin{algorithmic}[0]
\State Step 8. $k=0, Vc=A_\mathcal{S}*M*b^\delta$;
\State Step 9. $y_{k+1}^\delta=Q*y_k^\delta+Vc$;
\State Step 10. $k=k+1$, if $k\le K_{\max}$, go to Step 9; otherwise, output the numerical solution $y_k^\delta$.
\end{algorithmic}
\end{algorithm}

%
%
%
%

\section{Numerical tests}\label{section5}
\subsection{Model Problem 1}
\begin{align}\label{Tanabe.problem.1}
  \left (
     \begin{array}{rrrr}
         1.0 &  3.0 &  2.0 & -1.0\\
         1.0 &  2.0 & -1.0 & -2.0\\
         1.0 & -1.0 &  2.0 &  3.0\\
         2.0 &  1.0 &  1.0 &  1.0\\
         5.0 &  5.0 &  4.0 &  1.0\\
         4.0 & -1.0 &  5.0 &  7.0
      \end{array}
  \right  )x=
  \left (
     \begin{array}{r}
        5.0\\
        0.0\\
        5.0\\
        5.0\\
       15.0\\
       15.0
     \end{array}
  \right )
\end{align}

Model Problem 1\upcite{Tanabe1971} is a consistent over-determined linear system, the true solution is $x=(1,1,1,1)^T$. For the perturbed case, we take
the right-hand side
\begin{align}\label{perturbation_definition}
  b_i^\delta=b_i+\delta \max\limits_{i}(|b_i|),\qquad i=1,\ldots,5
\end{align}
in \eqref{linear.system.perturbation} with $\delta=0.1$ and $0.3$, respectively. The absolute and relative errors of the right-hand side for Model Problem 1 are listed in Table \ref{errror_level_for_model_problem1}.
\begin{table}[!htb]
  \centering
 \caption{The absolute and relative errors of the right-hand side for Model Problem 1}
 \setlength{\tabcolsep}{1.5cm}{
  \begin{tabular}{ccc}
    \hline
    \text{Menu}                      &\text{$\delta=0.1$} &\text{$\delta=0.3$}  \\
    \hline
    $\|b^\delta-b\|_2$               & 3.6742             & 11.0227                \\
    $\|b^\delta-b\|_2/\|b\|_2$       & 0.1604             &  0.4811                  \\
    \hline
  \end{tabular}}
\label{errror_level_for_model_problem1}
\end{table}

\begin{figure}[hbt]
  \centering
  \subfigure[$K=100,\delta=0$]{
  \begin{minipage}[ht]{.32\linewidth}
      \includegraphics[width=\linewidth]{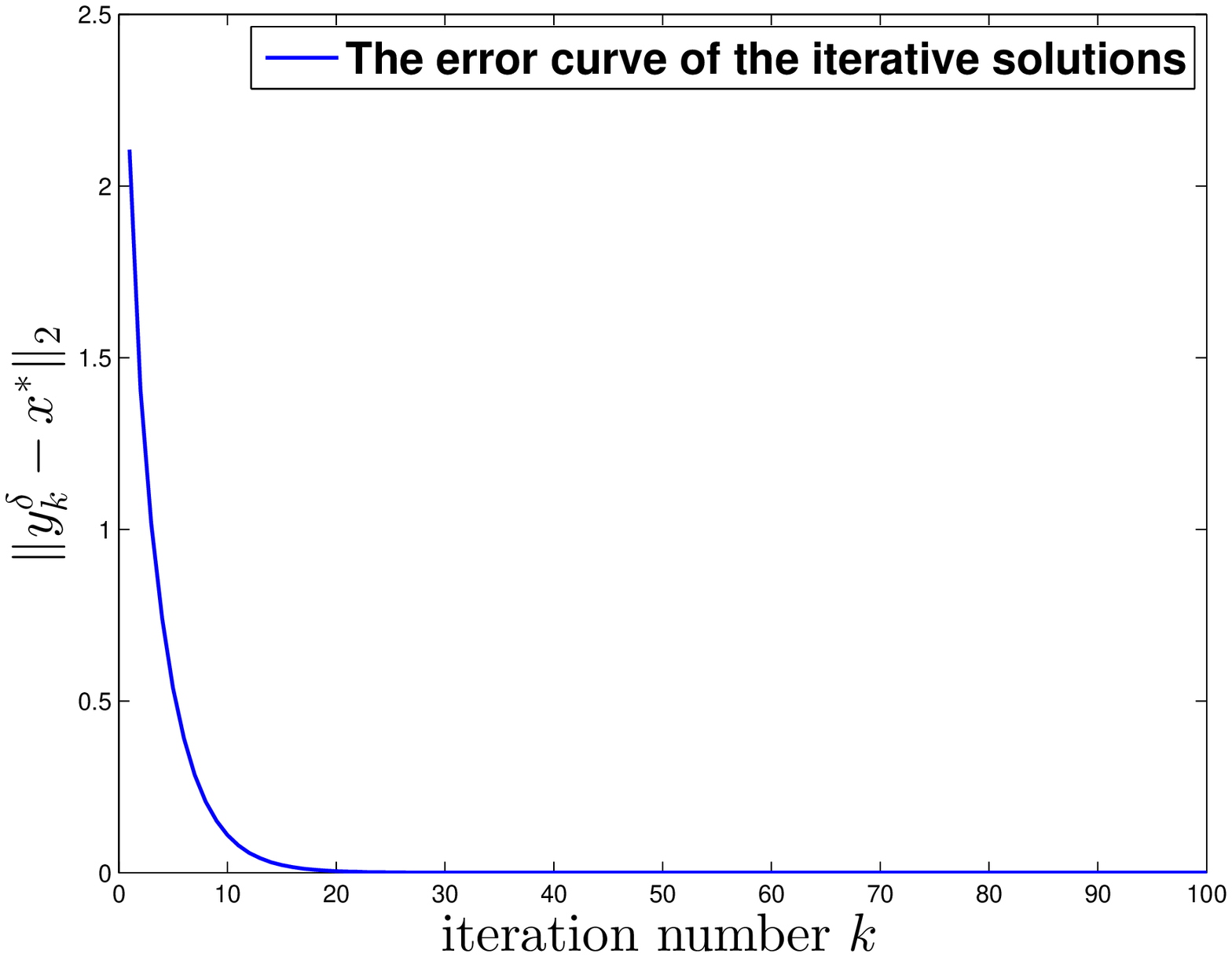}
  \end{minipage}}
  \subfigure[$K=100,\delta=0.1$]{
  \begin{minipage}[h]{.32\linewidth}
     \includegraphics[width=\linewidth]{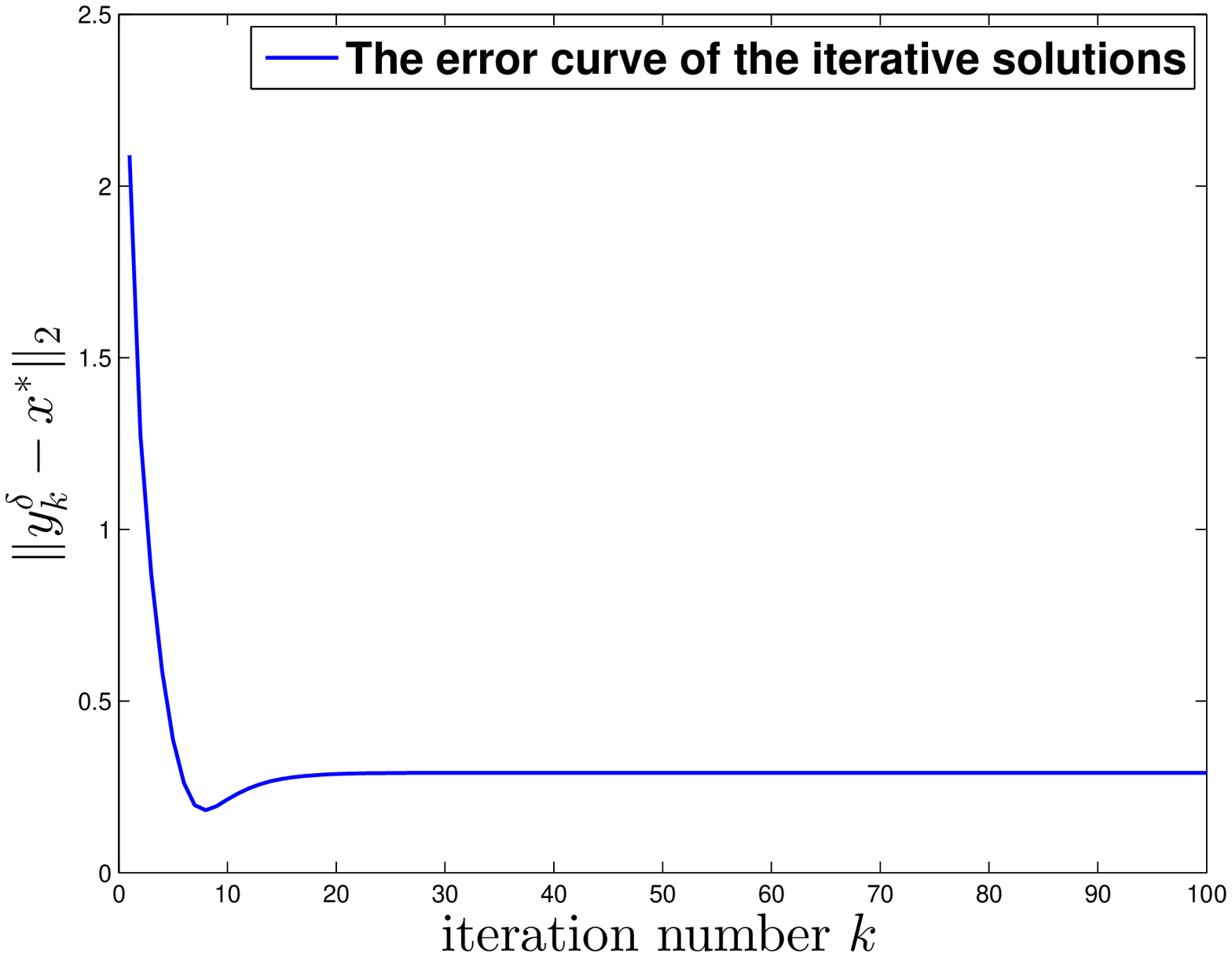}
    \end{minipage}}
  \subfigure[$K=100,\delta=0.3$]{
  \begin{minipage}[h]{.32\linewidth}
     \includegraphics[width=\linewidth]{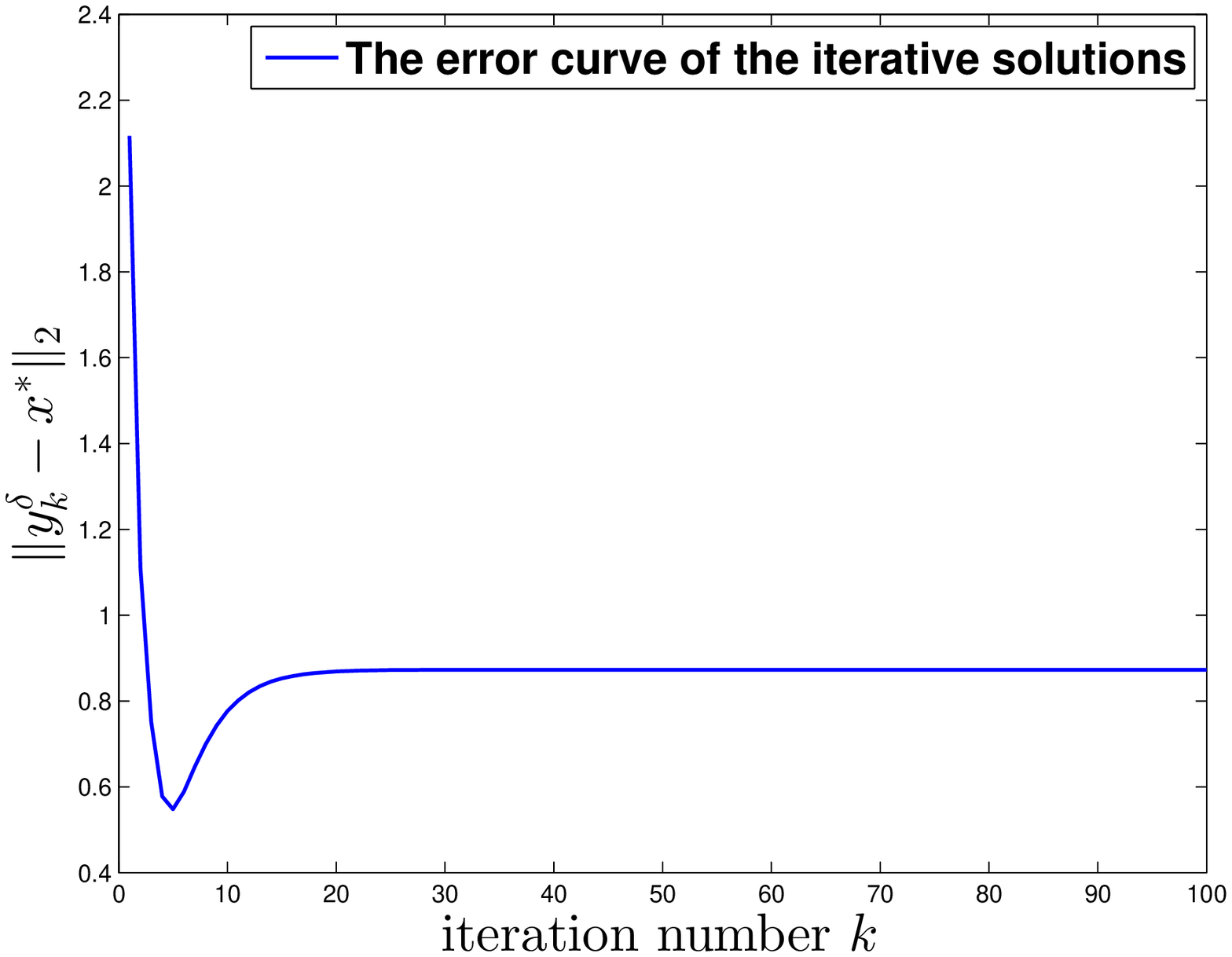}
  \end{minipage}}
       \caption{\footnotesize The iterative error curves of the Kaczmarz-Tanabe method for Model Problem 1}\label{Tanabe.error.curve}
\end{figure}

\begin{figure}[hbt]
  \centering
  \subfigure[$K=100,\delta=0$]{
  \begin{minipage}[ht]{.32\linewidth}
      \includegraphics[width=\linewidth]{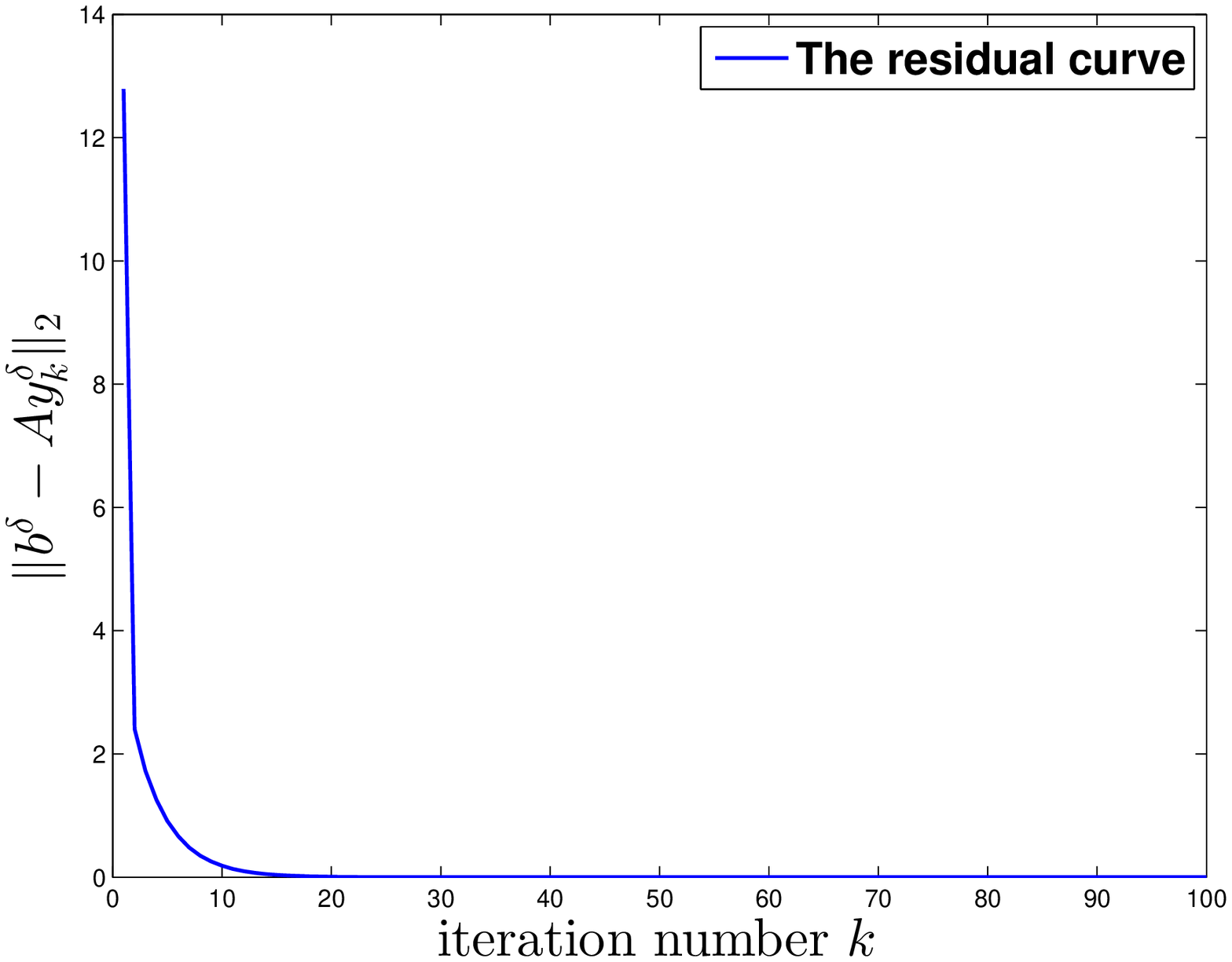}
  \end{minipage}}
  \subfigure[$K=100,\delta=0.1$]{
  \begin{minipage}[h]{.32\linewidth}
     \includegraphics[width=\linewidth]{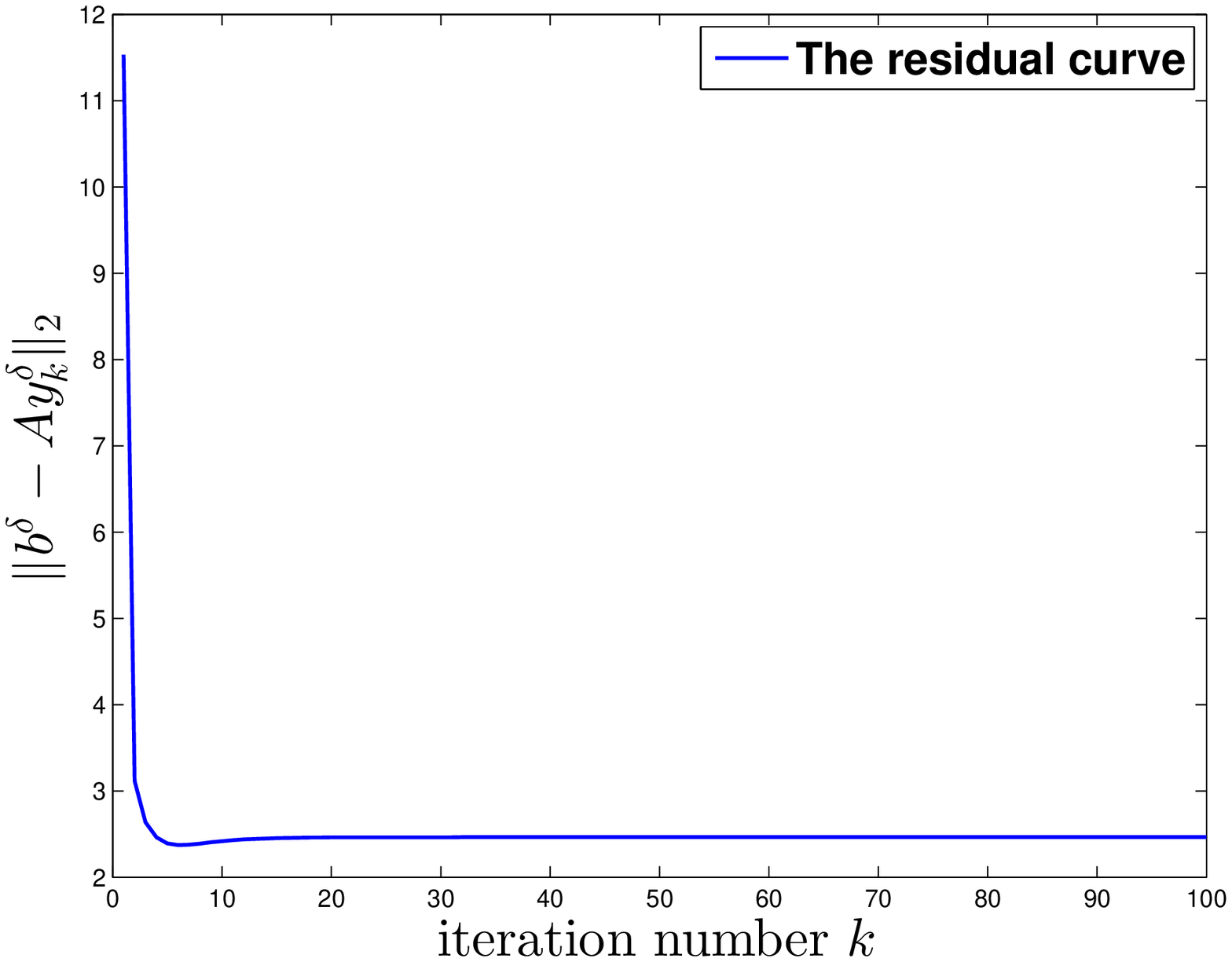}
    \end{minipage}}
  \subfigure[$K=100,\delta=0.3$]{
  \begin{minipage}[h]{.32\linewidth}
     \includegraphics[width=\linewidth]{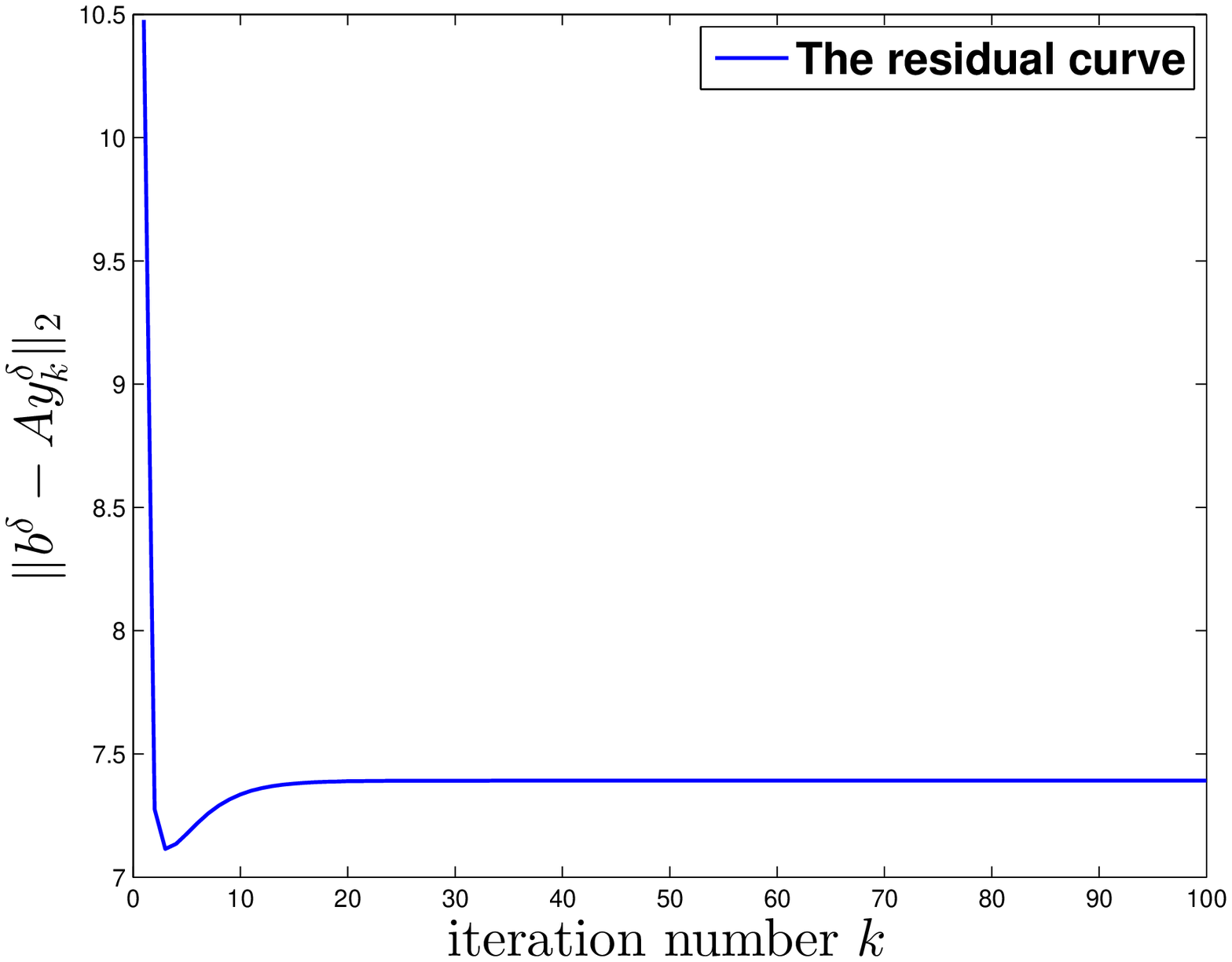}
  \end{minipage}}
       \caption{\footnotesize The residual curves of the Kaczmarz-Tanabe method for Model Problem 1}\label{Tanabe.residue.curve}
\end{figure}

Numerical results are presented in Figures \ref{Tanabe.error.curve} and \ref{Tanabe.residue.curve}. It is easy to see, from Figures \ref{Tanabe.error.curve}(a) and \ref{Tanabe.residue.curve}(a), that the iterative error and residual curves are decreasing monotonically for the exact case. However, all of these curves fluctuate under disturbance, that is, they are not monotonically decreasing.

\subsection{Model Problem 2}
\noindent Considering the following differential equation \cite{Elman82,Kamath88},
\begin{align}
  &-(e^{-xy}u_x)_x-(e^{xy}u_y)_y+\beta(x+y)u_y+[\beta(x+y)u]_y+\frac{1}{1+x+y}u=g,\quad (x,y)\in [0,1]\times [0,1],\label{partial_differential_equation}\\
  &u=xe^{xy}\sin(\pi x)\sin(\pi y), \beta=10000.0.\label{true_solution}
\end{align}
It is a classical example and often used to illustrate the parallel methods because its discrete form is a diagonal form linear system \cite{Kamath88}.

In this problem, the right-hand function $g(x,y)$ needs to be calculated from the true solution $u(x,y)$ in \eqref{true_solution}, and its analytic expression is
\begin{align*}
   g(x,y)=e^{-xy}(yu_x-u_{xx})-e^{xy}(u_{yy}+xu_y)+2\beta(x+y)u_y+(\beta+\frac{1}{1+x+y})u.
\end{align*}
The linear system can be derived from the discretization of the partial differential equation \eqref{partial_differential_equation} by the difference method. For this problem, we take the dimension of the discretization $n=32$. The perturbed right-hand side is defined as \eqref{perturbation_definition}

\begin{figure}[hbt]
  \centering
  \subfigure[$K=200,\delta=0$]{
  \begin{minipage}[ht]{.32\linewidth}
      \includegraphics[width=\linewidth]{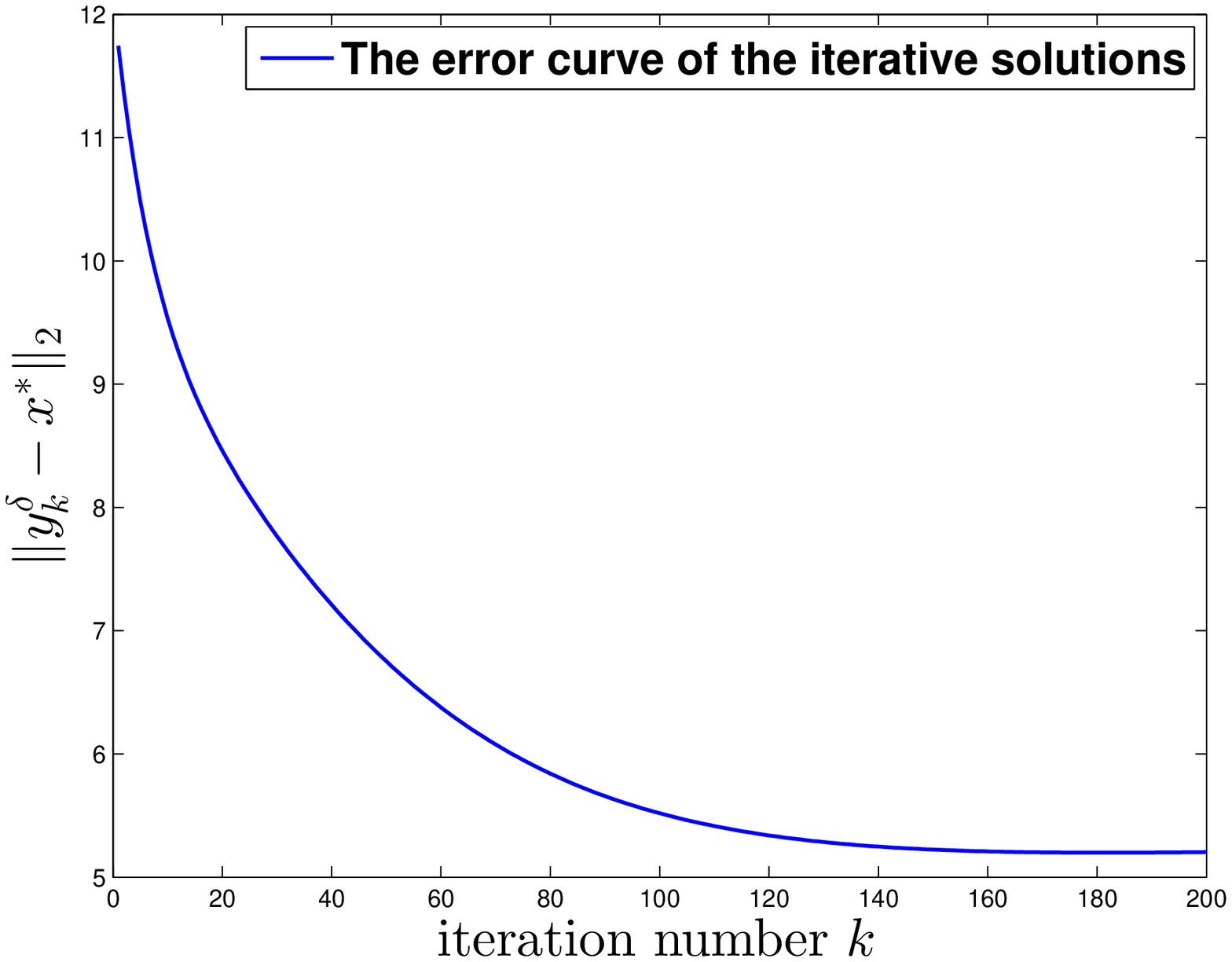}
  \end{minipage}}
  \subfigure[$K=200,\delta=0.1$]{
  \begin{minipage}[ht]{.32\linewidth}
      \includegraphics[width=\linewidth]{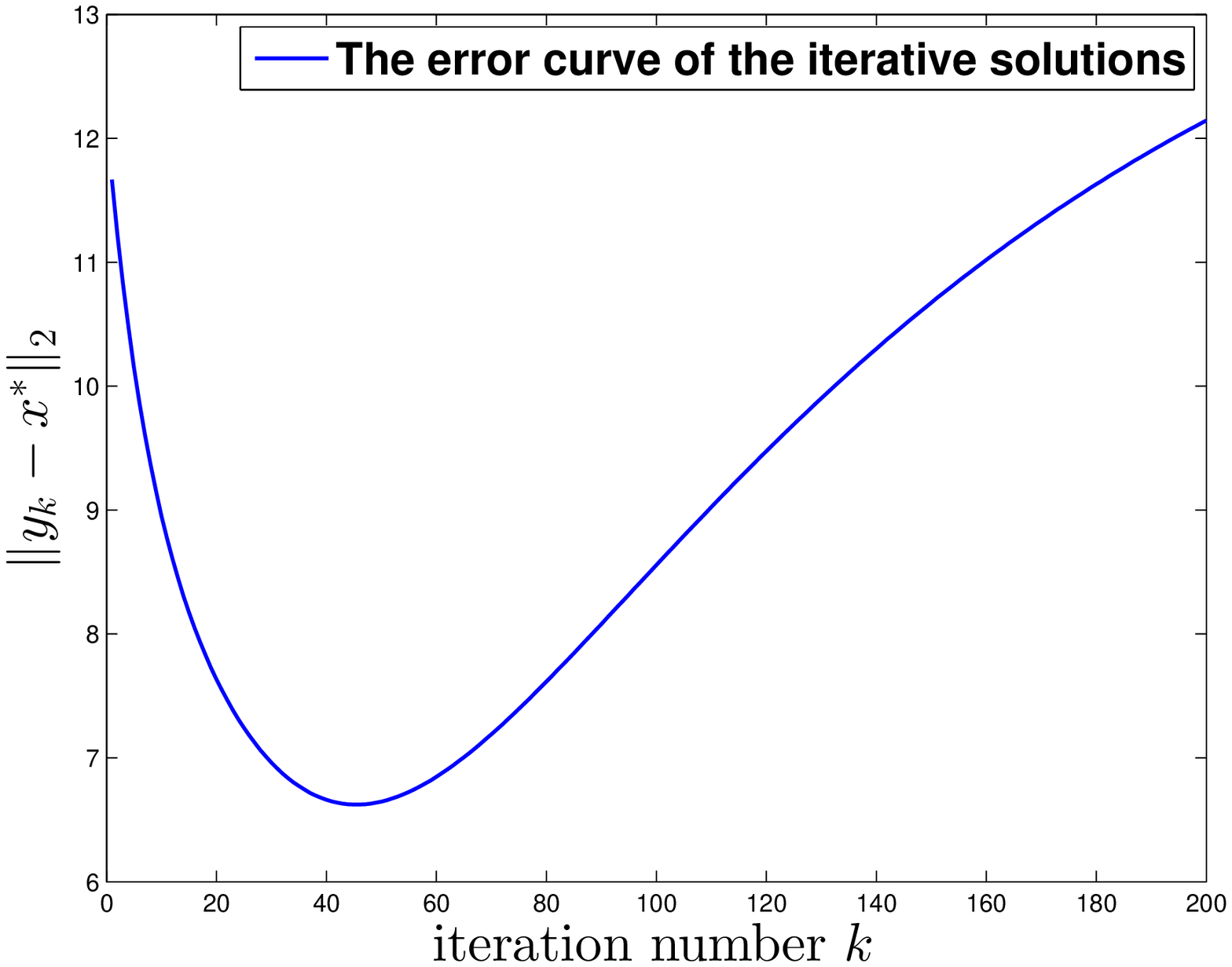}
  \end{minipage}}
  \subfigure[$K=200,\delta=0.3$]{
  \begin{minipage}[ht]{.32\linewidth}
      \includegraphics[width=\linewidth]{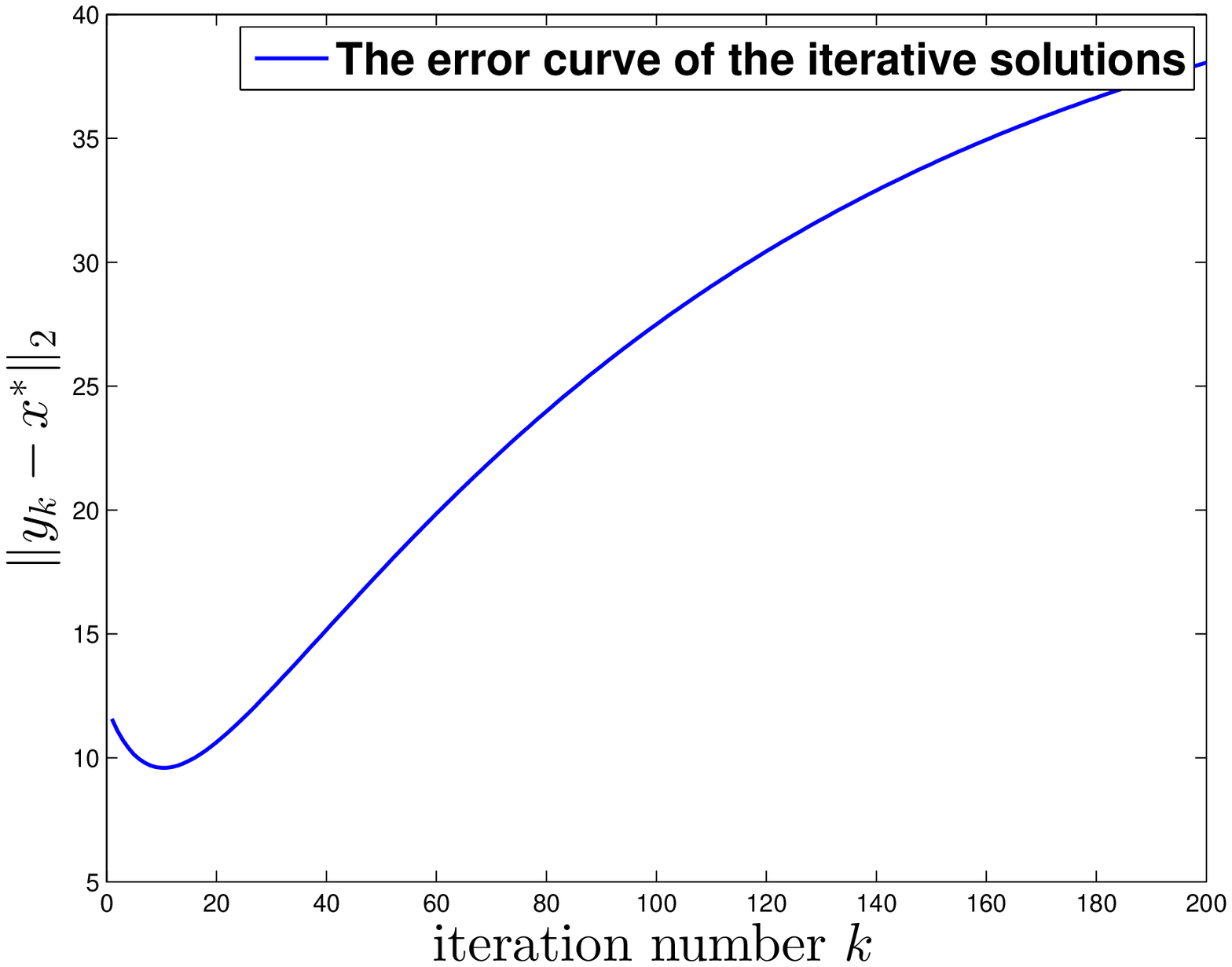}
  \end{minipage}}
  \caption{\footnotesize The iterative error curves of the Kaczmarz-Tanabe method for Model Problem 2}\label{error_model2}
\end{figure}

\begin{figure}[hbt]
  \centering
  \subfigure[$K=200,\delta=0$]{
  \begin{minipage}[ht]{.32\linewidth}
      \includegraphics[width=\linewidth]{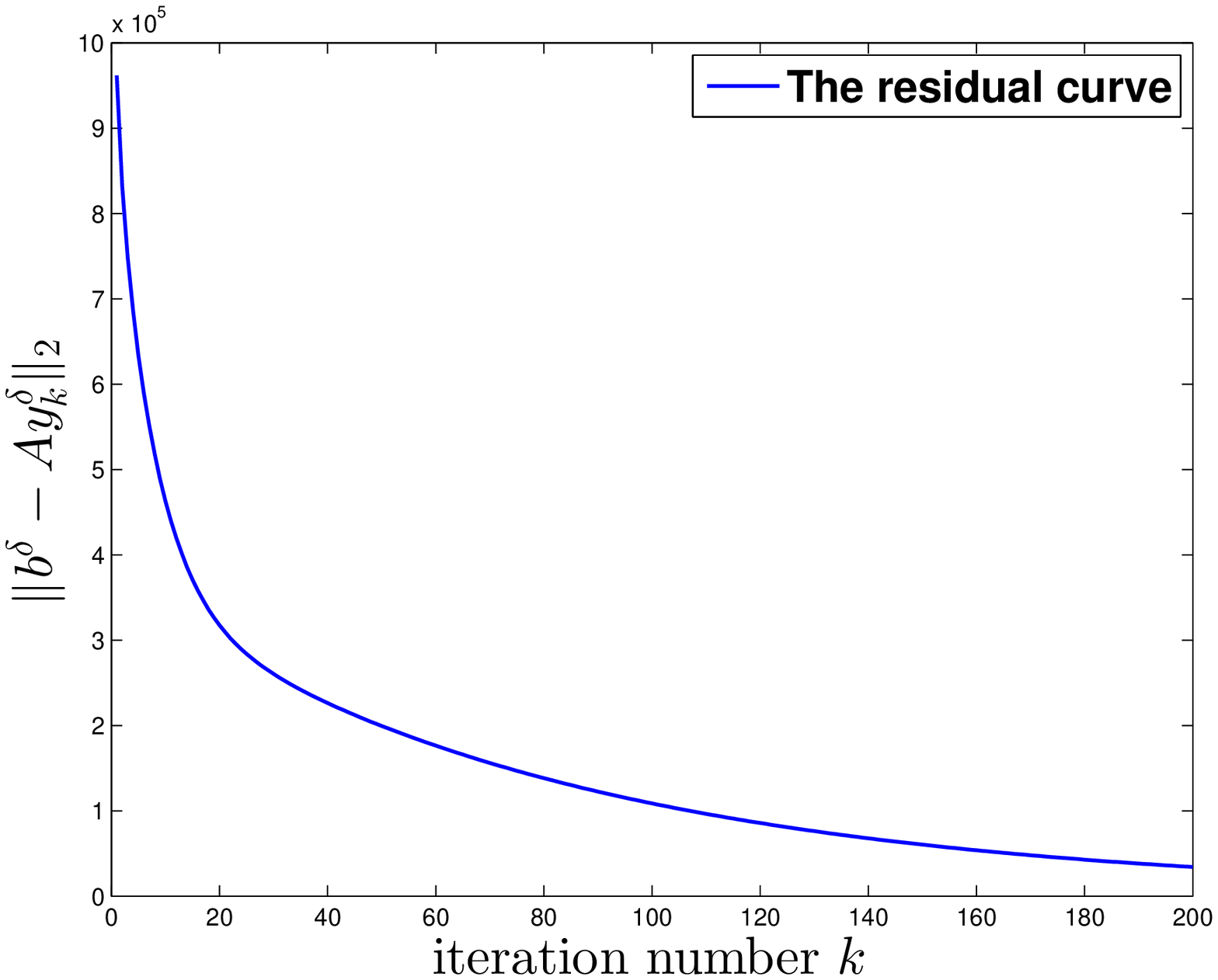}
  \end{minipage}}
  \subfigure[$K=200,\delta=0.1$]{
  \begin{minipage}[ht]{.32\linewidth}
      \includegraphics[width=\linewidth]{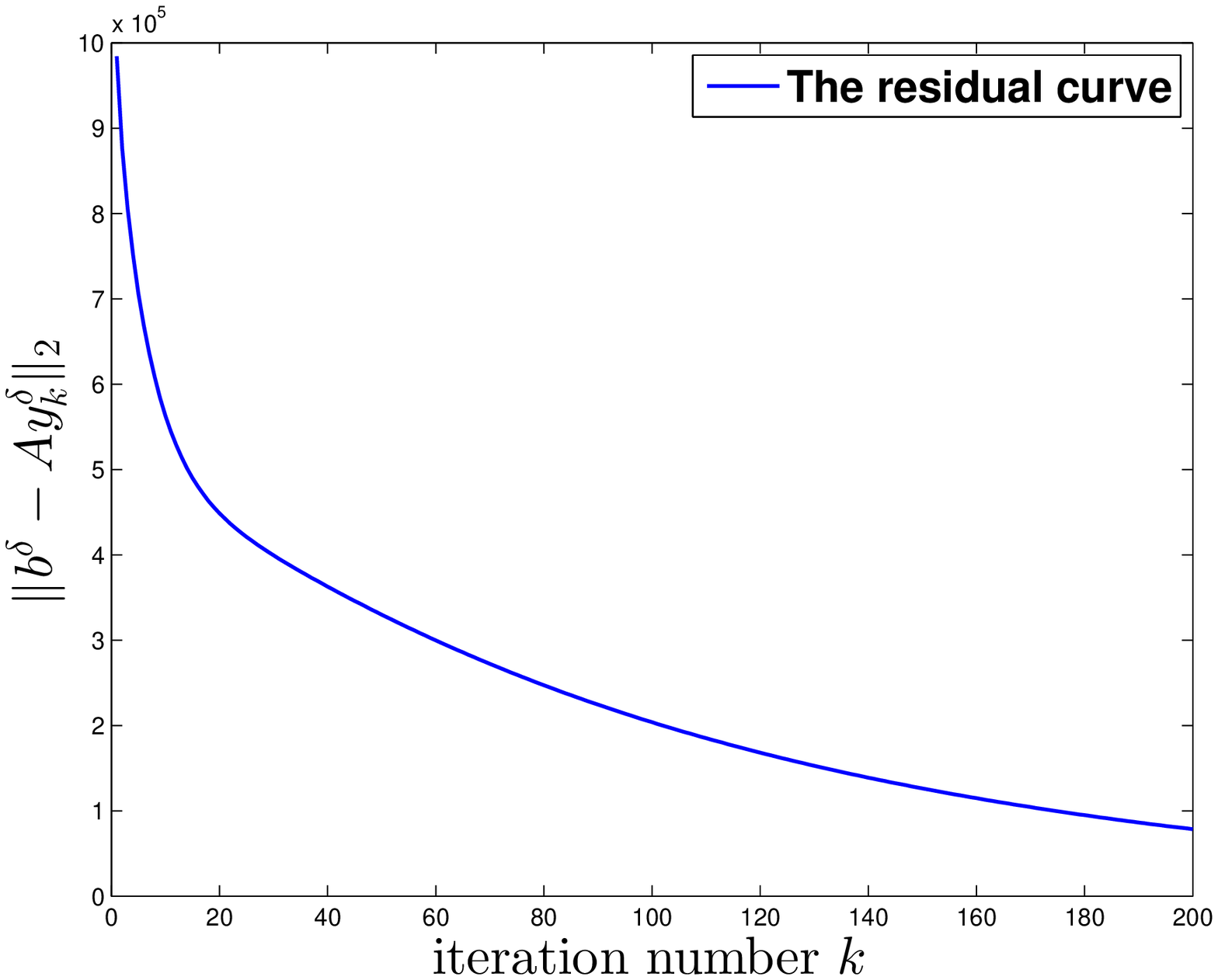}
  \end{minipage}}
  \subfigure[$K=200,\delta=0.3$]{
  \begin{minipage}[ht]{.32\linewidth}
      \includegraphics[width=\linewidth]{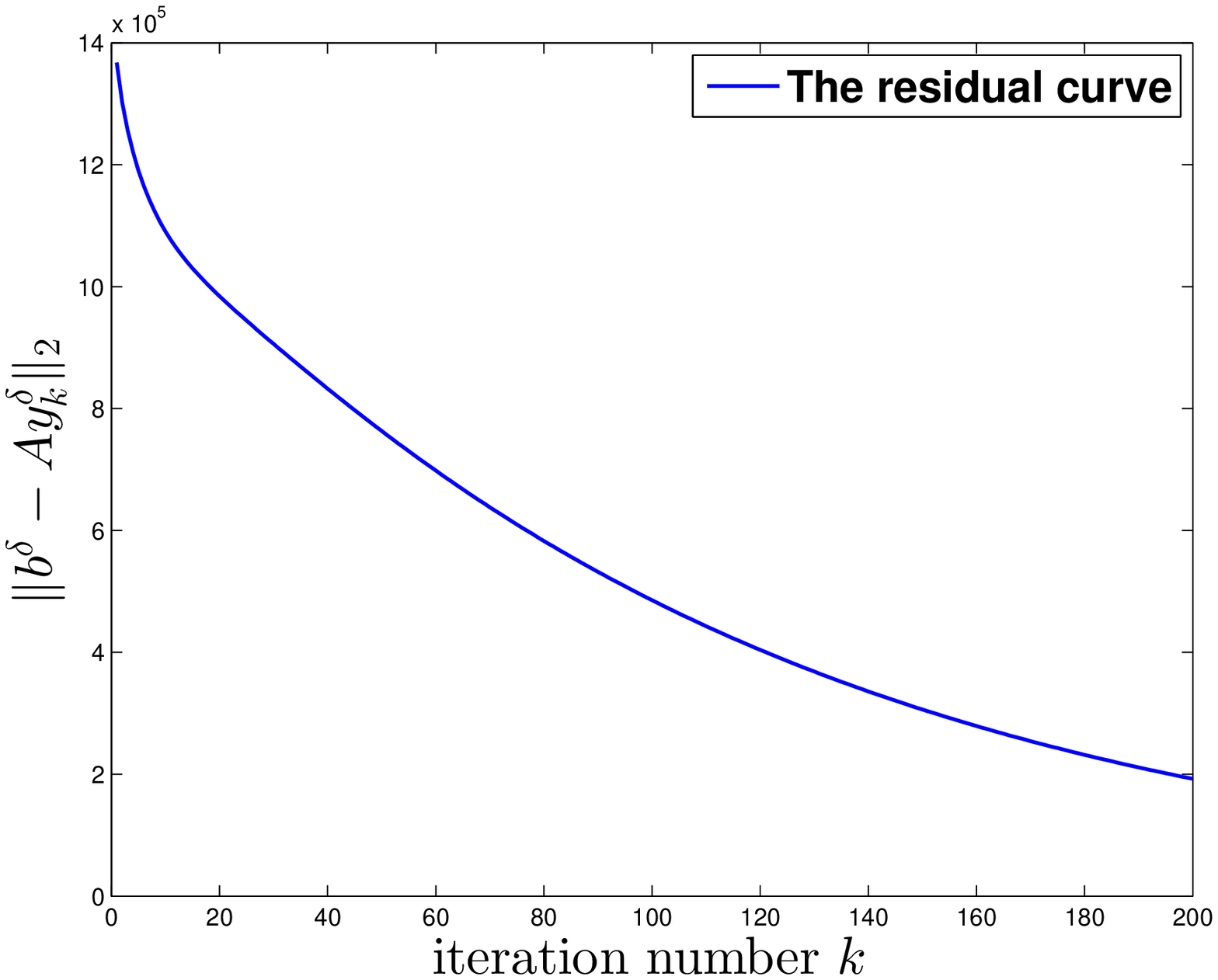}
  \end{minipage}}
  \caption{\footnotesize The residual curves of the Kaczmarz-Tanabe method for Model Problem 2}\label{residue_model2}
\end{figure}

\begin{figure}[hbt]
  \centering
  \subfigure[True image]{
  \begin{minipage}[ht]{.22\linewidth}
      \includegraphics[width=\linewidth]{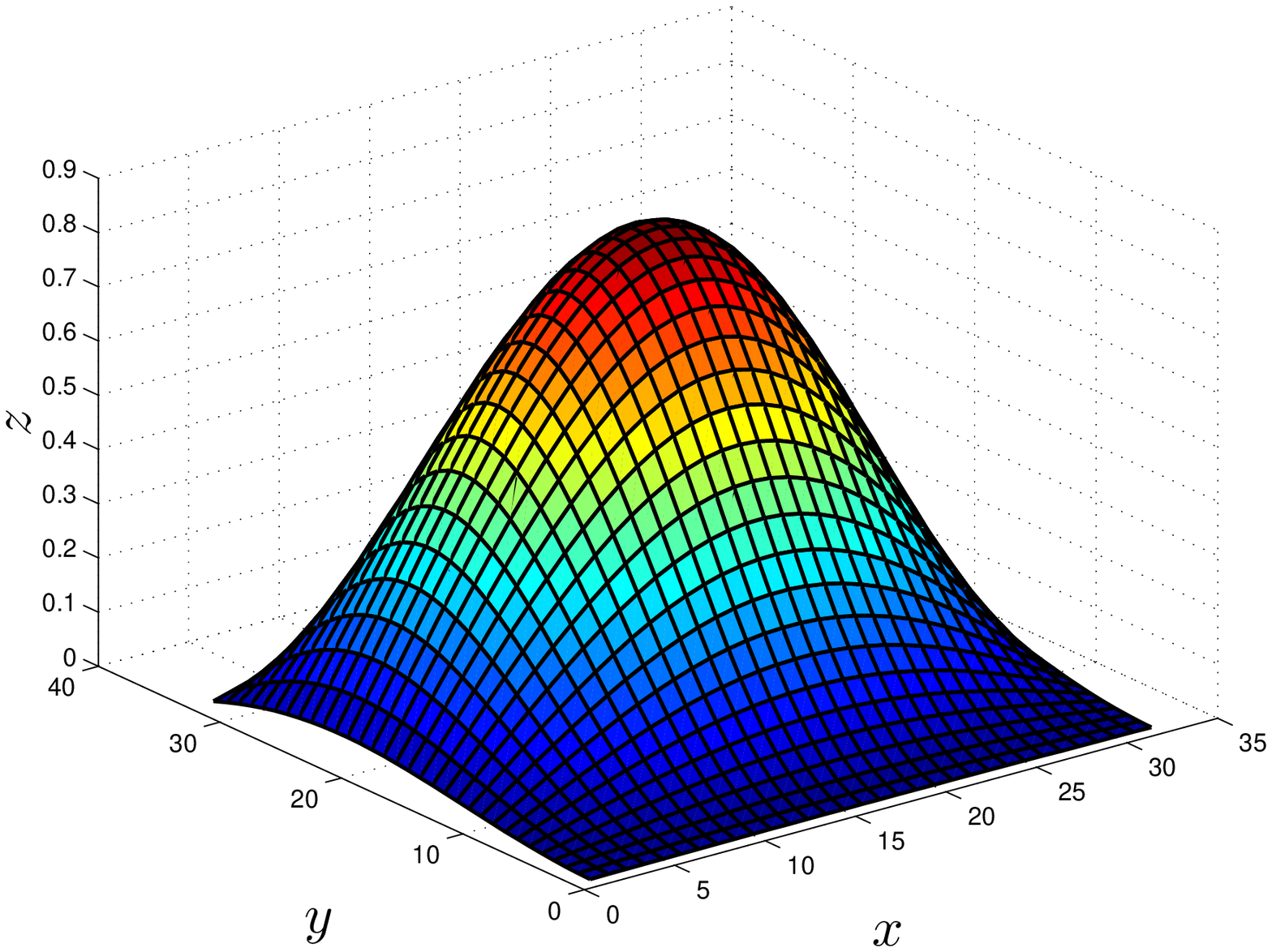}
  \end{minipage}}
  \subfigure[$K=200,\delta=0$]{
  \begin{minipage}[hbt]{.22\linewidth}
      \includegraphics[width=\linewidth]{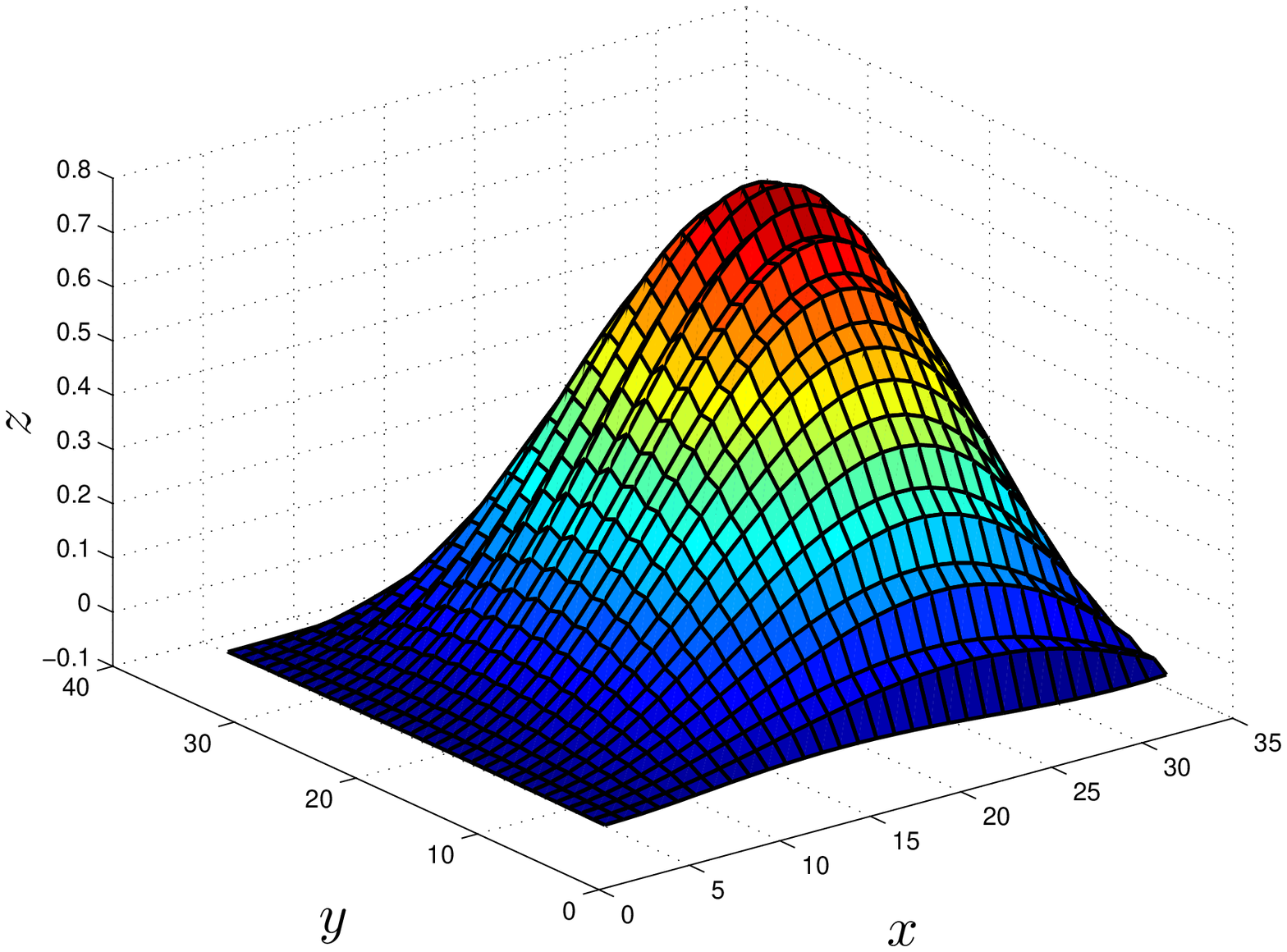}
  \end{minipage}}
   \subfigure[$K=200,\delta=0.1$]{
  \begin{minipage}[hbt]{.22\linewidth}
      \includegraphics[width=\linewidth]{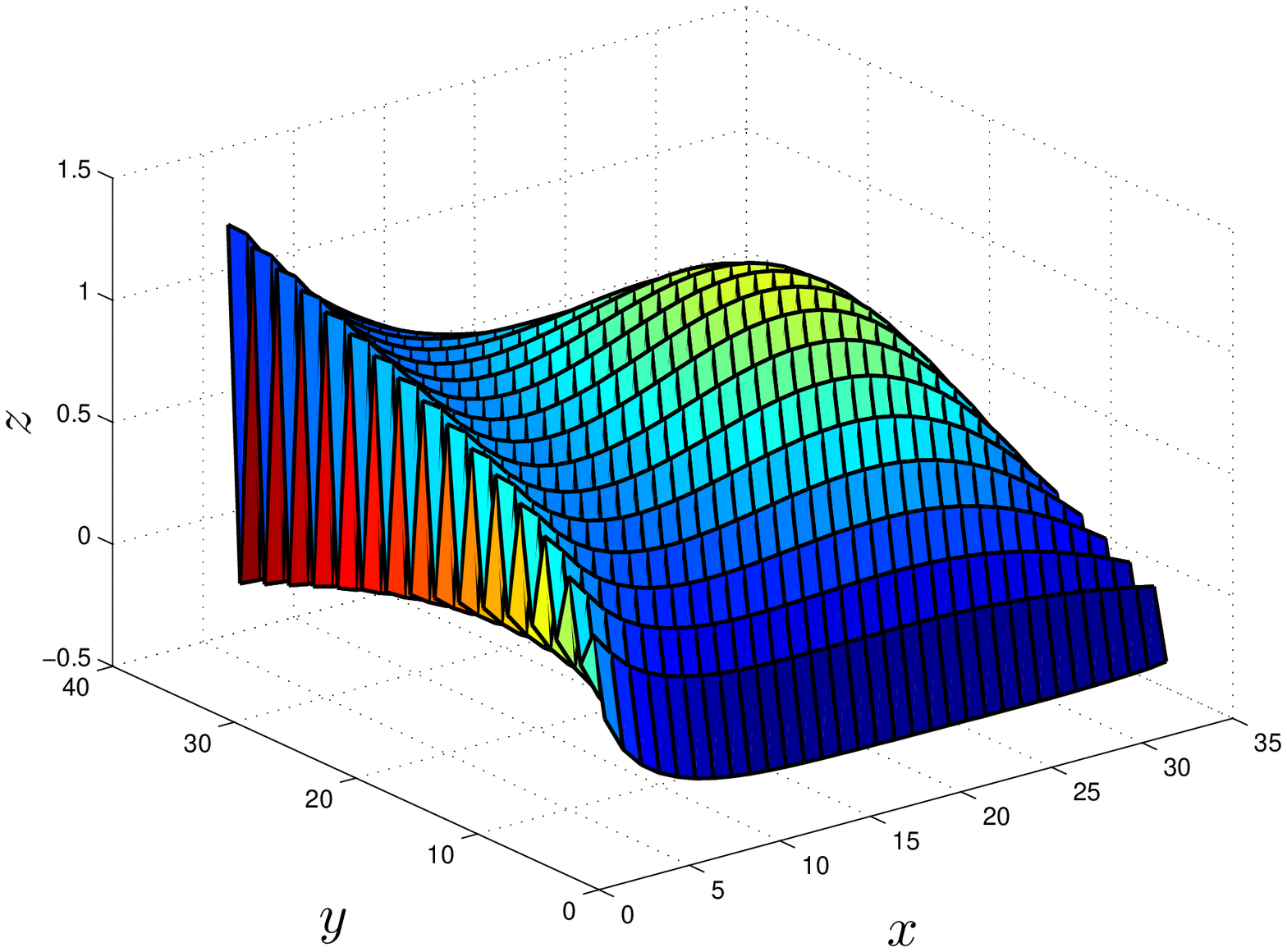}
  \end{minipage}}
     \subfigure[$K=200,\delta=0.3$]{
  \begin{minipage}[hbt]{.22\linewidth}
      \includegraphics[width=\linewidth]{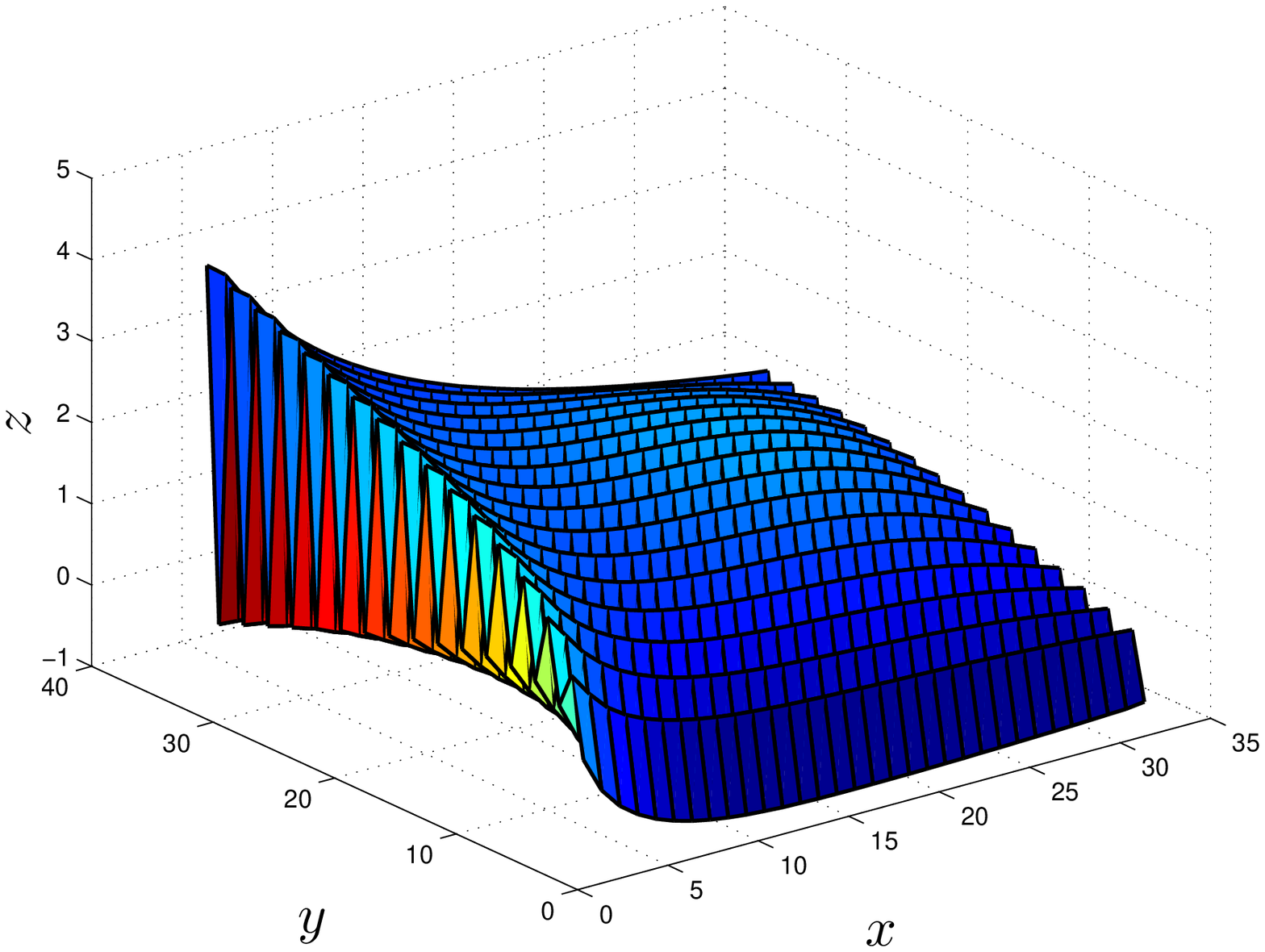}
  \end{minipage}}
  \caption{\footnotesize True image and reconstruction images of Model Problem 2}\label{image_model2}
\end{figure}

Numerical results are presented in Figures \ref{error_model2}$\sim$\ref{image_model2}. In Figures \ref{error_model2}(a) and \ref{residue_model2}(a), we see that the iterative error and residual curves are decreasing about the iteration number $k$.
Meanwhile, in Figures \ref{error_model2}(b)(c) and \ref{residue_model2}(b)(c), although the residual curves are monotonically decreasing, the iterative error curves have strong volatility. Figure \ref{image_model2} shows the comparison between the real image and the numerical image at noisy level $\delta=0.1$ and $0.3$ when the maximal iteration number $K=200$. The absolute and relative errors of the right-hand side for Model Problem 2 are listed in Table \ref{errror_level_for_model_problem2}.
\begin{table}[!htb]
  \centering
 \caption{The absolute and relative errors of the right-hand side for Model Problem 2}
 \setlength{\tabcolsep}{1.5cm}{
  \begin{tabular}{ccc}
    \hline
    \text{Menu}                      &\text{$\delta=0.1$} &\text{$\delta=0.3$}  \\
    \hline
    $\|b^\delta-b\|_2$               & 3.6611e+05         &  1.0983e+06                \\
    $\|b^\delta-b\|_2/\|b\|_2$       & 0.3778             &  1.1333                  \\
    \hline
  \end{tabular}}
\label{errror_level_for_model_problem2}
\end{table}

\subsection{Model Problem 3. Head Phantom}

\noindent In computed tomography, the distribution of some physical parameter(such as absorption intensities) at the cross-section of the object need to be reconstructed from the projection data such as medical diagnosis--the distribution of the absorption intensities of tissue slice need to be reconstructed from x-ray data. The computed tomography system attributes to a linear system $Ax=b$, where $A$ is a projected system of computed tomography, $b$ is scanning data, $x$ is unknown intensity image of an object. In the general case, the system is overdetermined.

The linear system is generated from the subroutine 'parallel' in ARTool package\upcite{Hansen2018AIRtools}, and there are 36 projective angles at equal intervals in $[0,2\pi]$ and 75 equi-spaced parallel rays per angle. The phantom is discretized into $50\times 50$ pixels. and the dimension of $A$ is $2700\times 2500$. We take the right-hand side
\begin{align}\label{perturbed.right.for.headphantom}
  b_i^\delta=b_i+\delta\max\limits_{i}|b_i|, \qquad i=1,\ldots,m.
\end{align}
in \eqref{linear.system.perturbation} with $\delta=0, 0.01, 0.02$ and $0.05$, respectively. The corresponding absolute and relative errors are listed in Table \ref{errror_level_1_for_model_problem3}.
\begin{table}[!htb]
  \centering
 \caption{The absolute and relative errors of the right-hand side $b^\delta$ in \eqref{perturbed.right.for.headphantom} for Model Problem 3}
  \setlength{\tabcolsep}{0.7cm}{
  \begin{tabular}{ccccc}
    \hline
    \text{Menu}                      & \text{$\delta=0$} & \text{$\delta=0.01$} & \text{$\delta=0.02$} & \text{$\delta=0.05$} \\
    \hline
    $\|b^\delta-b\|_2$               &       0          &  6.9109               &  13.8218             &   34.5544       \\
    $\|b^\delta-b\|_2/\|b\|_2$       &       0          &  0.023                 &  0.046               &   0.115           \\
    \hline
  \end{tabular}}
\label{errror_level_1_for_model_problem3}
\end{table}

\begin{figure}[!hbt]
  \centering
  \subfigure[$\delta=0$]{
  \begin{minipage}[ht]{.22\linewidth}
      \includegraphics[width=\linewidth]{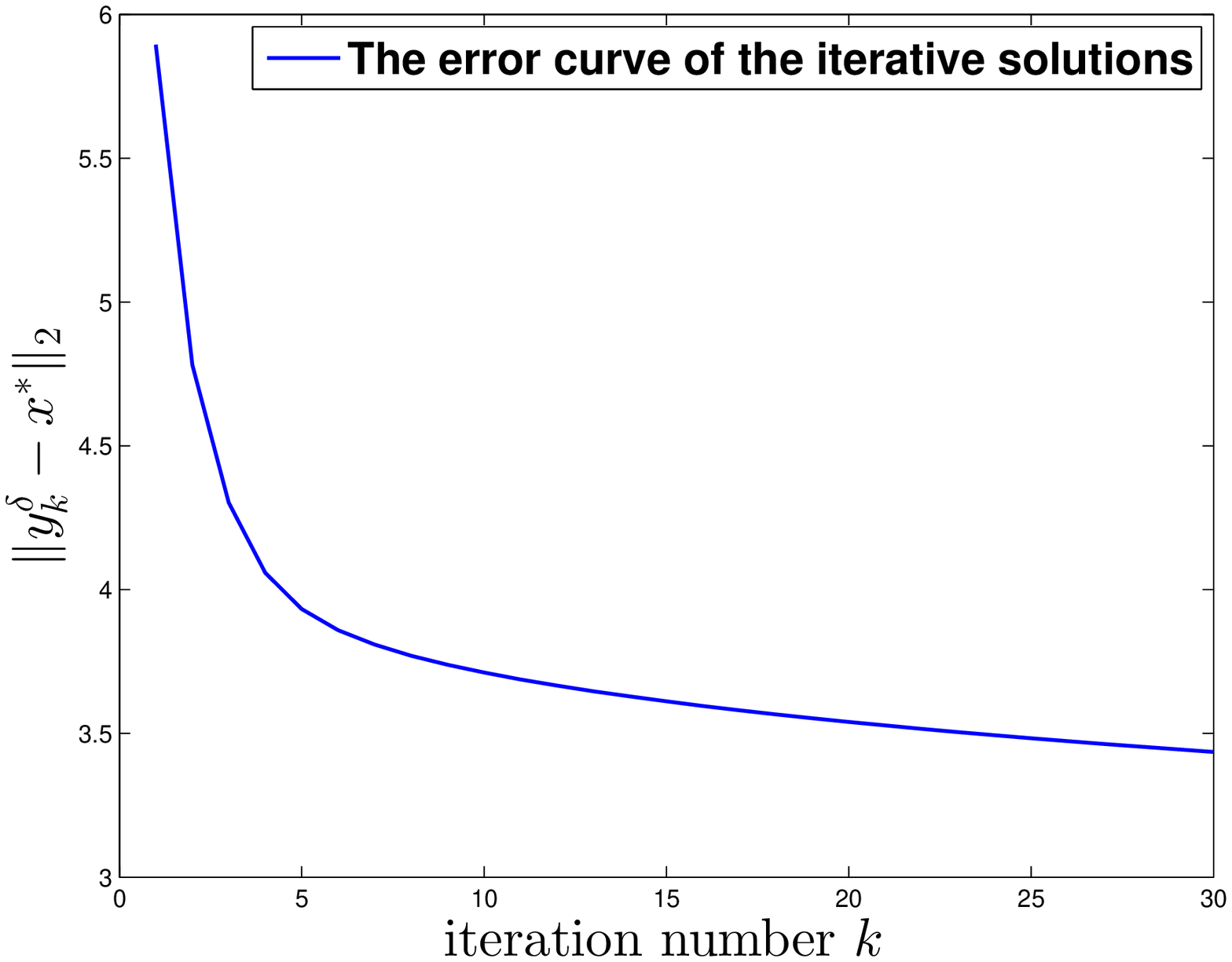}
  \end{minipage}}
  \subfigure[$\delta=0.01$]{
  \begin{minipage}[ht]{.22\linewidth}
      \includegraphics[width=\linewidth]{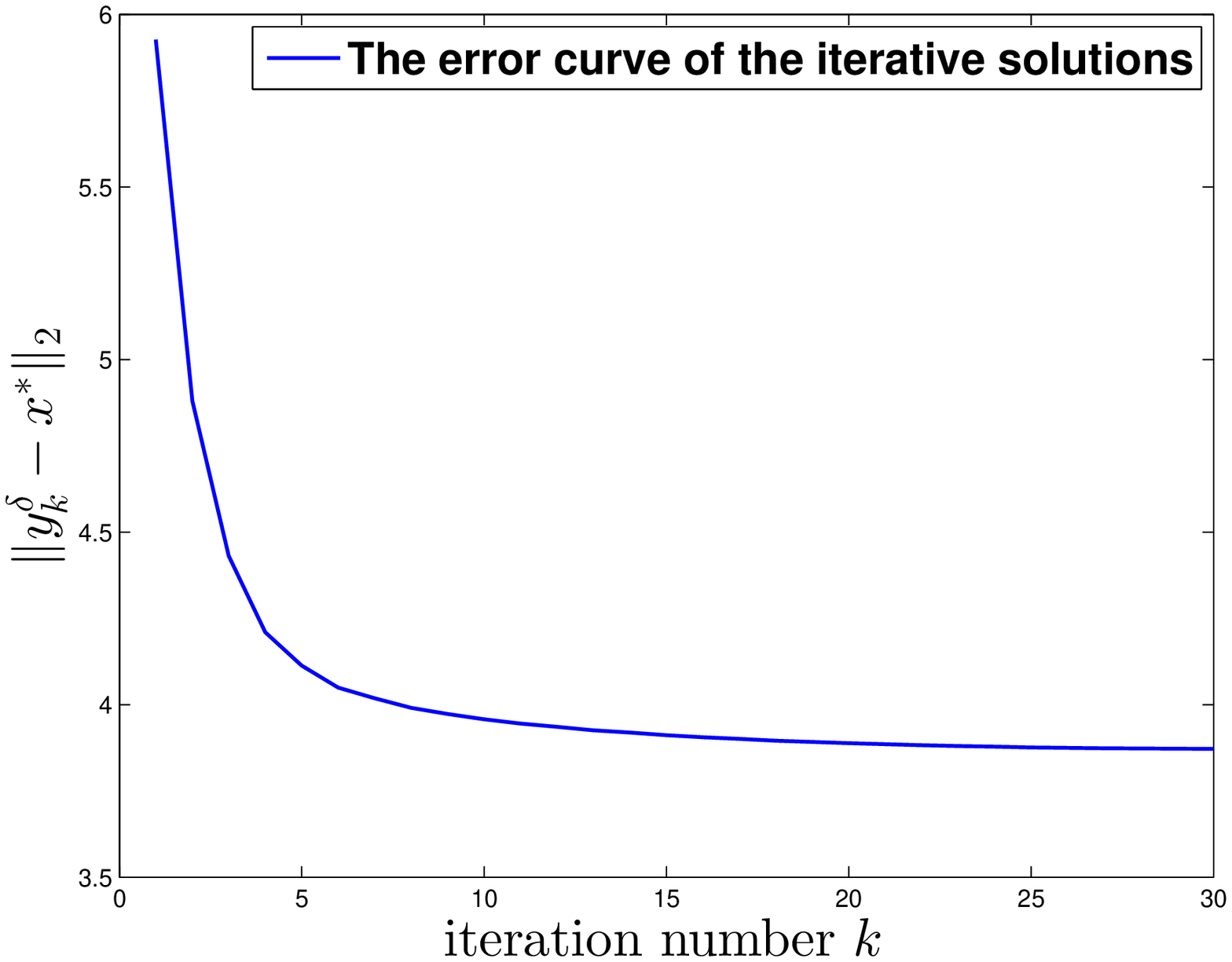}
  \end{minipage}}
  \subfigure[$\delta=0.02$]{
  \begin{minipage}[ht]{.22\linewidth}
      \includegraphics[width=\linewidth]{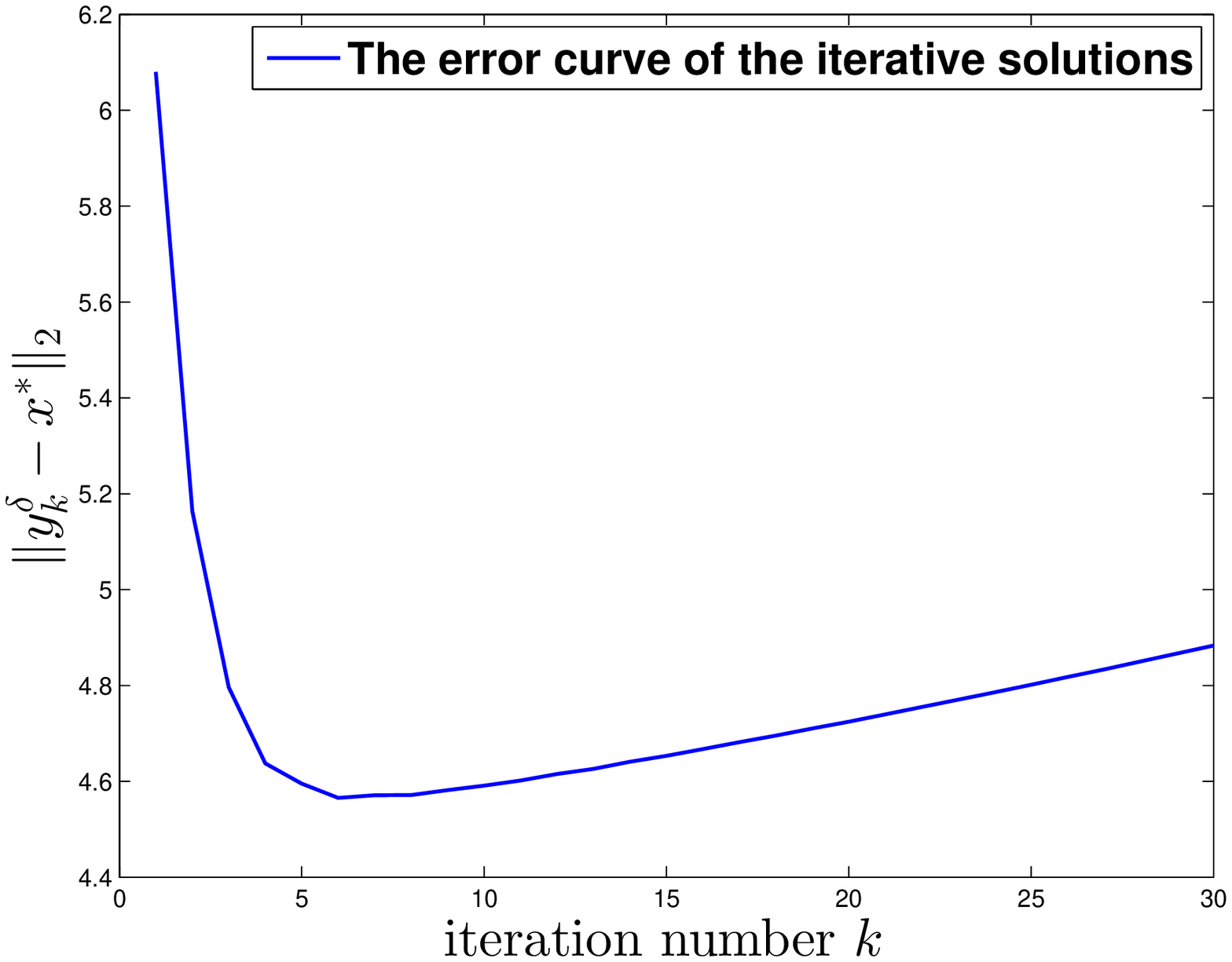}
  \end{minipage}}
  \subfigure[$\delta=0.05$]{
  \begin{minipage}[ht]{.22\linewidth}
      \includegraphics[width=\linewidth]{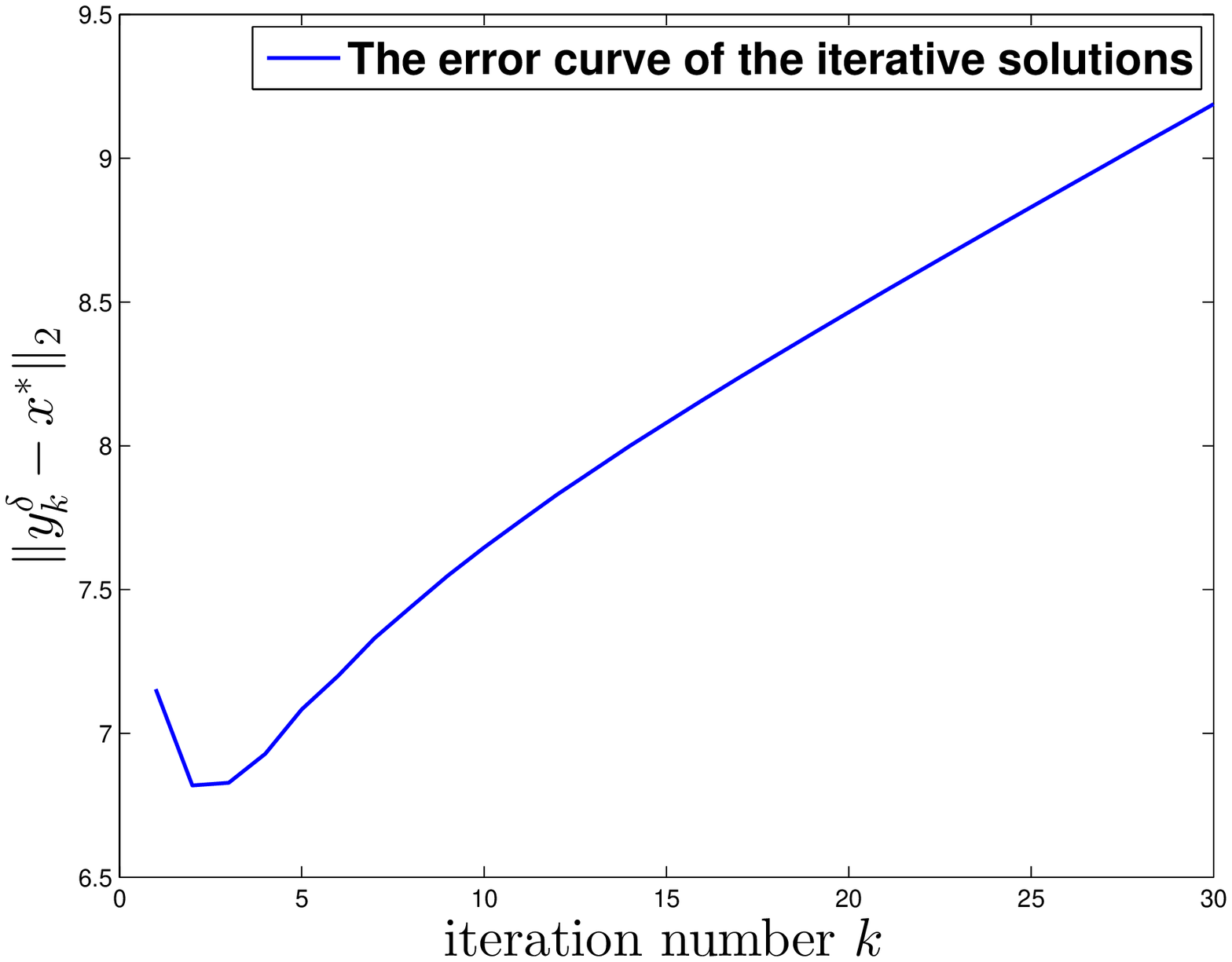}
  \end{minipage}}
  \caption{\footnotesize The error curves of the Kaczmarz-Tanabe method for Model Problem 3}\label{the_periodic_Kaczmarz_method_errors}
\end{figure}

\begin{figure}[!hbt]
  \centering
  \subfigure[$\delta=0$]{
  \begin{minipage}[ht]{.22\linewidth}
      \includegraphics[width=\linewidth]{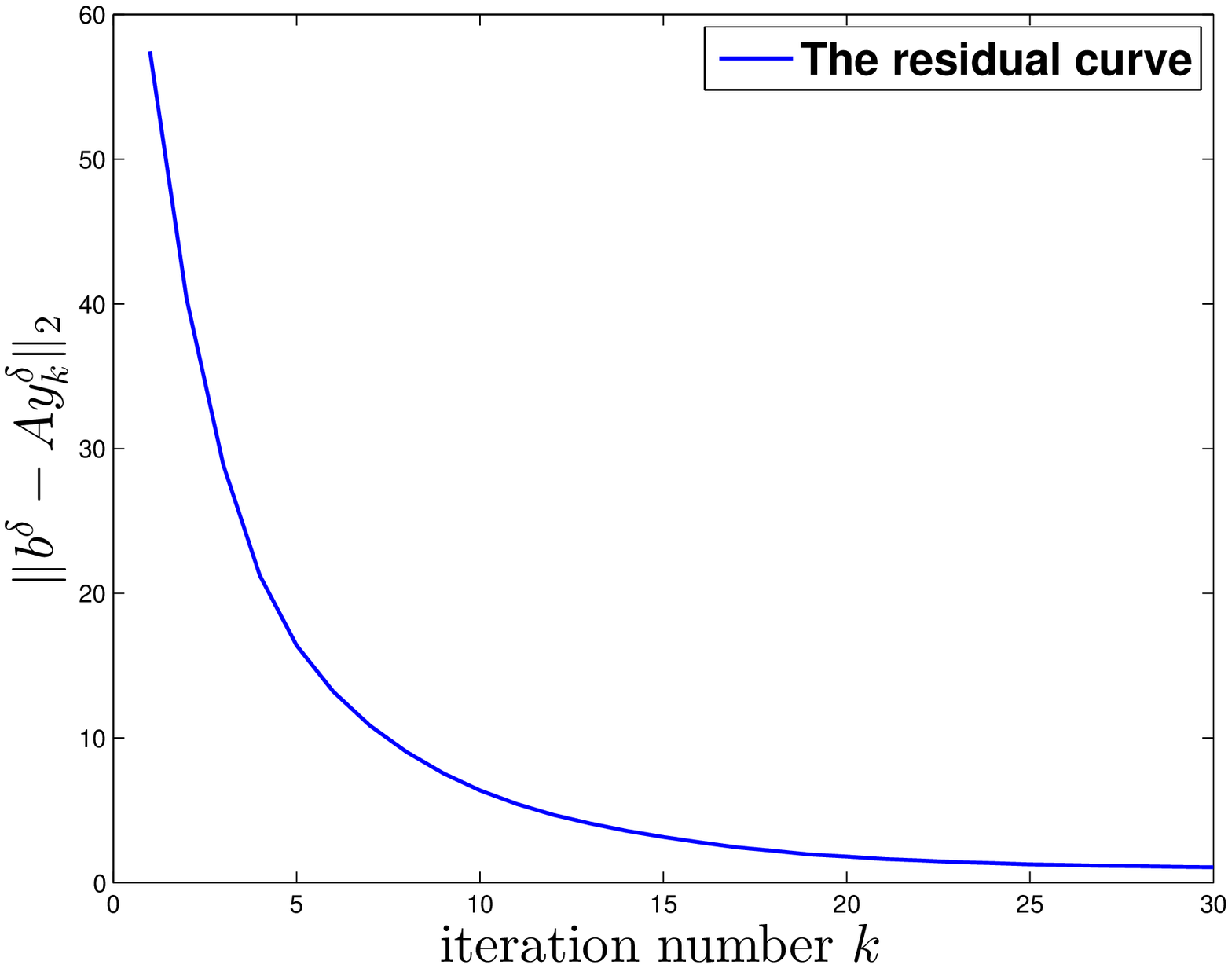}
  \end{minipage}}
  \subfigure[$\delta=0.01$]{
  \begin{minipage}[ht]{.22\linewidth}
      \includegraphics[width=\linewidth]{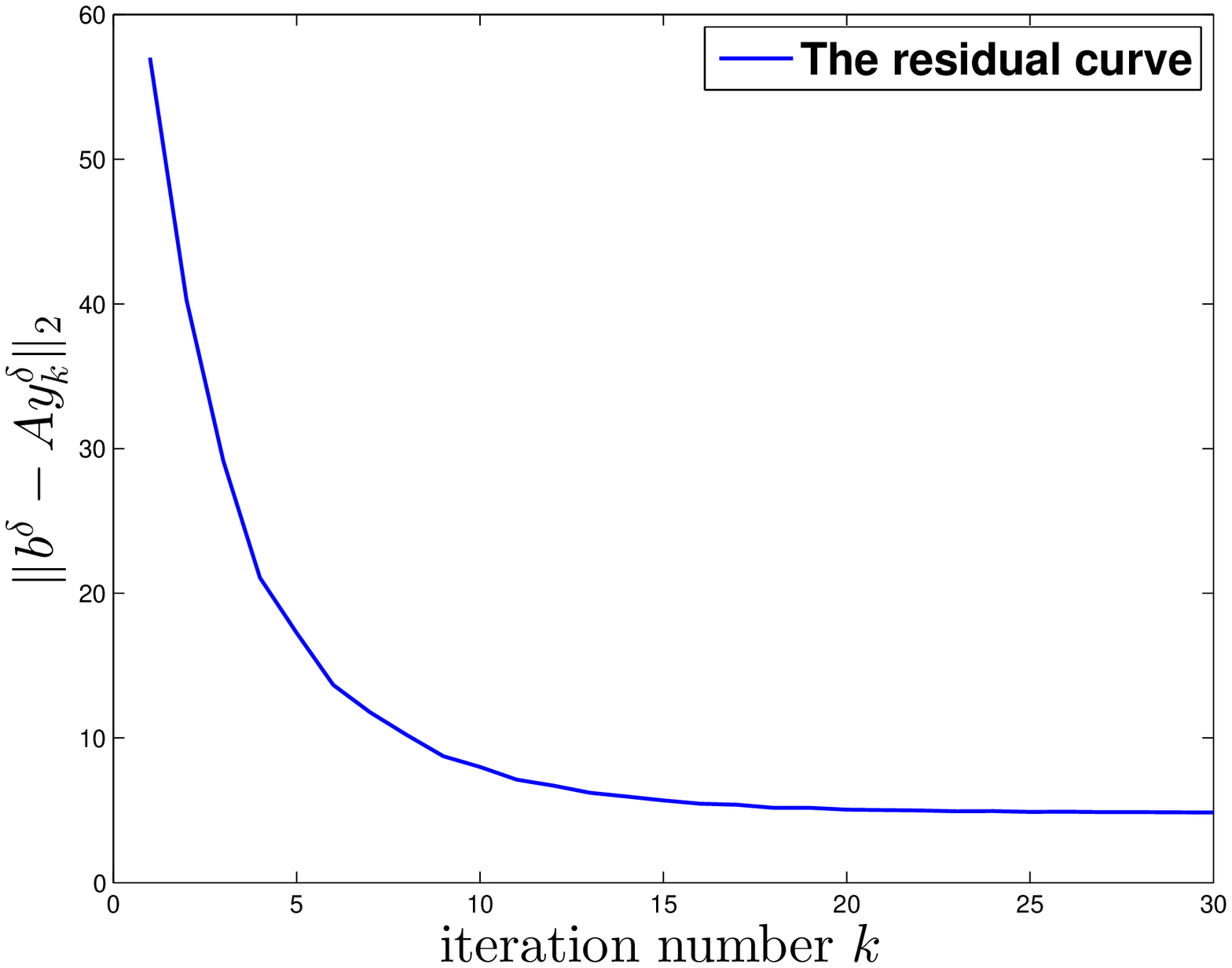}
  \end{minipage}}
  \subfigure[$\delta=0.02$]{
  \begin{minipage}[ht]{.22\linewidth}
      \includegraphics[width=\linewidth]{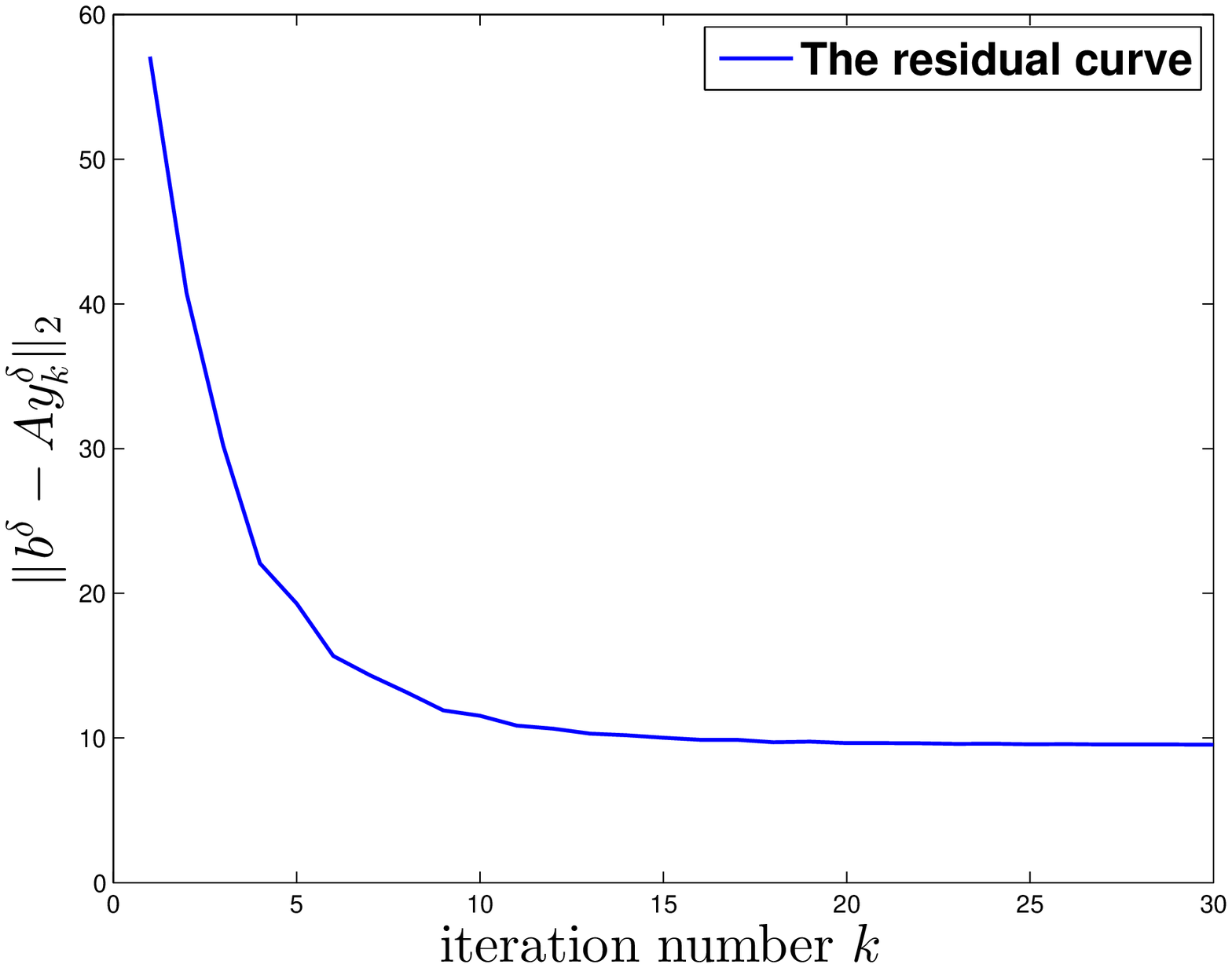}
  \end{minipage}}
  \subfigure[$\delta=0.05$]{
  \begin{minipage}[ht]{.22\linewidth}
      \includegraphics[width=\linewidth]{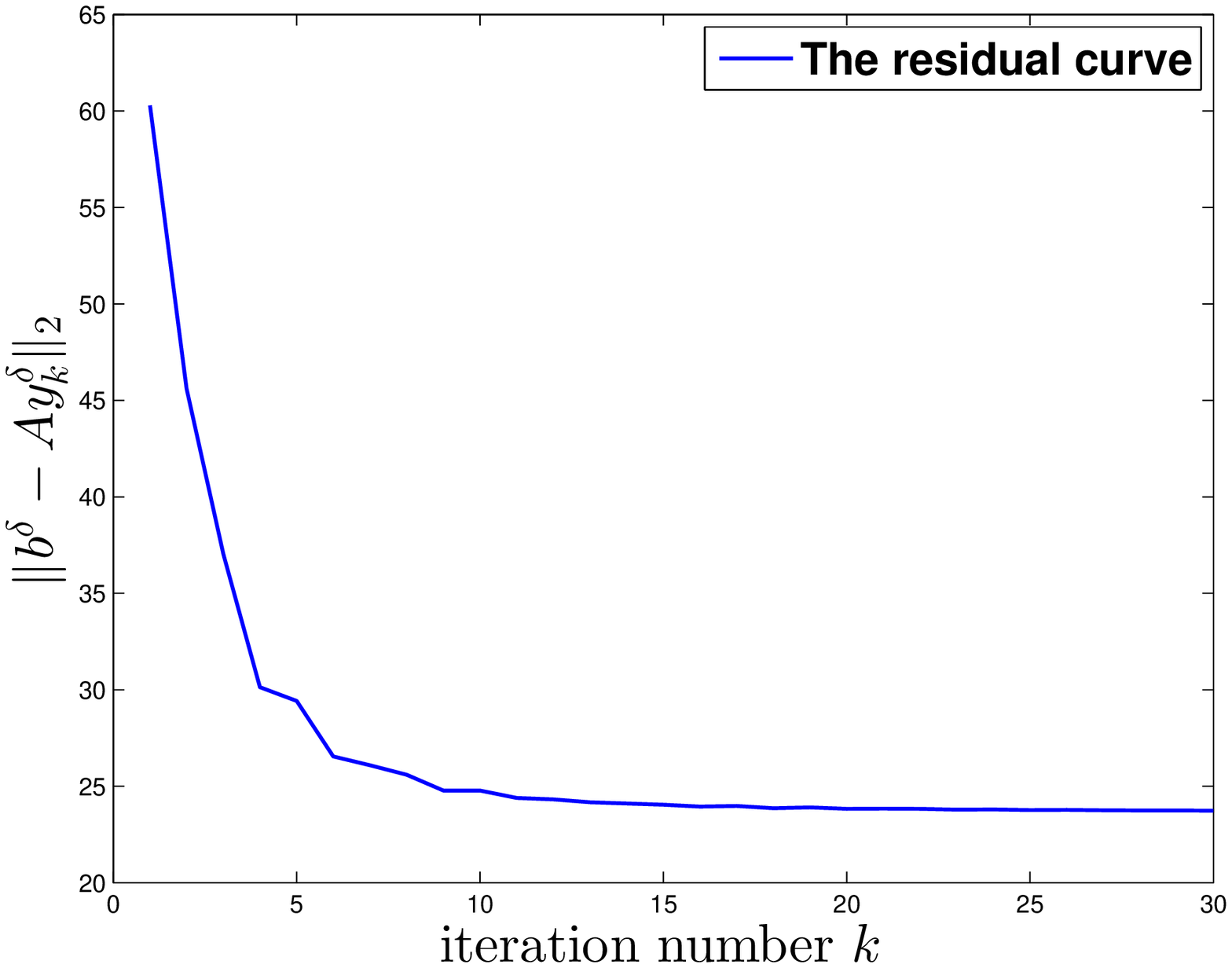}
  \end{minipage}}
  \caption{\footnotesize The residual curves of the Kaczmarz-Tanabe method for Model Problem 3}\label{the_periodic_Kaczmarz_method_residues}
\end{figure}

\begin{figure}[!hbt]
  \centering
  \subfigure[$\delta=0$]{
  \begin{minipage}[ht]{.22\linewidth}
      \includegraphics[width=\linewidth]{figures/Phantom/Kaczmarz.Tanabe/phantom_numerical_image_0.eps}
  \end{minipage}}
   \subfigure[$\delta=0.01$]{
  \begin{minipage}[ht]{.22\linewidth}
      \includegraphics[width=\linewidth]{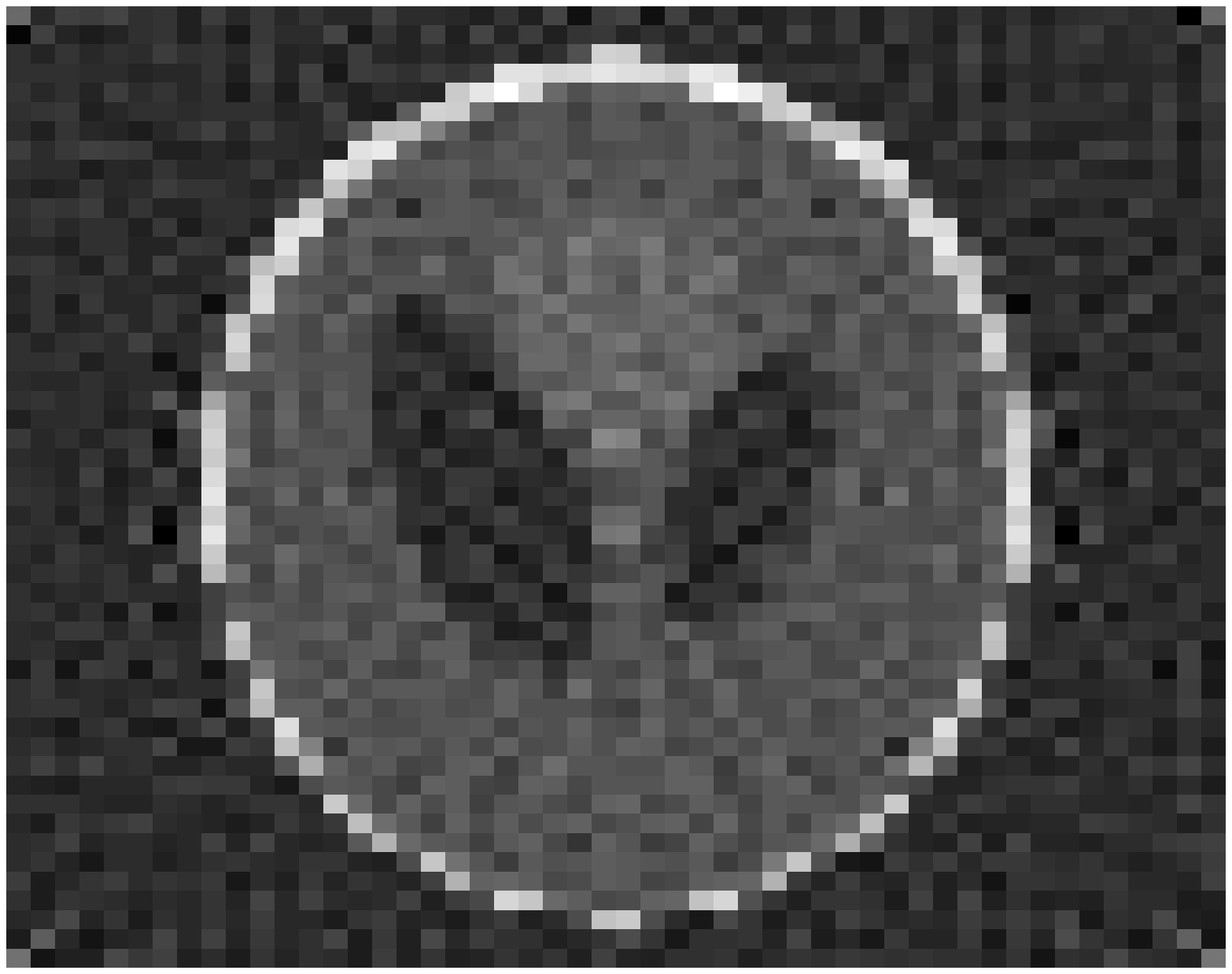}
  \end{minipage}}
  \subfigure[$\delta=0.02$]{
  \begin{minipage}[ht]{.22\linewidth}
      \includegraphics[width=\linewidth]{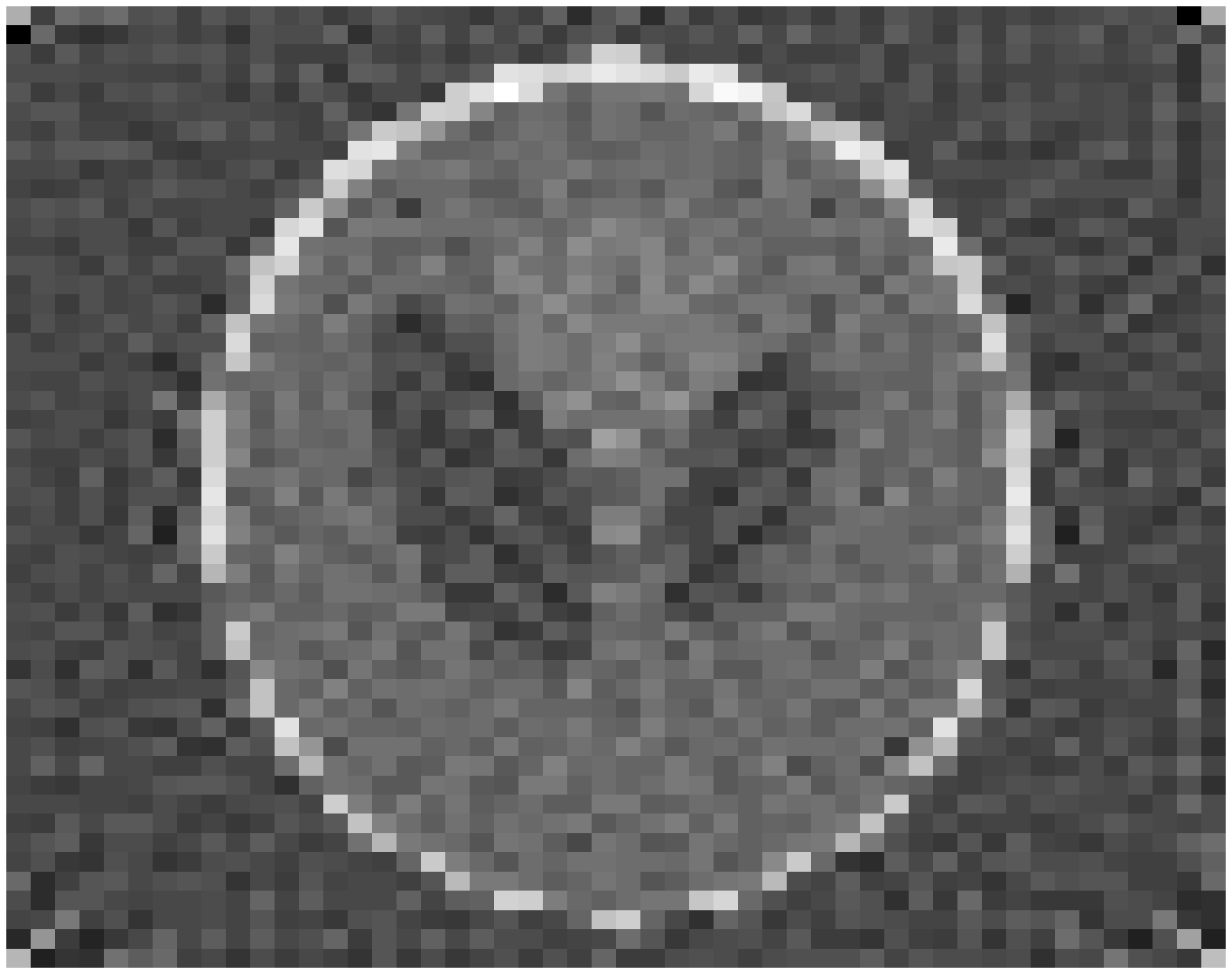}
  \end{minipage}}
  \subfigure[$\delta=0.05$]{
  \begin{minipage}[ht]{.22\linewidth}
      \includegraphics[width=\linewidth]{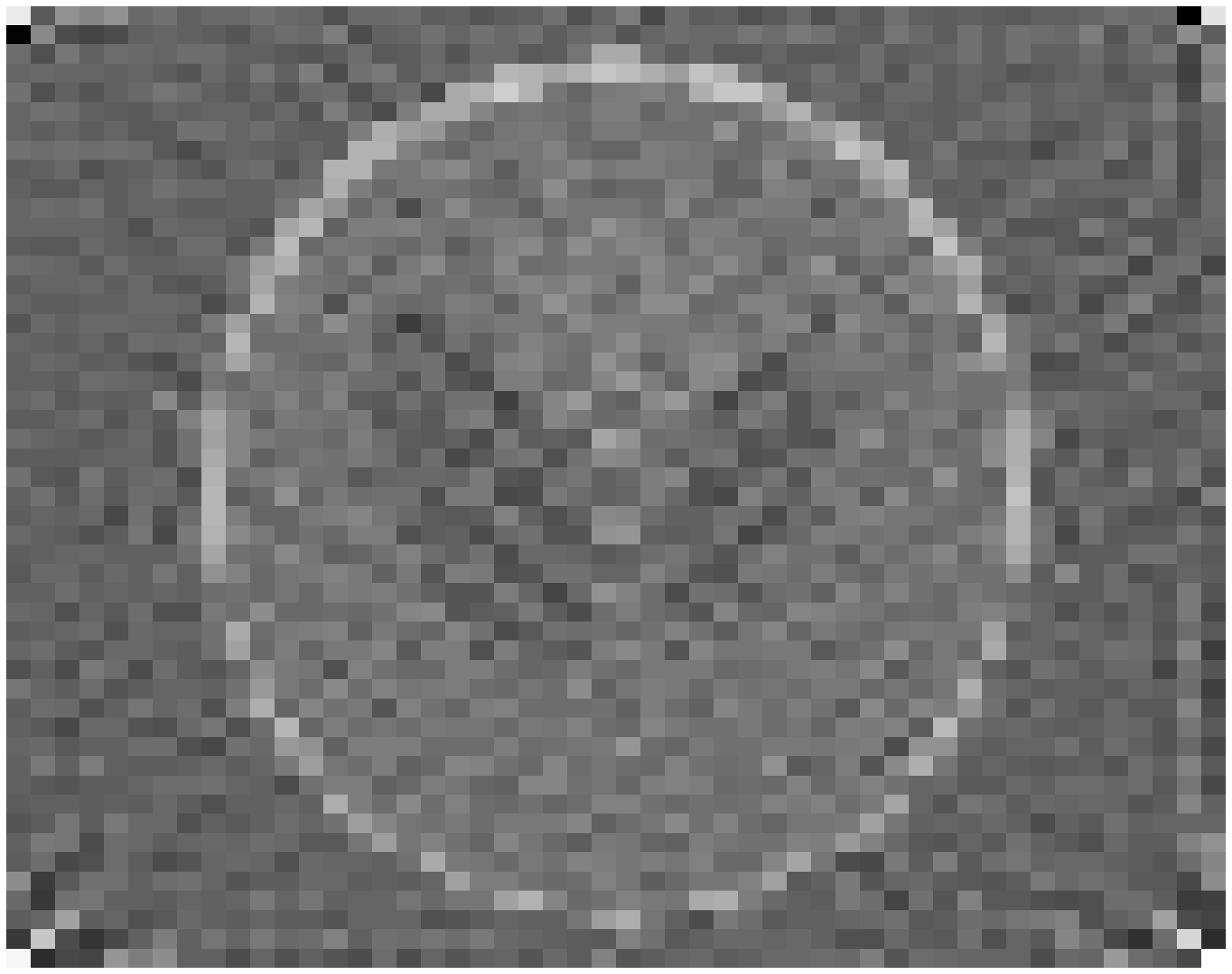}
  \end{minipage}}
  \caption{\footnotesize Numerical images of the Kaczmarz-Tanabe method for Model Problem 3}\label{the_Kaczmarz_Tanabe_method_images}
\end{figure}
\begin{figure}[!hbt]
  \centering
  \subfigure[$\delta=0$]{
  \begin{minipage}[ht]{.22\linewidth}
      \includegraphics[width=\linewidth]{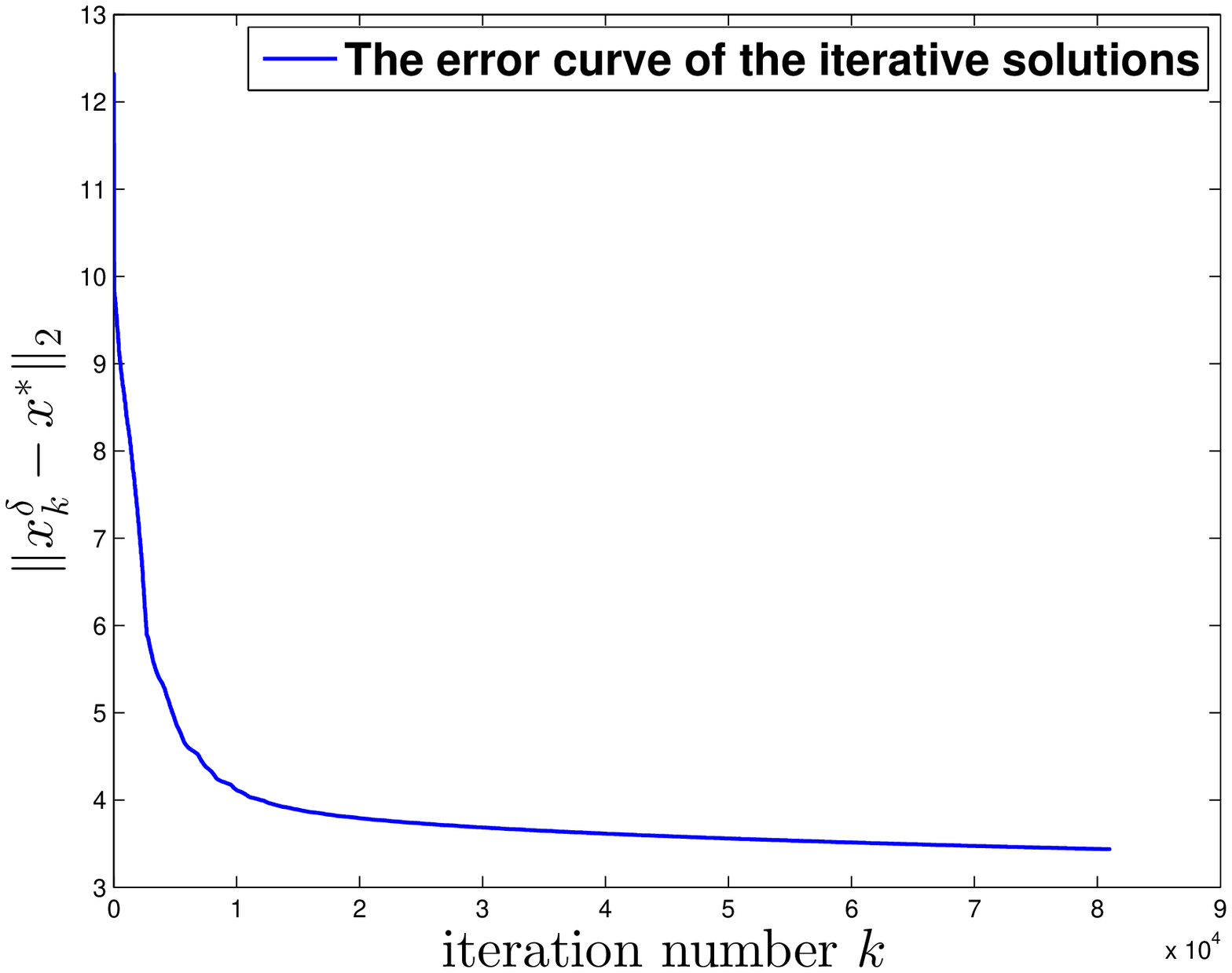}
  \end{minipage}}
  \subfigure[$\delta=0.01$]{
  \begin{minipage}[ht]{.22\linewidth}
      \includegraphics[width=\linewidth]{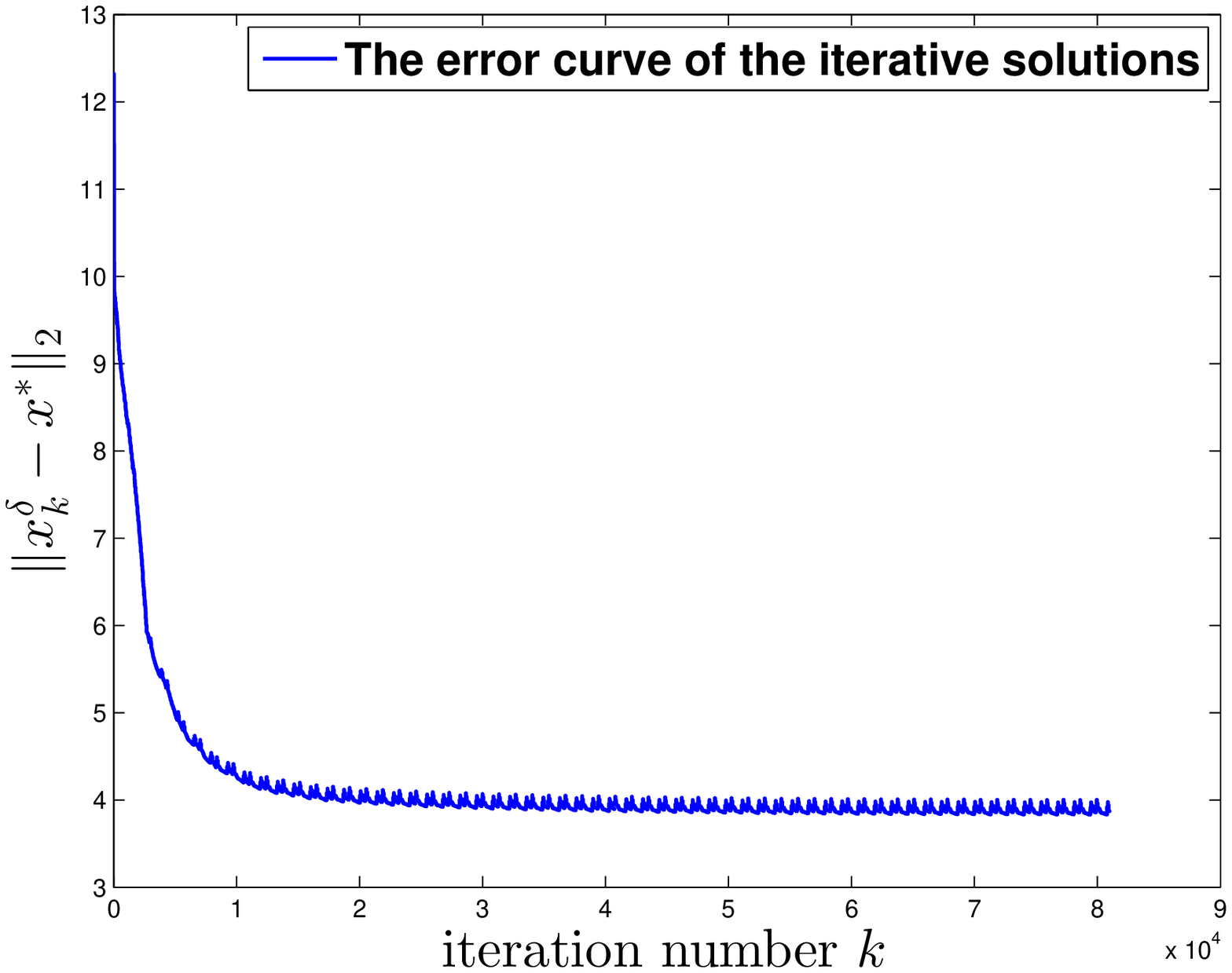}
  \end{minipage}}
    \subfigure[$\delta=0.05$]{
  \begin{minipage}[ht]{.22\linewidth}
      \includegraphics[width=\linewidth]{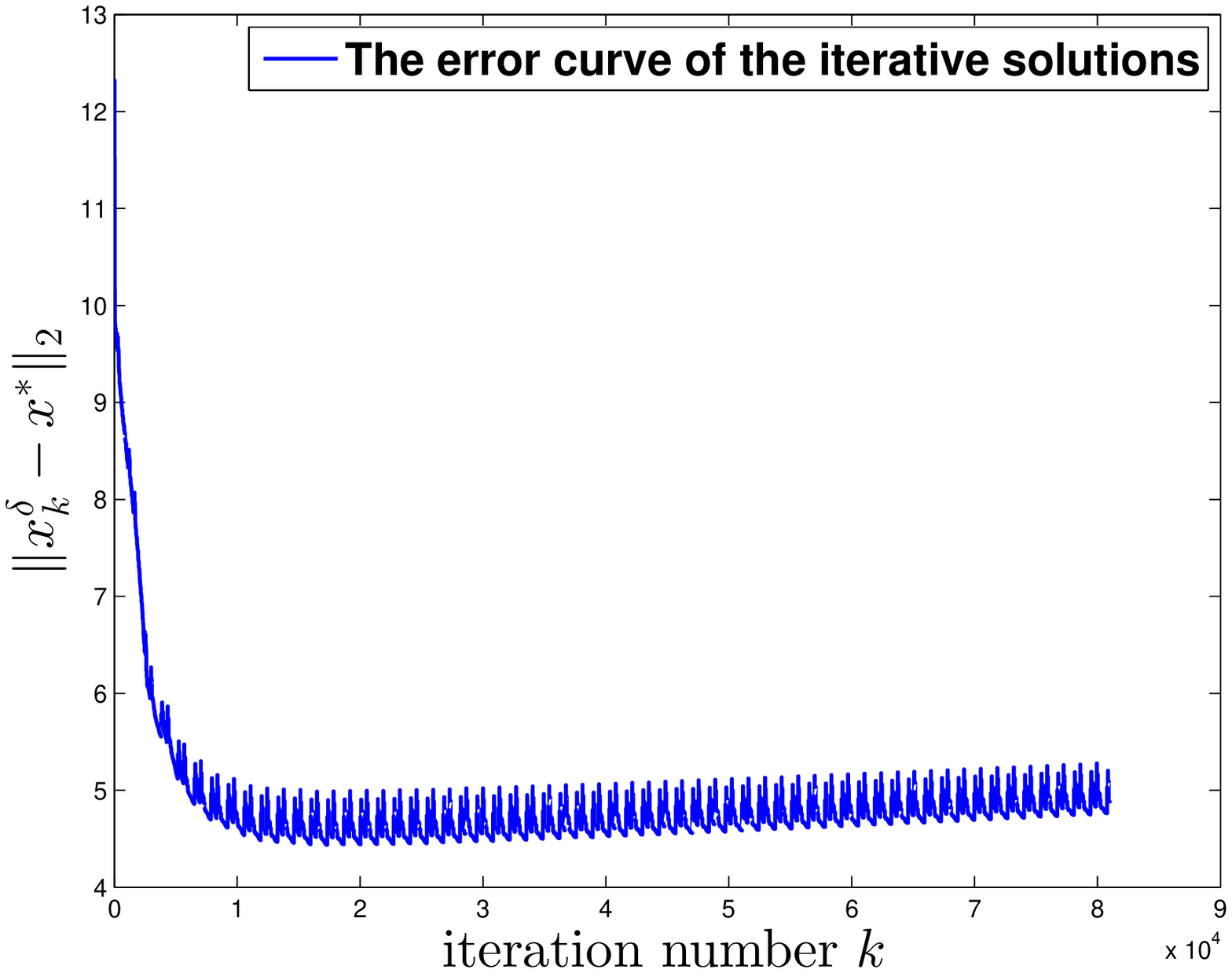}
  \end{minipage}}
  \subfigure[$\delta=0.1$]{
  \begin{minipage}[ht]{.22\linewidth}
      \includegraphics[width=\linewidth]{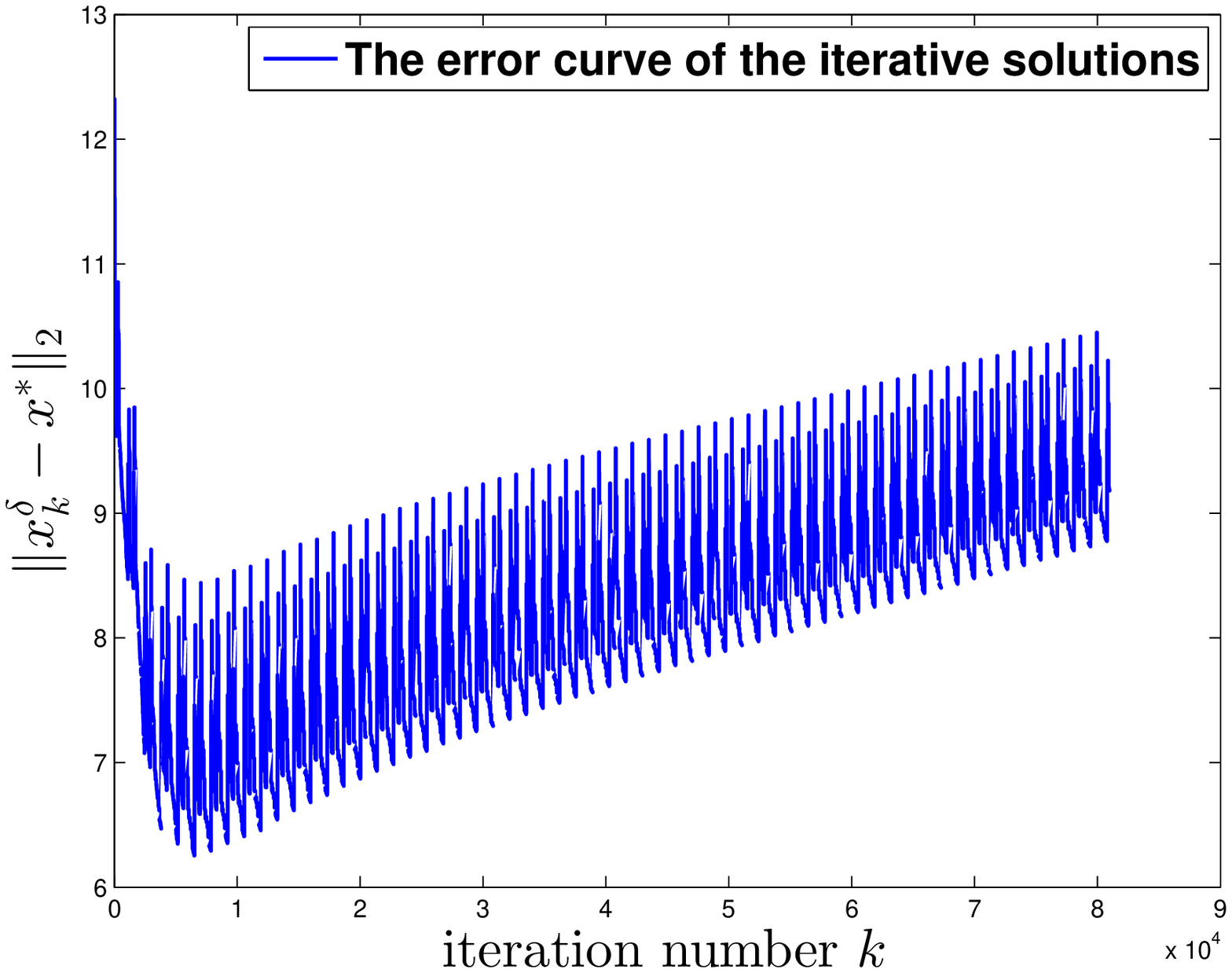}
  \end{minipage}}
  \caption{\footnotesize The error curves of the Kaczmarz method for Model Problem 3}\label{the_Kaczmarz_method_error}
\end{figure}
\begin{figure}[!hbt]
  \centering
 \subfigure[$\delta=0$]{
  \begin{minipage}[ht]{.22\linewidth}
      \includegraphics[width=\linewidth]{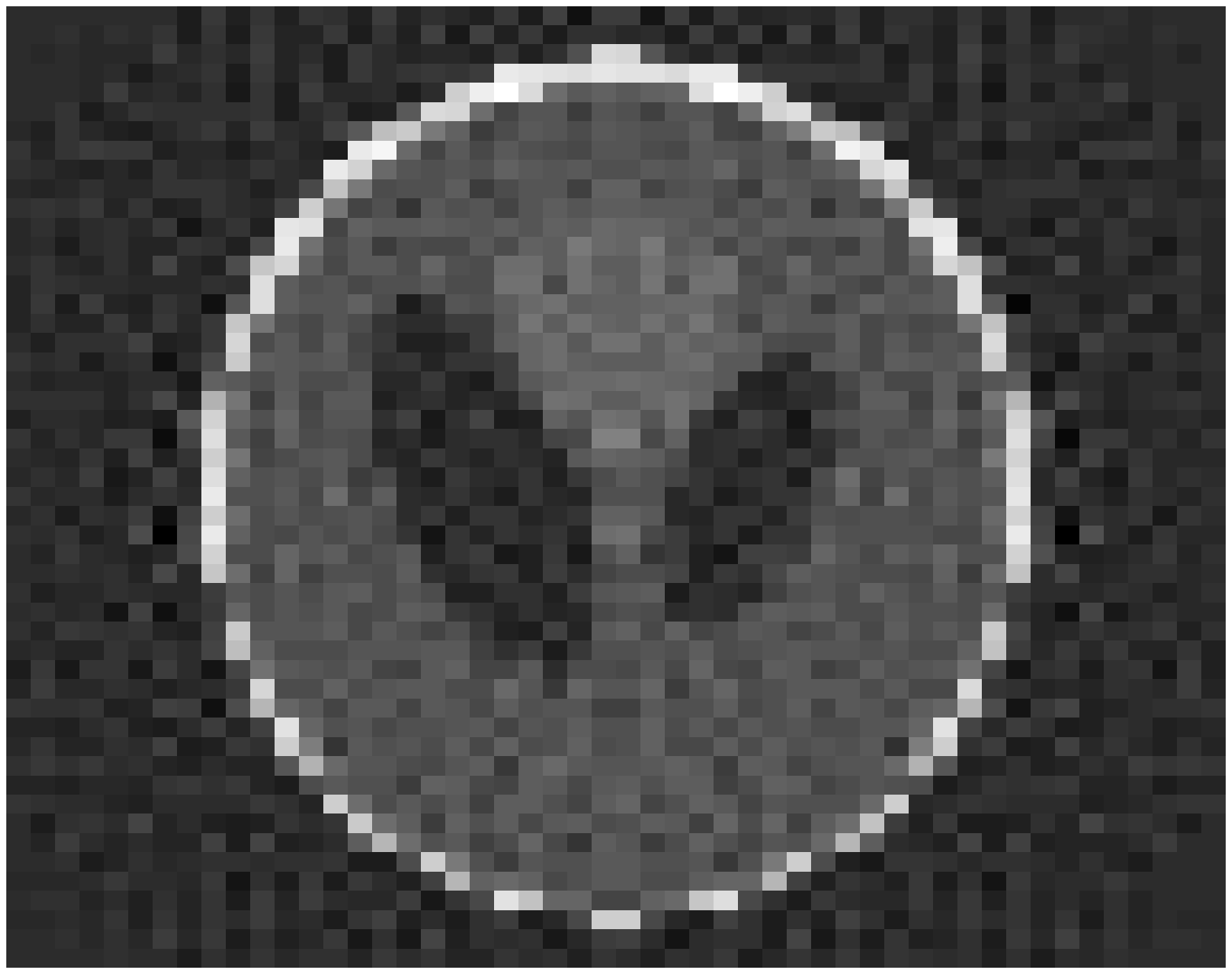}
  \end{minipage}}
  \subfigure[$\delta=0.01$]{
  \begin{minipage}[ht]{.22\linewidth}
      \includegraphics[width=\linewidth]{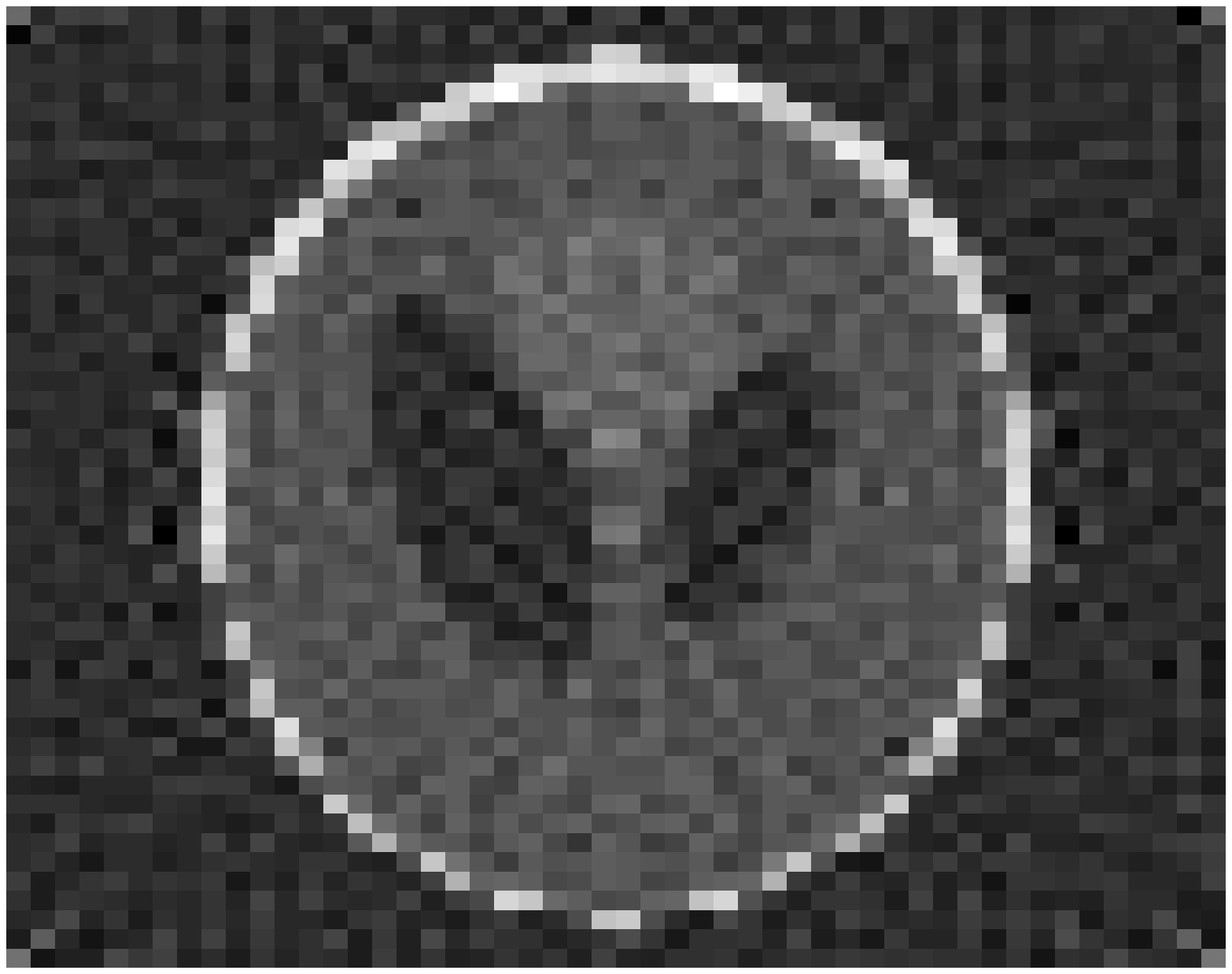}
  \end{minipage}}
    \subfigure[$\delta=0.02$]{
  \begin{minipage}[ht]{.22\linewidth}
      \includegraphics[width=\linewidth]{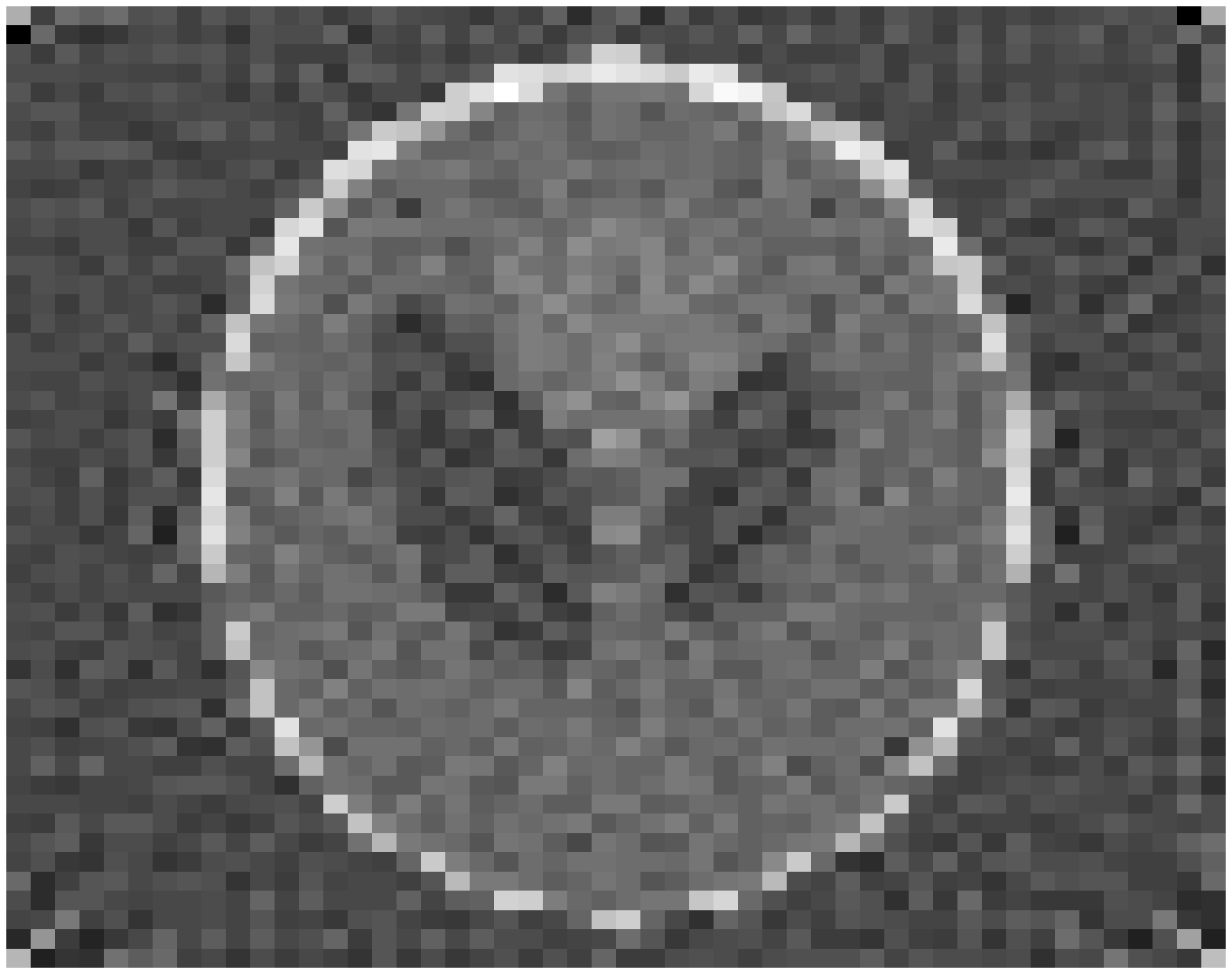}
  \end{minipage}}
  \subfigure[$\delta=0.05$]{
  \begin{minipage}[ht]{.22\linewidth}
      \includegraphics[width=\linewidth]{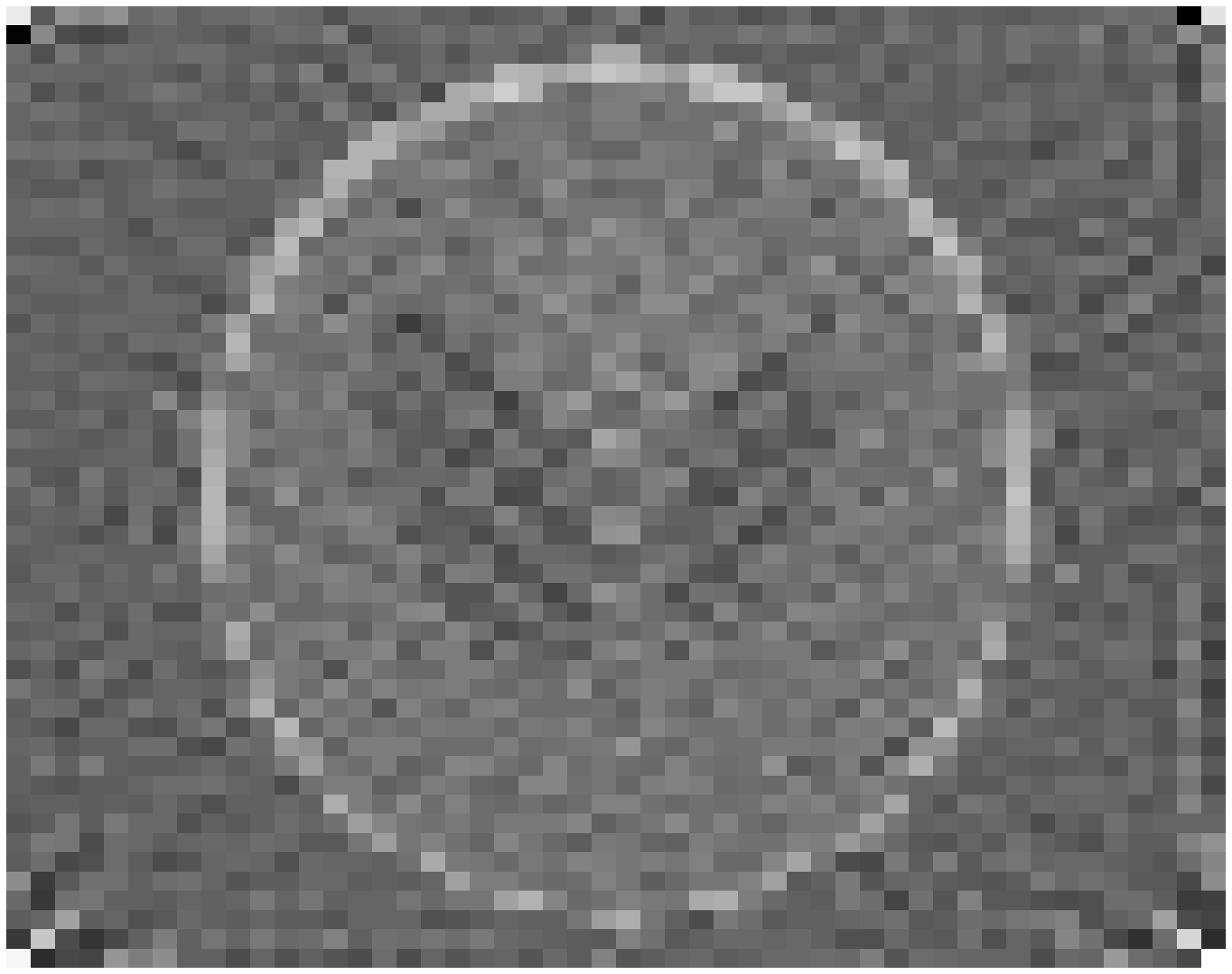}
  \end{minipage}}
\caption{\footnotesize Numerical images of the Kaczmarz method for Model Problem 3}\label{the_Kaczmarz_method_images}
\end{figure}

Figures \ref{the_periodic_Kaczmarz_method_errors} $\sim$ \ref{the_Kaczmarz_Tanabe_method_images} are iterative error curves, residual curves and numerical images of the Kaczmarz-Tanabe method, respectively.
As a comparison, we also give the figures of the Kaczmarz method as shown in Figures \ref{the_Kaczmarz_method_error} $\sim$ \ref{the_Kaczmarz_method_images}. In fact, Figure \ref{the_Kaczmarz_Tanabe_method_images} and Figure \ref{the_Kaczmarz_method_images} are the same. The maximal iteration number of the Kaczmarz-Tanabe method is $K_{\max}=30$, and the maximal iteration number of the Kaczmarz method is $81000$ in order to match the periods of the Kaczmarz-Tanabe method.

\begin{figure}[!hbt]
  \centering
  \subfigure[$\delta=0$]{
  \begin{minipage}[ht]{.22\linewidth}
      \includegraphics[width=\linewidth]{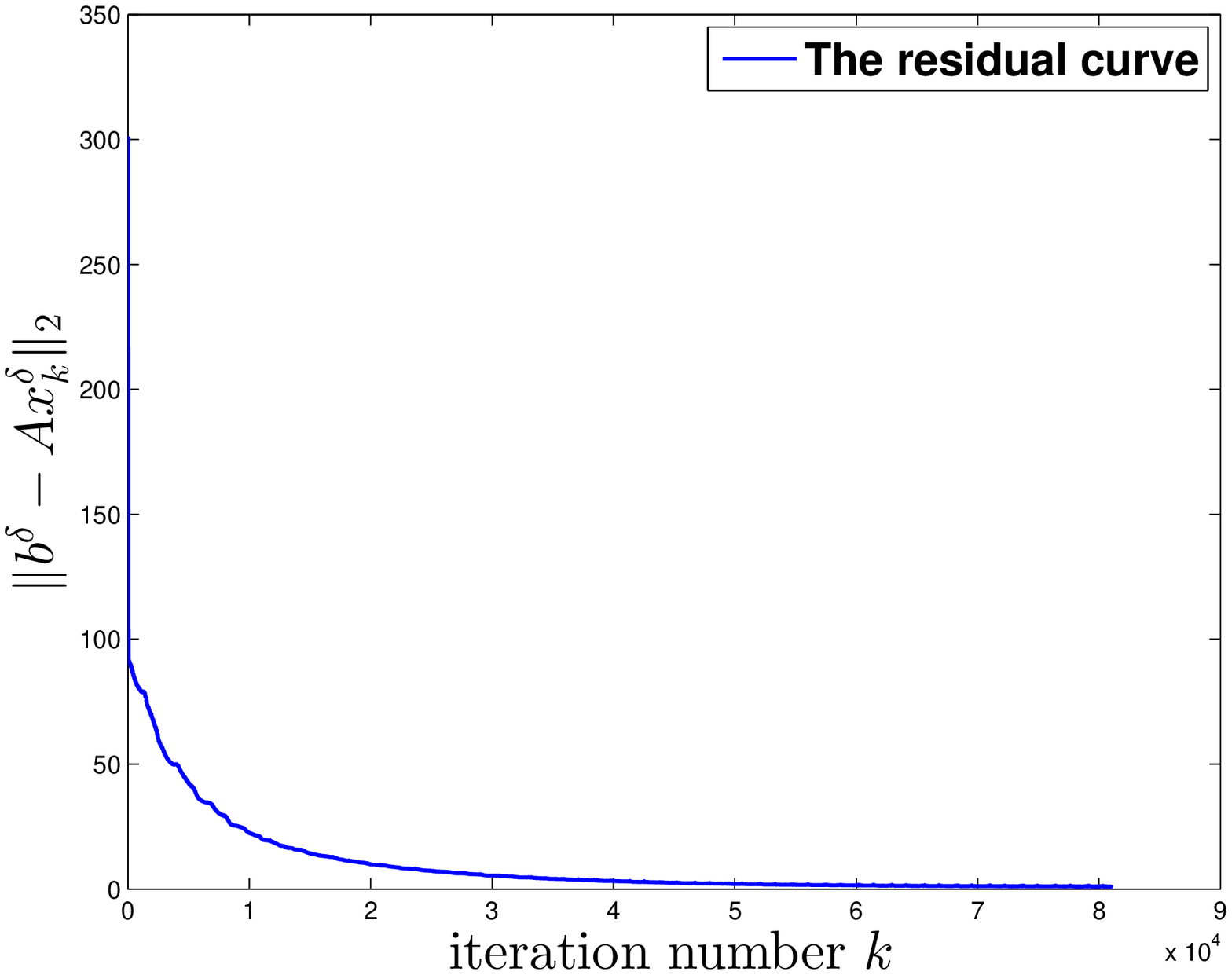}
  \end{minipage}}
  \subfigure[$\delta=0.01$]{
  \begin{minipage}[ht]{.22\linewidth}
      \includegraphics[width=\linewidth]{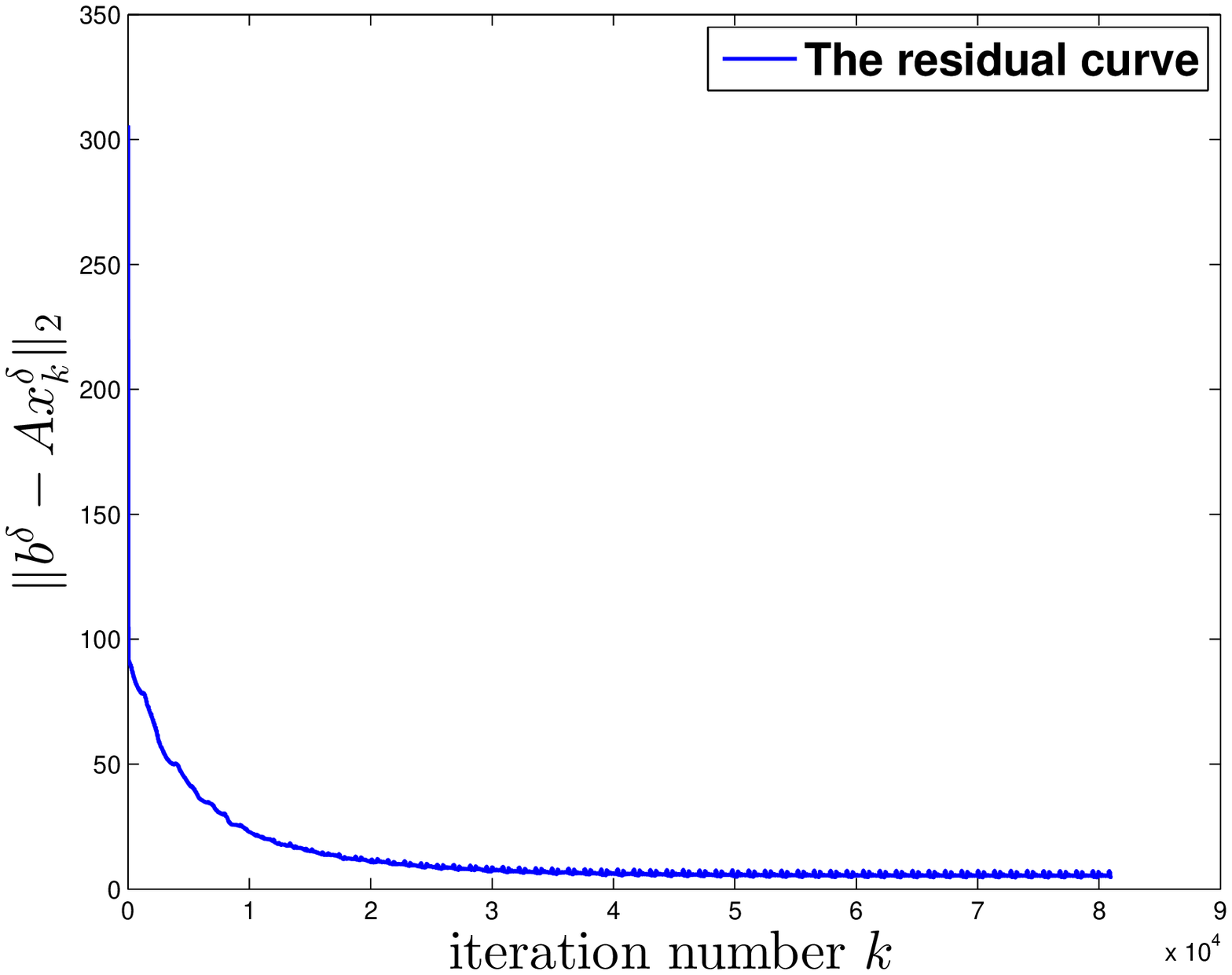}
  \end{minipage}}
    \subfigure[$\delta=0.02$]{
  \begin{minipage}[ht]{.22\linewidth}
      \includegraphics[width=\linewidth]{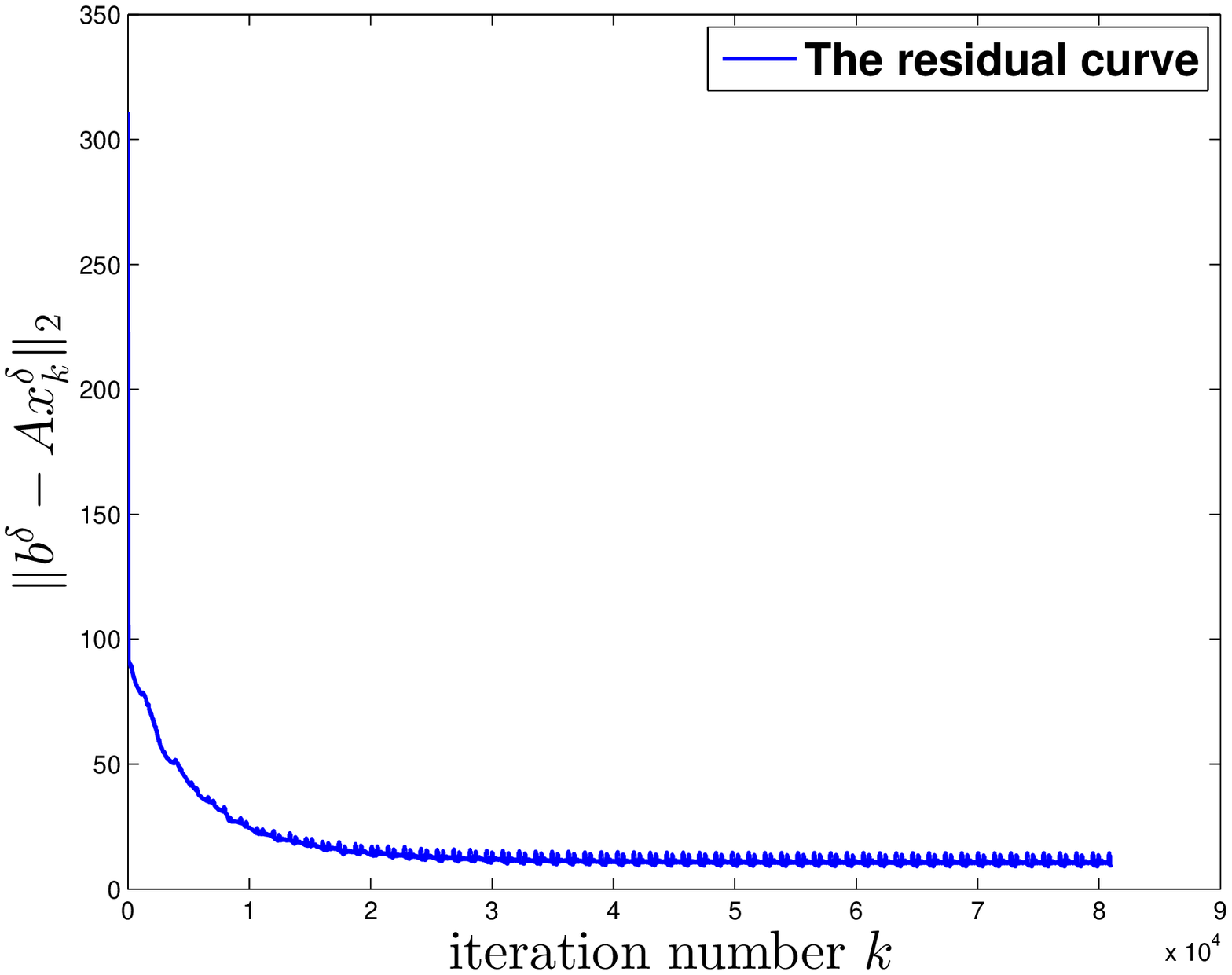}
  \end{minipage}}
  \subfigure[$\delta=0.05$]{
  \begin{minipage}[ht]{.22\linewidth}
      \includegraphics[width=\linewidth]{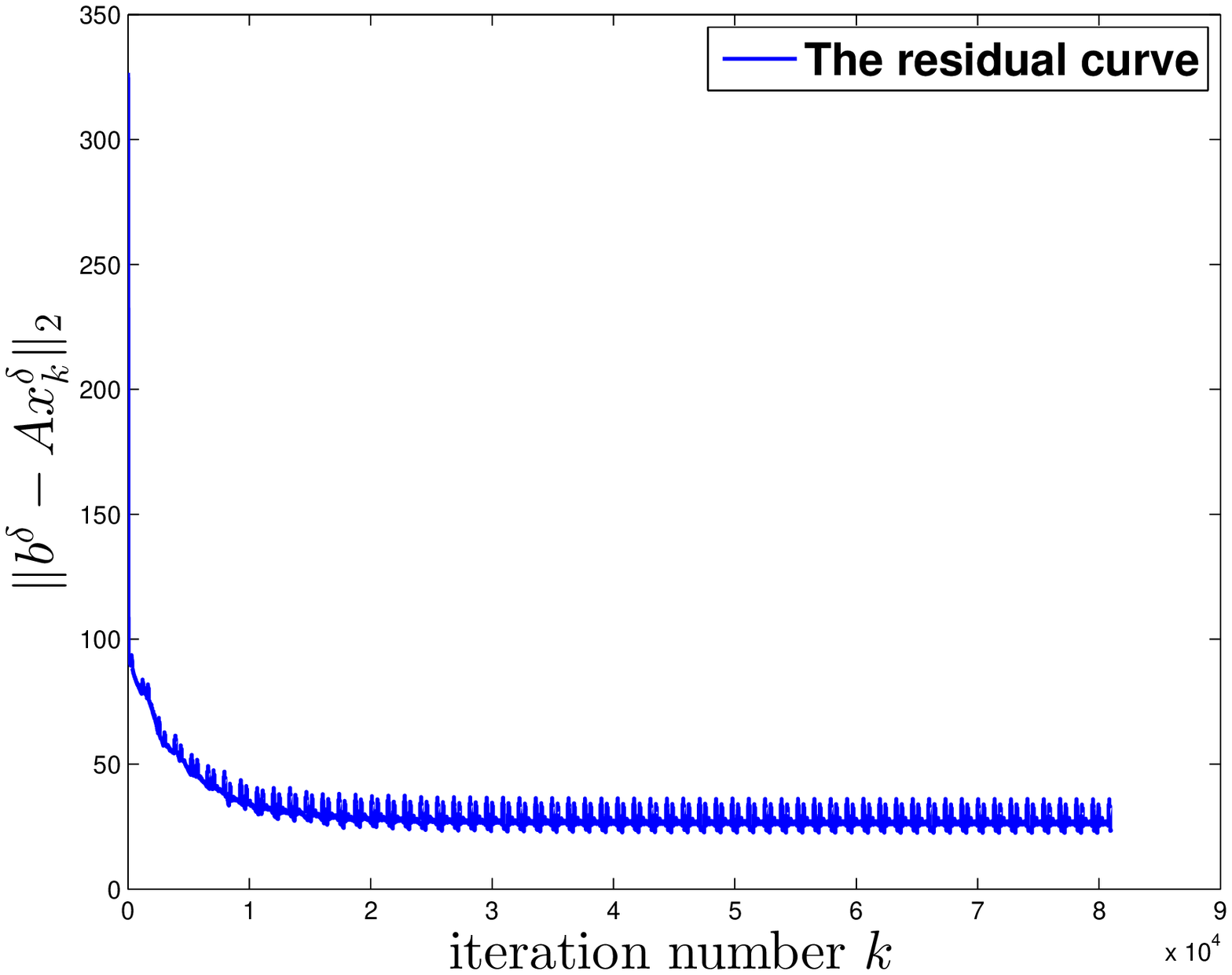}
  \end{minipage}}
  \caption{\footnotesize The residual curves of the Kaczmarz method for Model Problem 3}\label{the_Kaczmarz_method_residue}
\end{figure}
%

Comparing Figures \ref{the_Kaczmarz_method_error} $\sim$ \ref{the_Kaczmarz_method_residue}, we find that the iterative error and residual curves of the Kaczmarz method fluctuate violently, while the iterative error and residual curves of the Kaczmarz-Tanabe method are relatively smooth. Consequently, for perturbed cases, the properties of the Kaczmarz-Tanabe method are better than those of the Kaczmarz method.

\begin{figure}[!hbt]
  \centering
  \subfigure[$\eta=0$]{
  \begin{minipage}[ht]{.22\linewidth}
      \includegraphics[width=\linewidth]{figures/Phantom/Kaczmarz.Tanabe/Phantom_error_0.eps}
  \end{minipage}}
  \subfigure[$\eta=0.023$]{
  \begin{minipage}[ht]{.22\linewidth}
      \includegraphics[width=\linewidth]{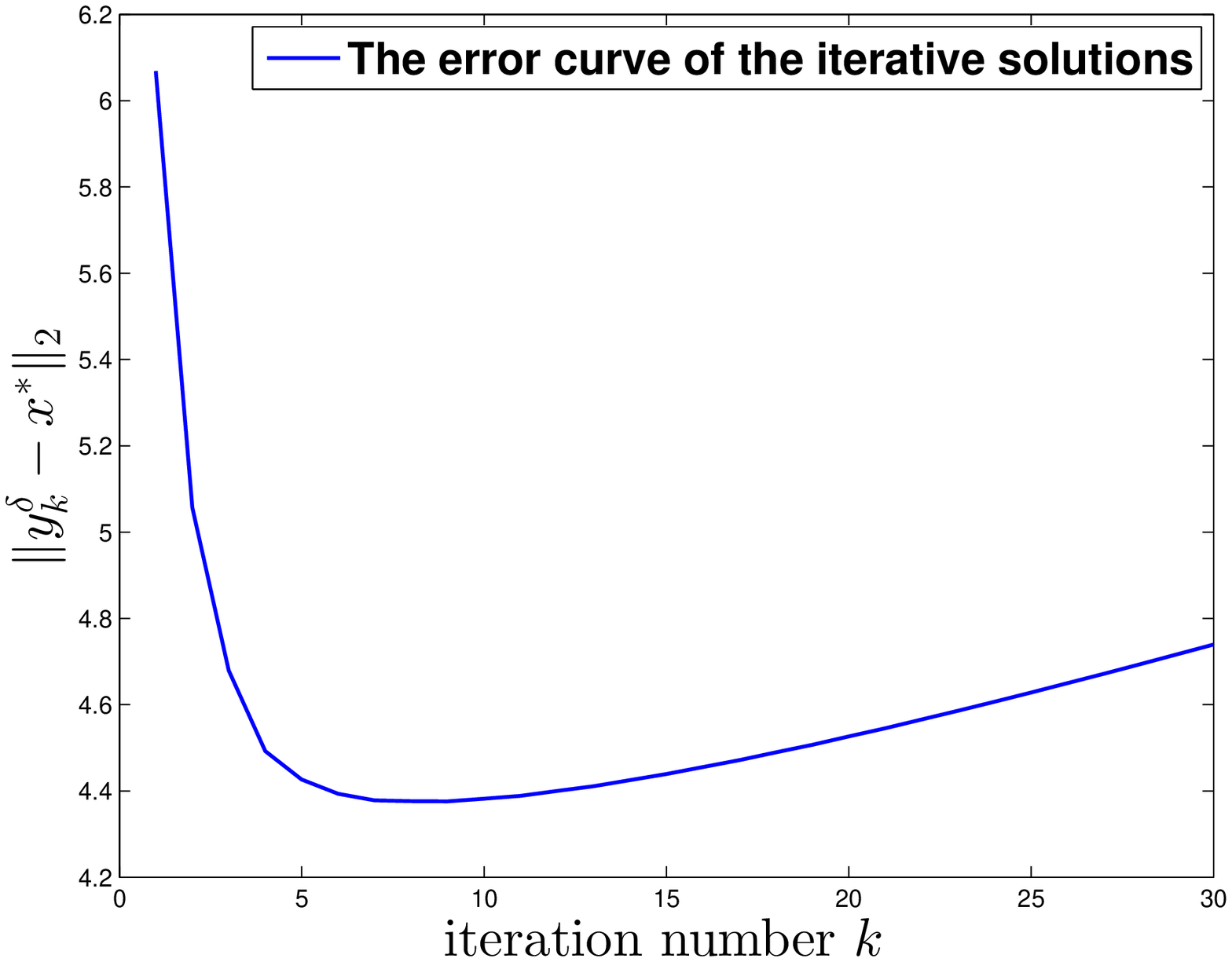}
  \end{minipage}}
    \subfigure[$\eta=0.046$]{
  \begin{minipage}[ht]{.22\linewidth}
      \includegraphics[width=\linewidth]{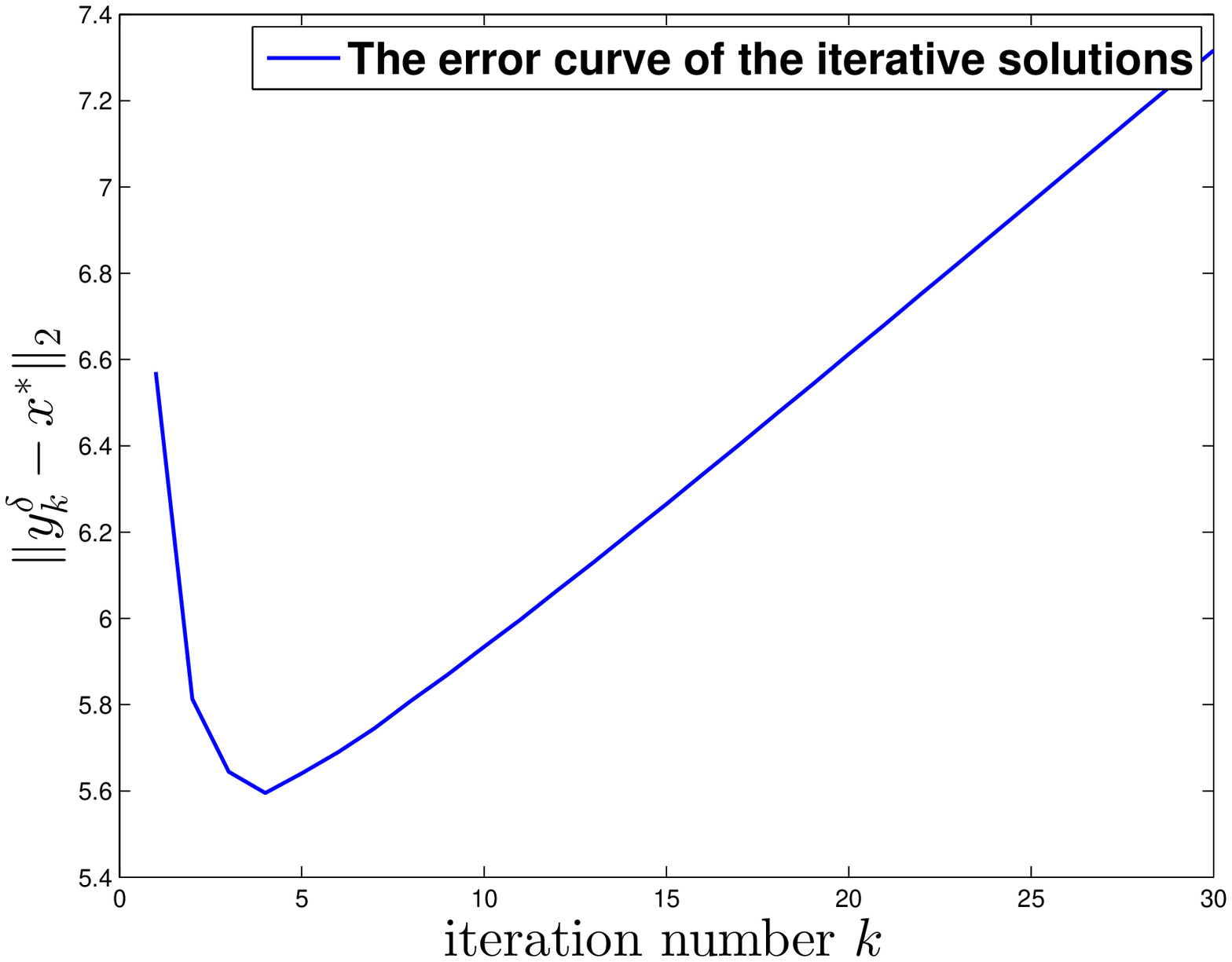}
  \end{minipage}}
  \subfigure[$\eta=0.115$]{
  \begin{minipage}[ht]{.22\linewidth}
      \includegraphics[width=\linewidth]{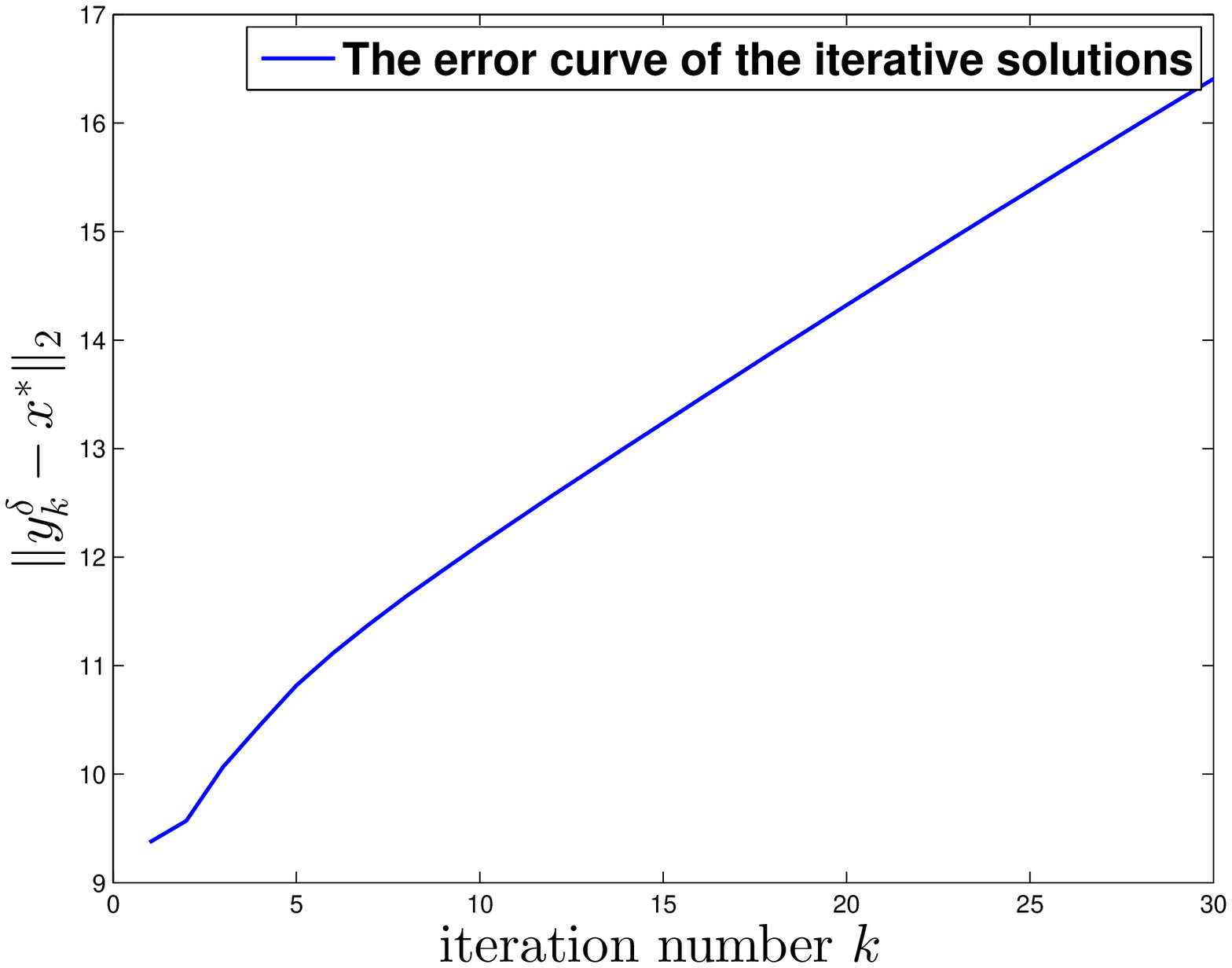}
  \end{minipage}}
  \caption{\footnotesize The error curves of the Kaczmarz-Tanabe method for Model Problem 3 with Gaussian noise}\label{the_Kaczmarz_Tanabe_method_error_Gaussian_noise}
\end{figure}

\begin{figure}[!hbt]
  \centering
 \subfigure[$\eta=0$]{
  \begin{minipage}[ht]{.22\linewidth}
      \includegraphics[width=\linewidth]{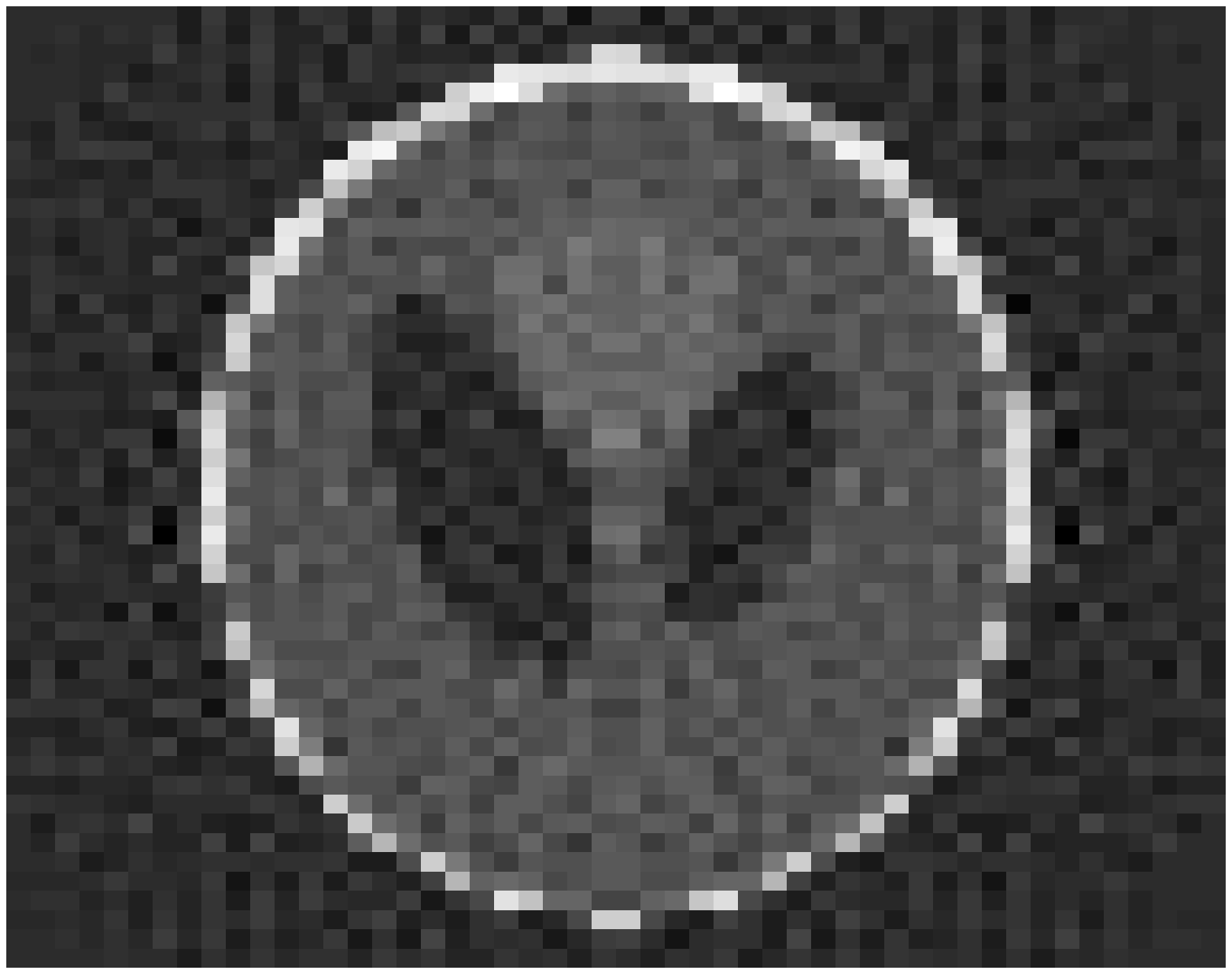}
  \end{minipage}}
  \subfigure[$\eta=0.023$]{
  \begin{minipage}[ht]{.22\linewidth}
      \includegraphics[width=\linewidth]{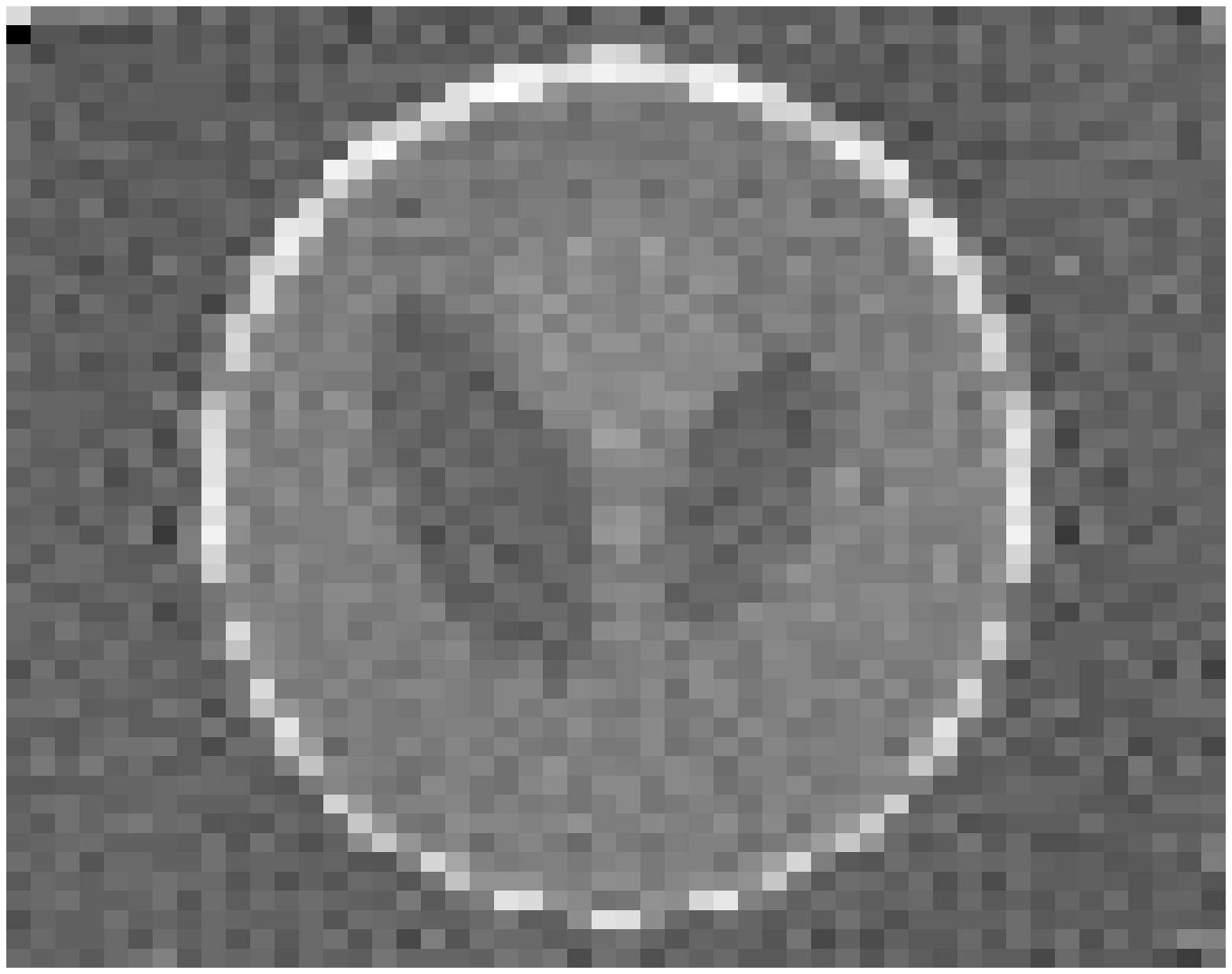}
  \end{minipage}}
    \subfigure[$\eta=0.046$]{
  \begin{minipage}[ht]{.22\linewidth}
      \includegraphics[width=\linewidth]{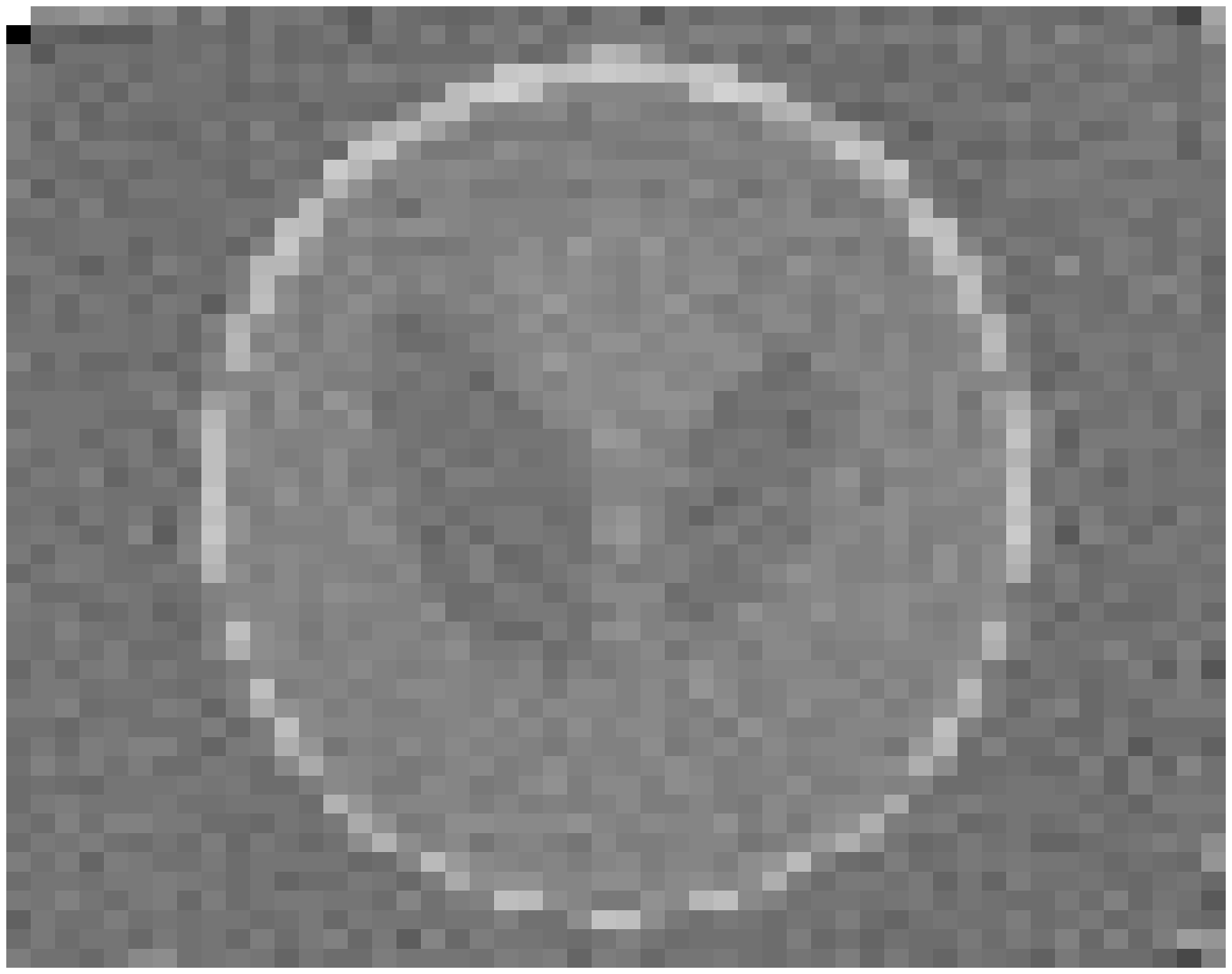}
  \end{minipage}}
  \subfigure[$\eta=0.115$]{
  \begin{minipage}[ht]{.22\linewidth}
      \includegraphics[width=\linewidth]{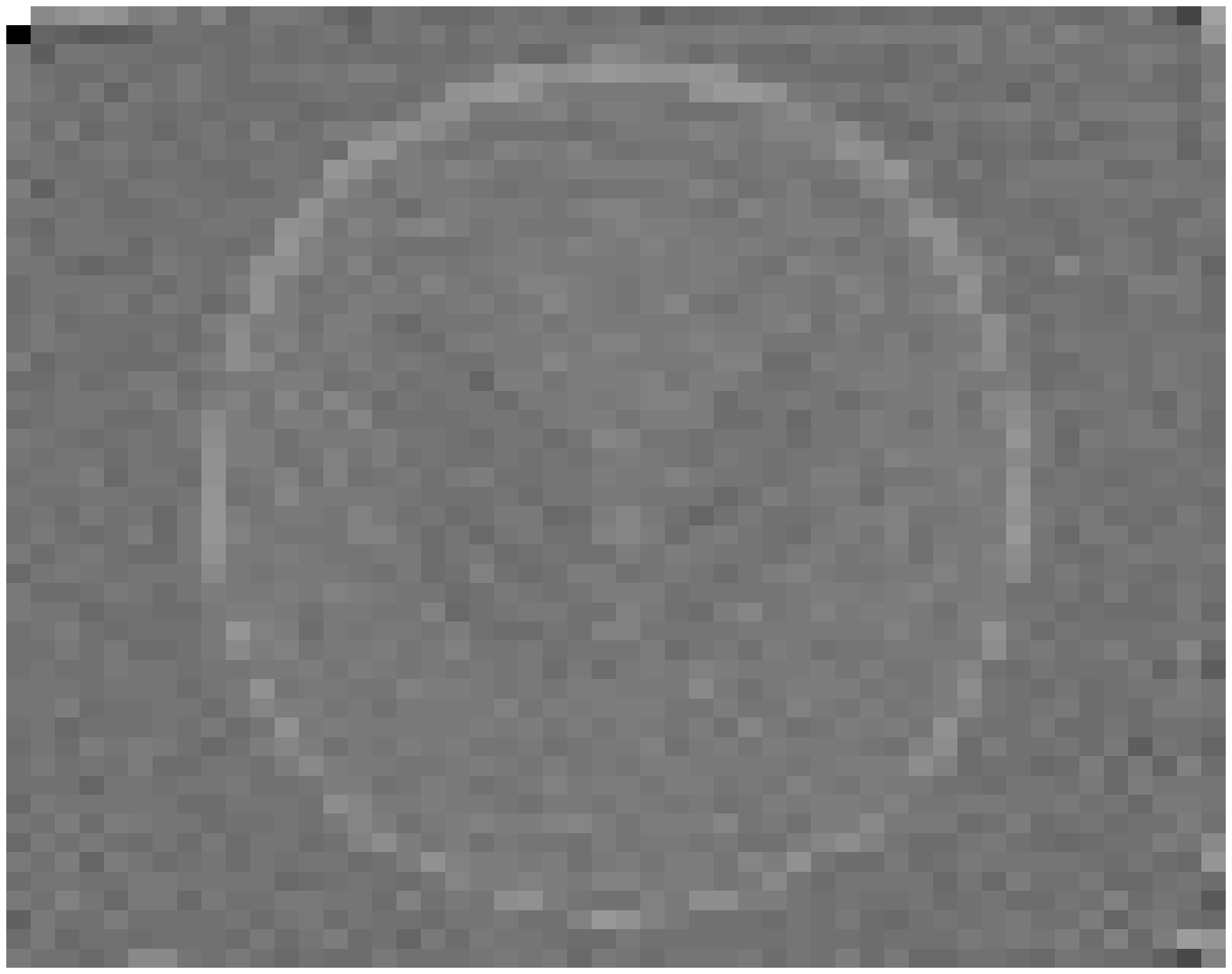}
  \end{minipage}}
\caption{\footnotesize Numerical images of the Kaczmarz-Tanabe method for Model Problem 3 with Gaussian noise}\label{the_Kaczmarz_Tanabe_method_Gaussian_noise_images}
\end{figure}

\begin{figure}[!hbt]
  \centering
  \subfigure[$\eta=0$]{
  \begin{minipage}[ht]{.22\linewidth}
      \includegraphics[width=\linewidth]{figures/Phantom/Kaczmarz/Phantom_error_0.eps}
  \end{minipage}}
  \subfigure[$\eta=0.023$]{
  \begin{minipage}[ht]{.22\linewidth}
      \includegraphics[width=\linewidth]{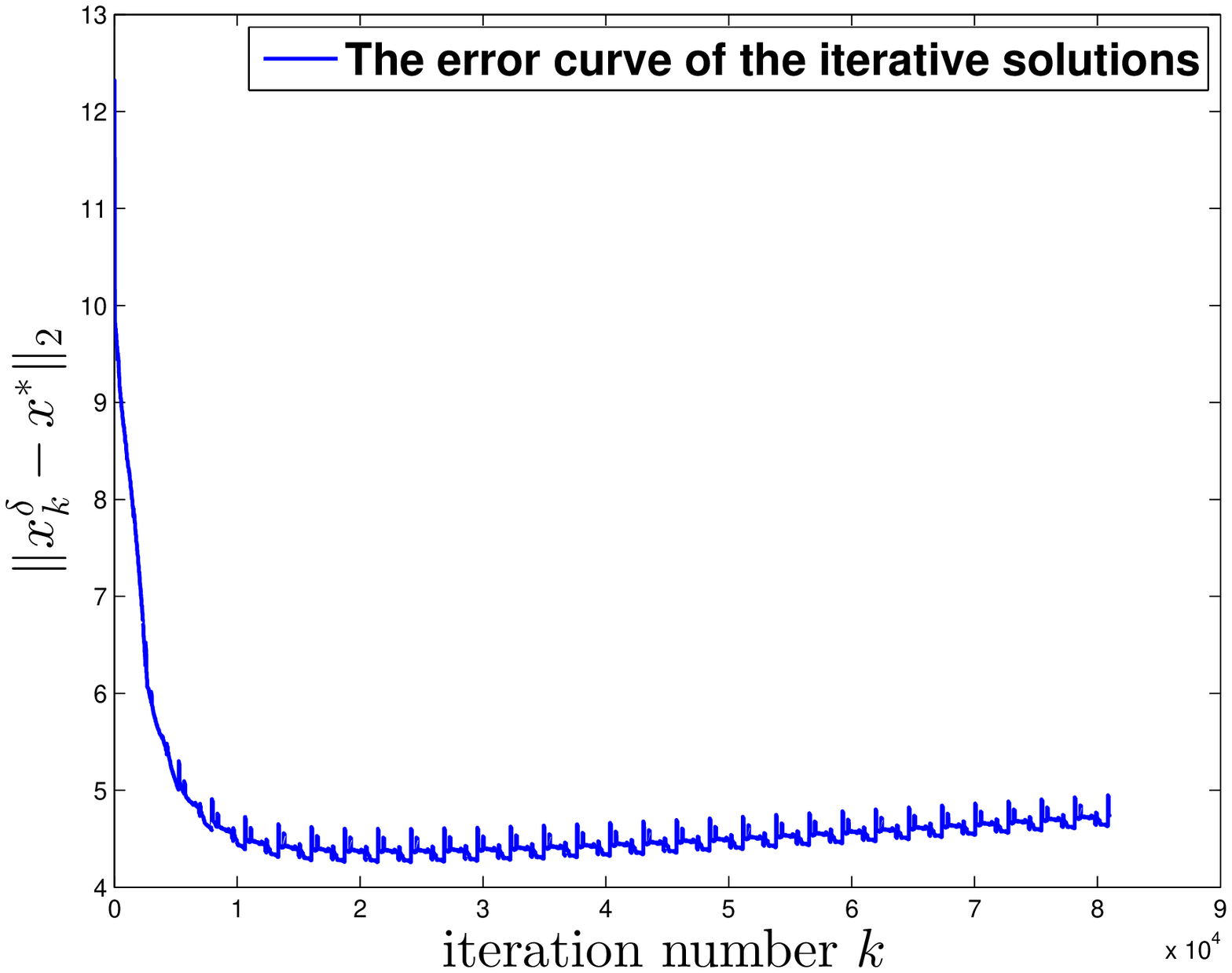}
  \end{minipage}}
    \subfigure[$\eta=0.046$]{
  \begin{minipage}[ht]{.22\linewidth}
      \includegraphics[width=\linewidth]{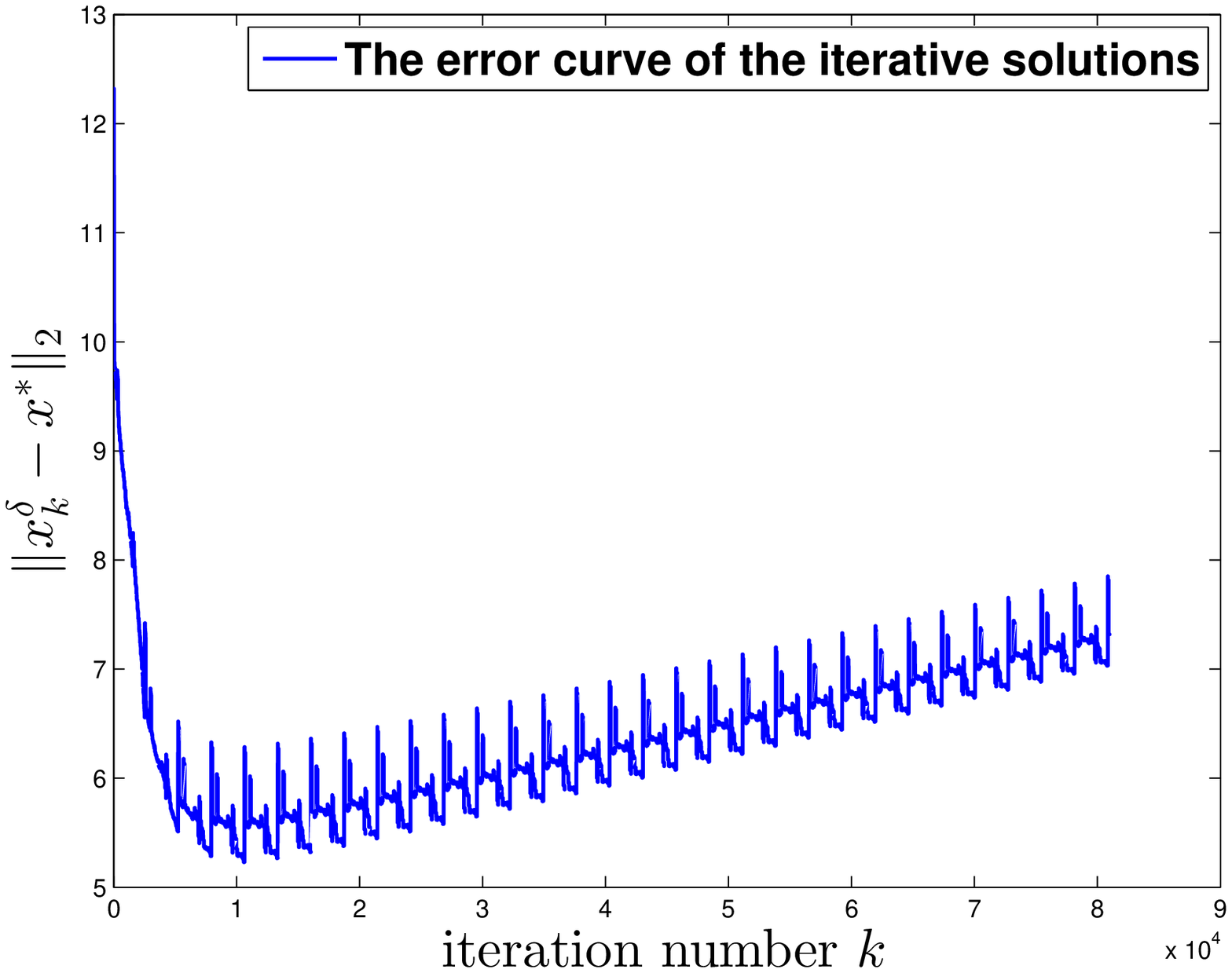}
  \end{minipage}}
  \subfigure[$\eta=0.115$]{
  \begin{minipage}[ht]{.22\linewidth}
      \includegraphics[width=\linewidth]{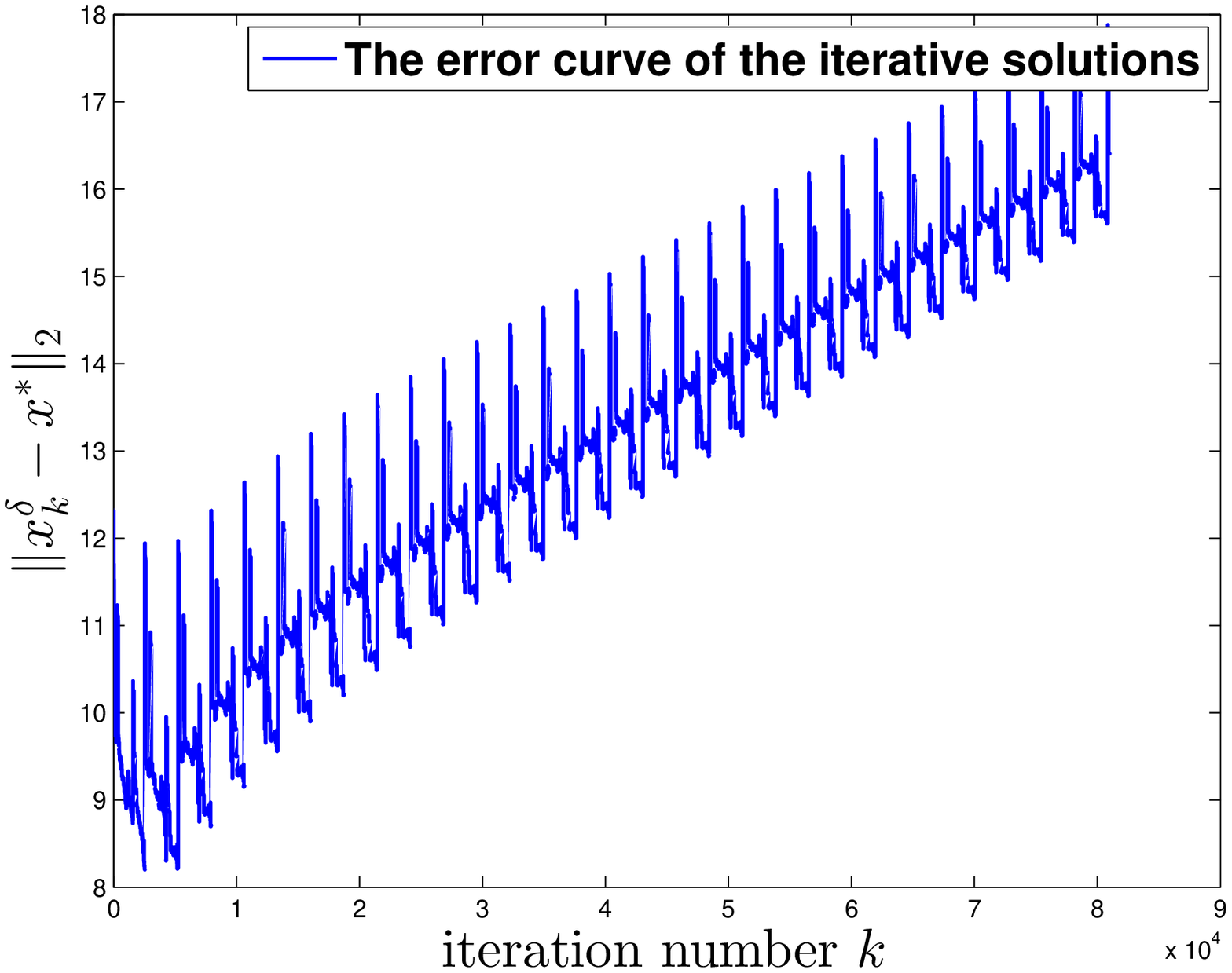}
  \end{minipage}}
  \caption{\footnotesize The error curves of the Kaczmarz method for Model Problem 3 with Gaussian noise}\label{the_Kaczmarz_method_error_Gaussian_noise}
\end{figure}

\begin{figure}[!hbt]
  \centering
 \subfigure[$\eta=0$]{
  \begin{minipage}[ht]{.22\linewidth}
      \includegraphics[width=\linewidth]{figures/Phantom/Kaczmarz/Phantom_numerical_image_0.eps}
  \end{minipage}}
  \subfigure[$\eta=0.023$]{
  \begin{minipage}[ht]{.22\linewidth}
      \includegraphics[width=\linewidth]{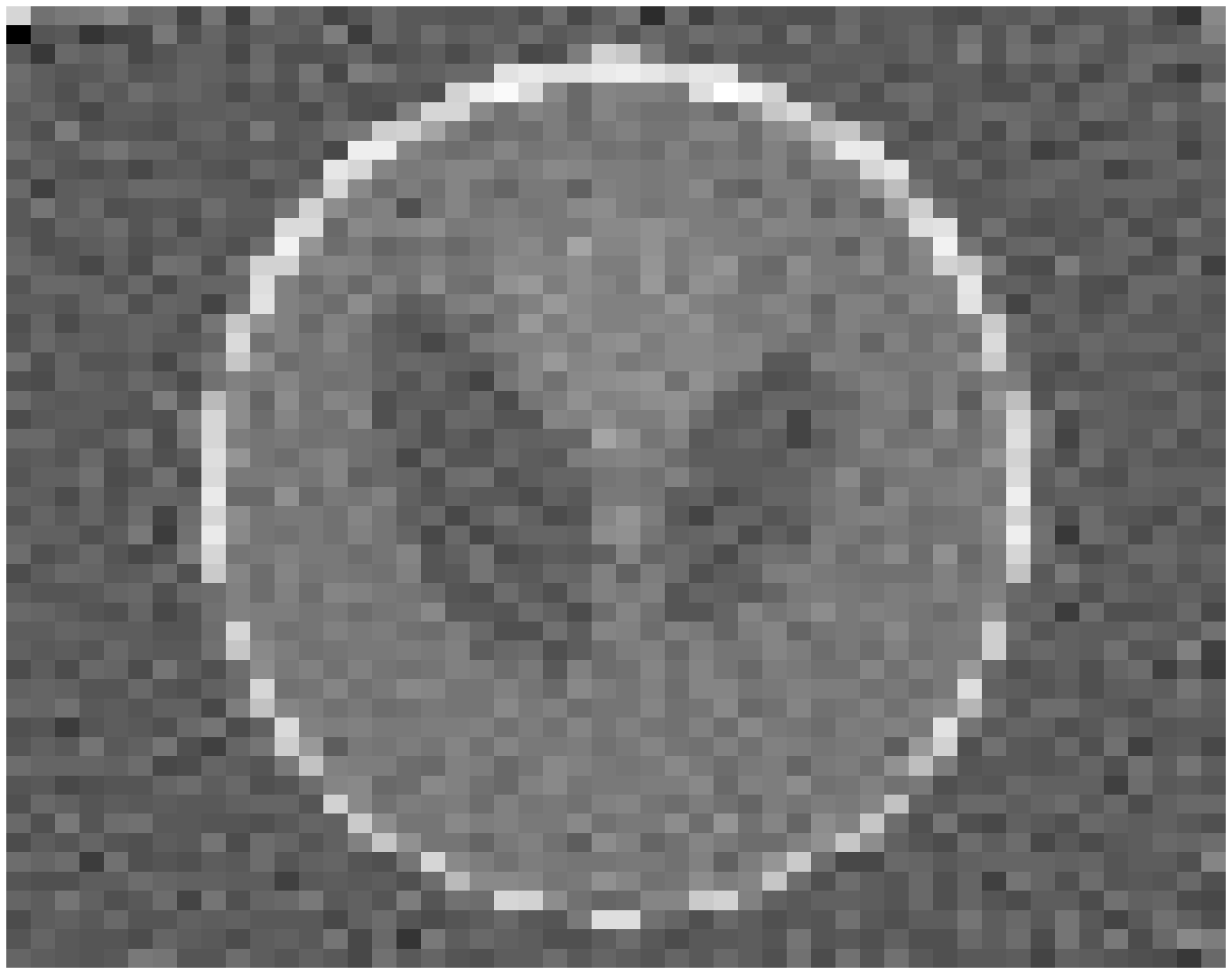}
  \end{minipage}}
    \subfigure[$\eta=0.046$]{
  \begin{minipage}[ht]{.22\linewidth}
      \includegraphics[width=\linewidth]{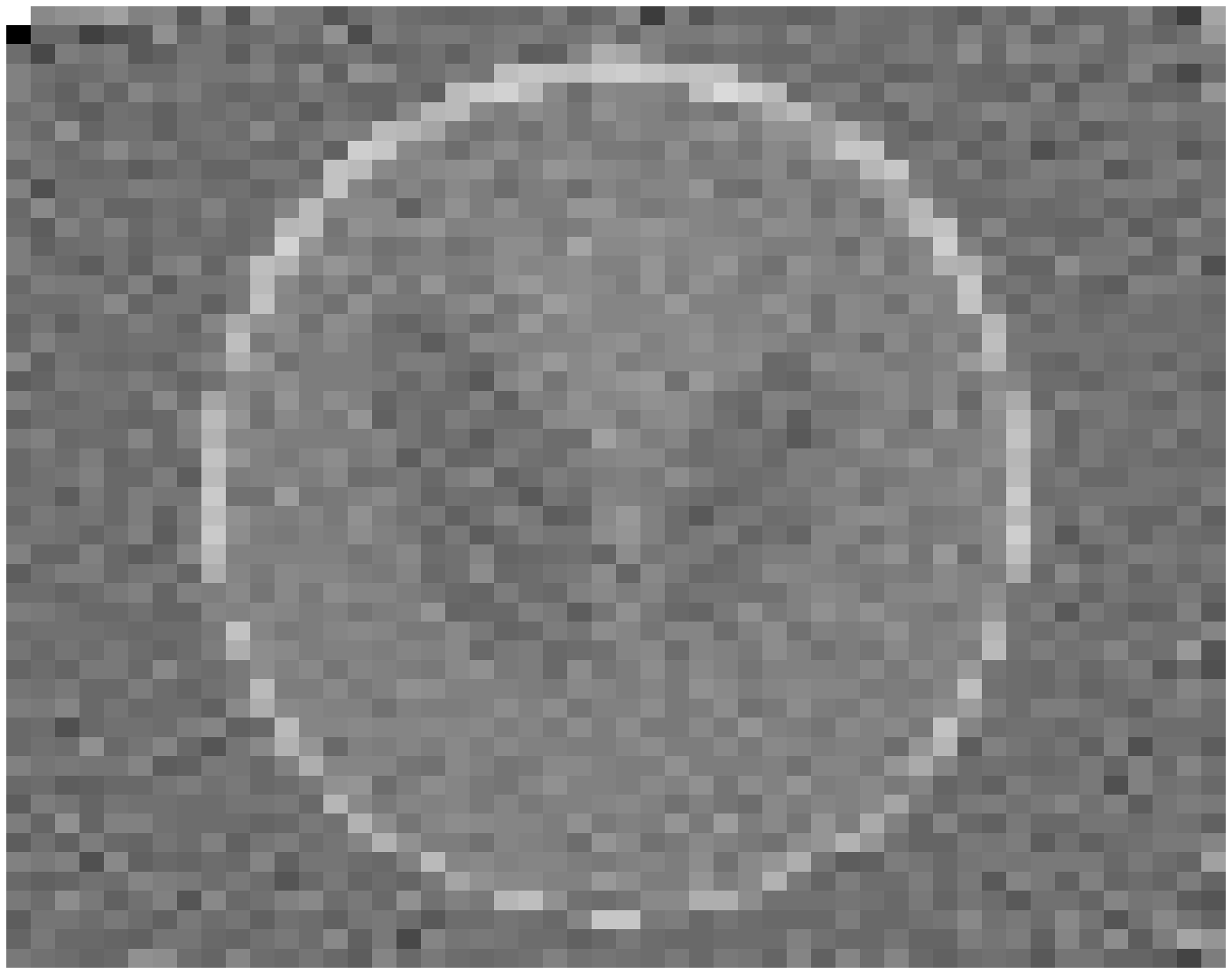}
  \end{minipage}}
  \subfigure[$\eta=0.115$]{
  \begin{minipage}[ht]{.22\linewidth}
      \includegraphics[width=\linewidth]{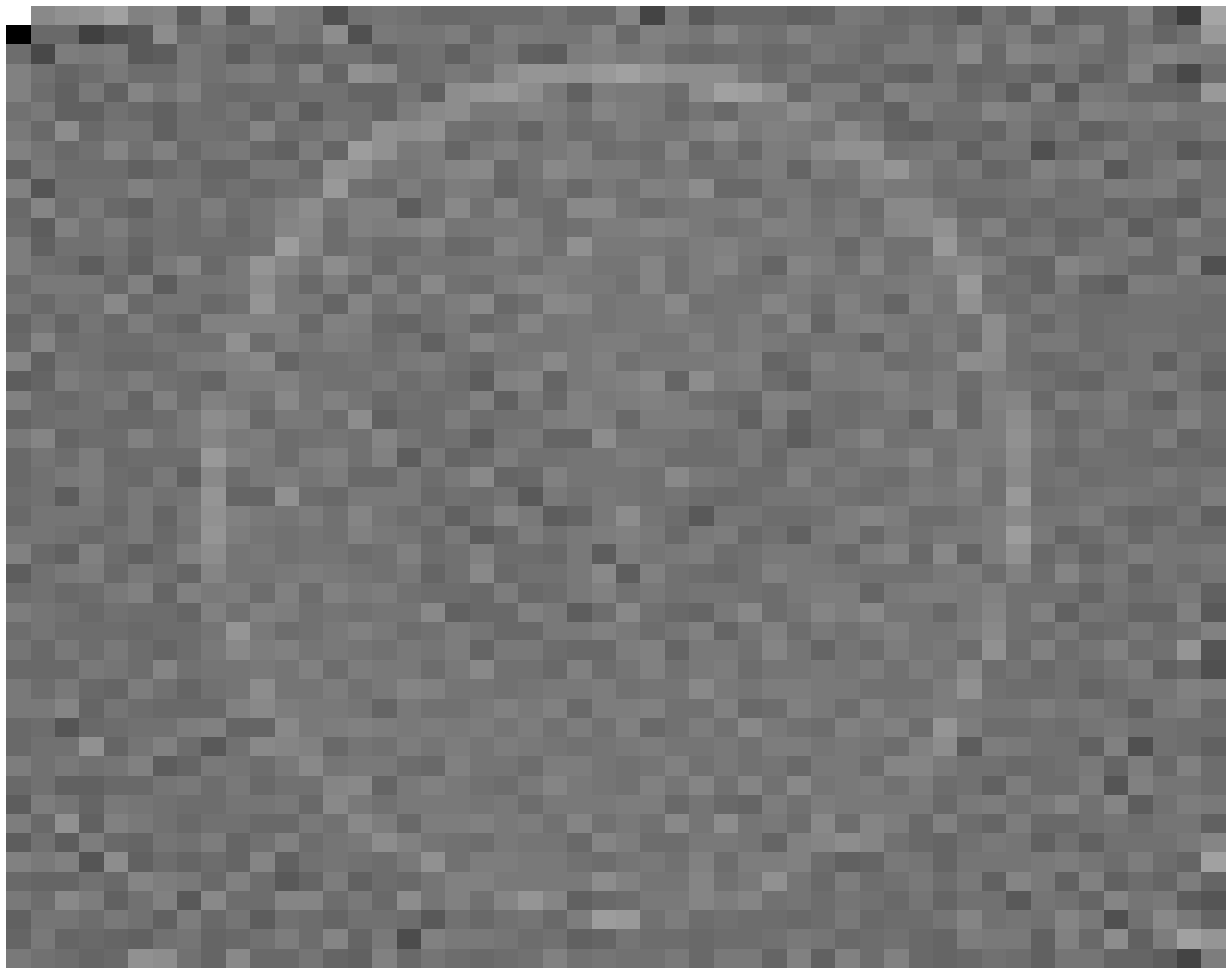}
  \end{minipage}}
\caption{\footnotesize Numerical images of the Kaczmarz method for Model Problem 3 with Gaussian noise}\label{the_Kaczmarz_method_Gaussian_noise_images}
\end{figure}

We next superimpose the Gaussian noise with mean value $0$ and variance $0, 0.023\|b\|, 0.046\|b\|$ and $0.115\|b\|$ on $b$, numerical results are shown in Figures \ref{the_Kaczmarz_Tanabe_method_error_Gaussian_noise} $\sim$ \ref{the_Kaczmarz_method_Gaussian_noise_images} and the corresponding numerical results are marked with $\eta=0, 0.023, 0.046$ and $0.115$. Figures \ref{the_Kaczmarz_Tanabe_method_error_Gaussian_noise} $\sim$ \ref{the_Kaczmarz_Tanabe_method_Gaussian_noise_images} are for the Kaczmarz-Tanabe method and Figures \ref{the_Kaczmarz_method_error_Gaussian_noise} $\sim$ \ref{the_Kaczmarz_method_Gaussian_noise_images} are for the Kaczmarz method. The absolute and relative errors are same to $\delta=0, 0.01, 0.02$ and $0.05$ in Table \ref{errror_level_1_for_model_problem3}, respectively.

\begin{figure}[!hbt]
  \centering
  \subfigure[$\eta=0$]{
  \begin{minipage}[ht]{.22\linewidth}
      \includegraphics[width=\linewidth]{figures/Phantom/Kaczmarz.Tanabe/Phantom_residue_0.eps}
  \end{minipage}}
  \subfigure[$\eta=0.023$]{
  \begin{minipage}[ht]{.22\linewidth}
      \includegraphics[width=\linewidth]{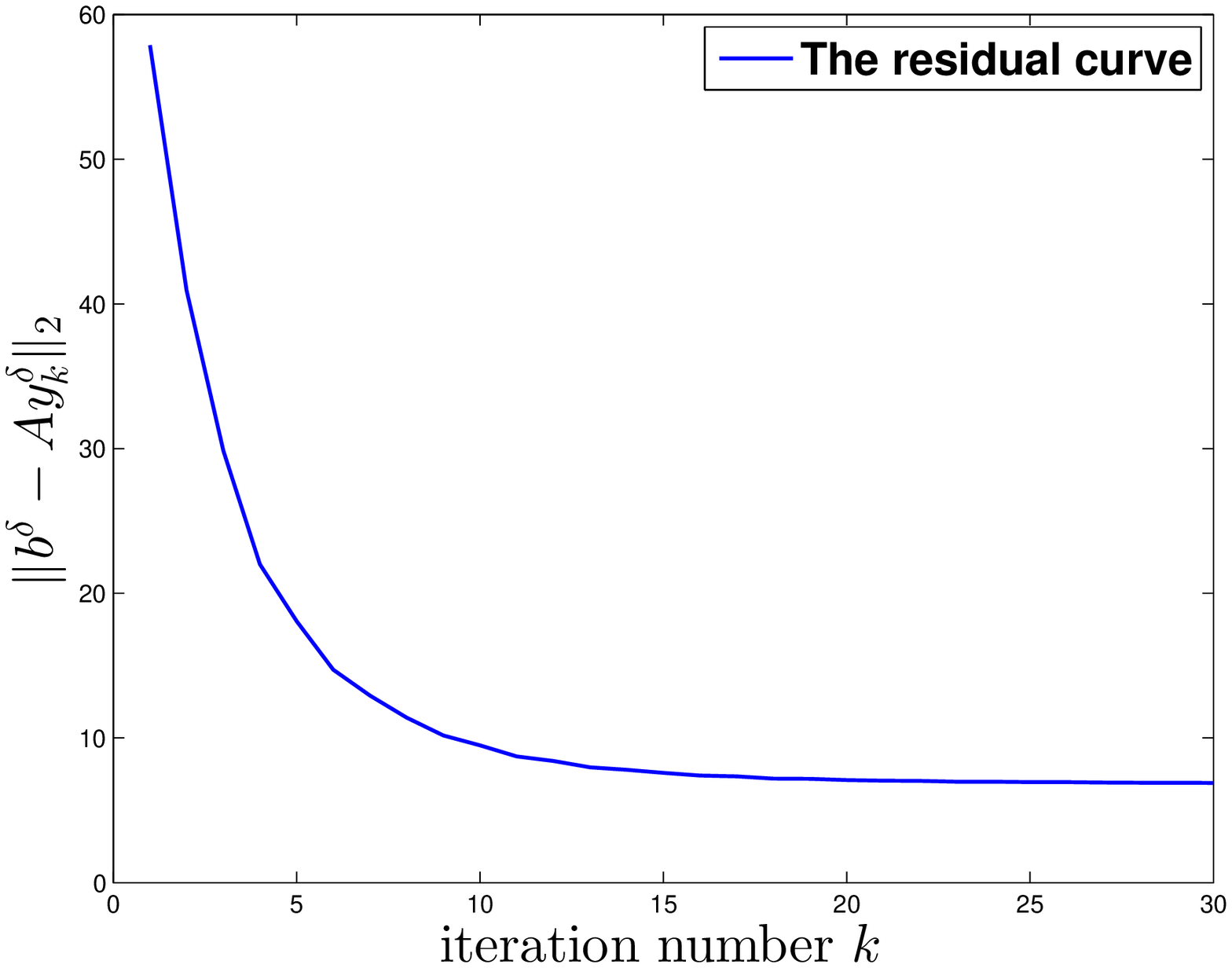}
  \end{minipage}}
    \subfigure[$\eta=0.046$]{
  \begin{minipage}[ht]{.22\linewidth}
      \includegraphics[width=\linewidth]{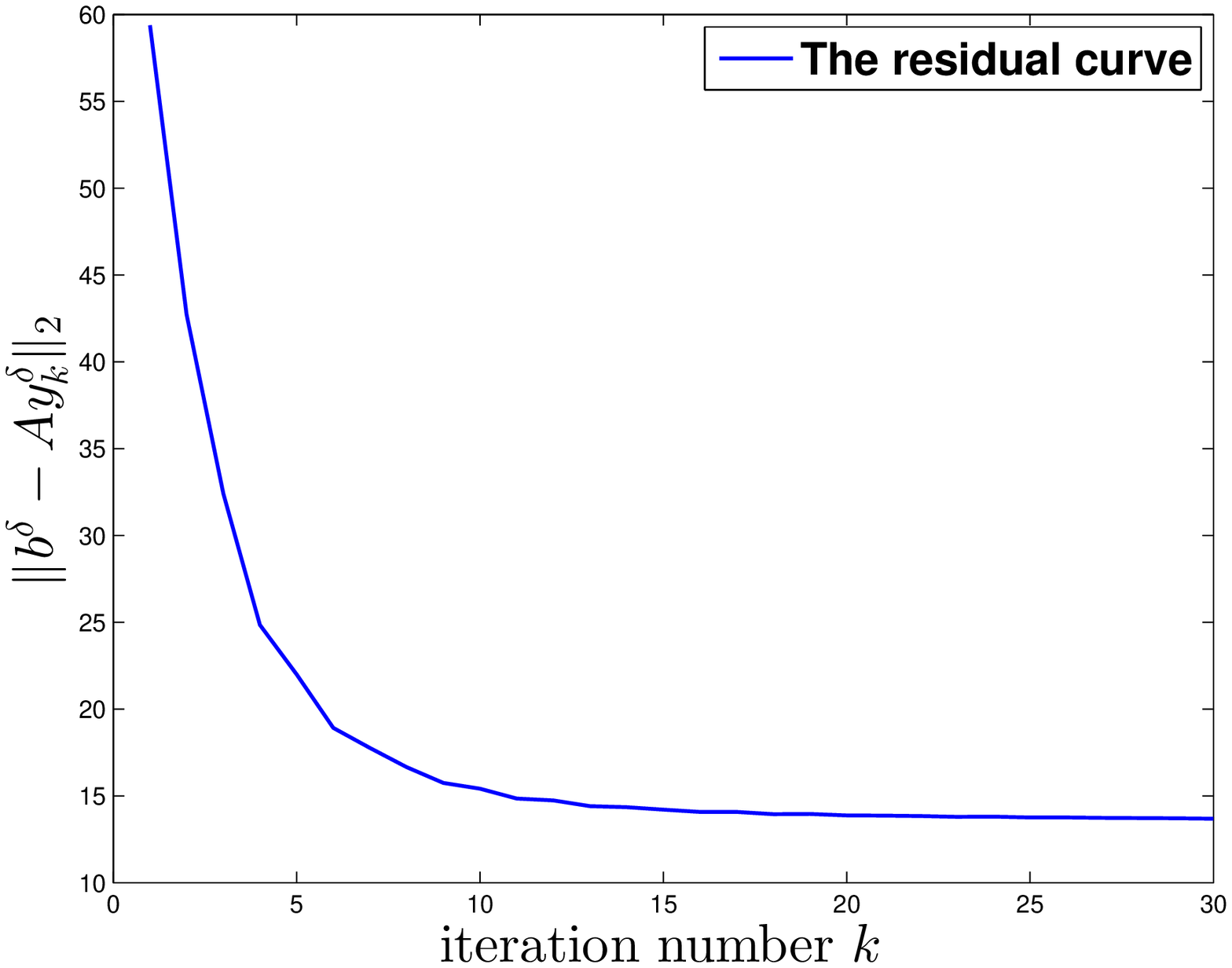}
  \end{minipage}}
  \subfigure[$\eta=0.115$]{
  \begin{minipage}[ht]{.22\linewidth}
      \includegraphics[width=\linewidth]{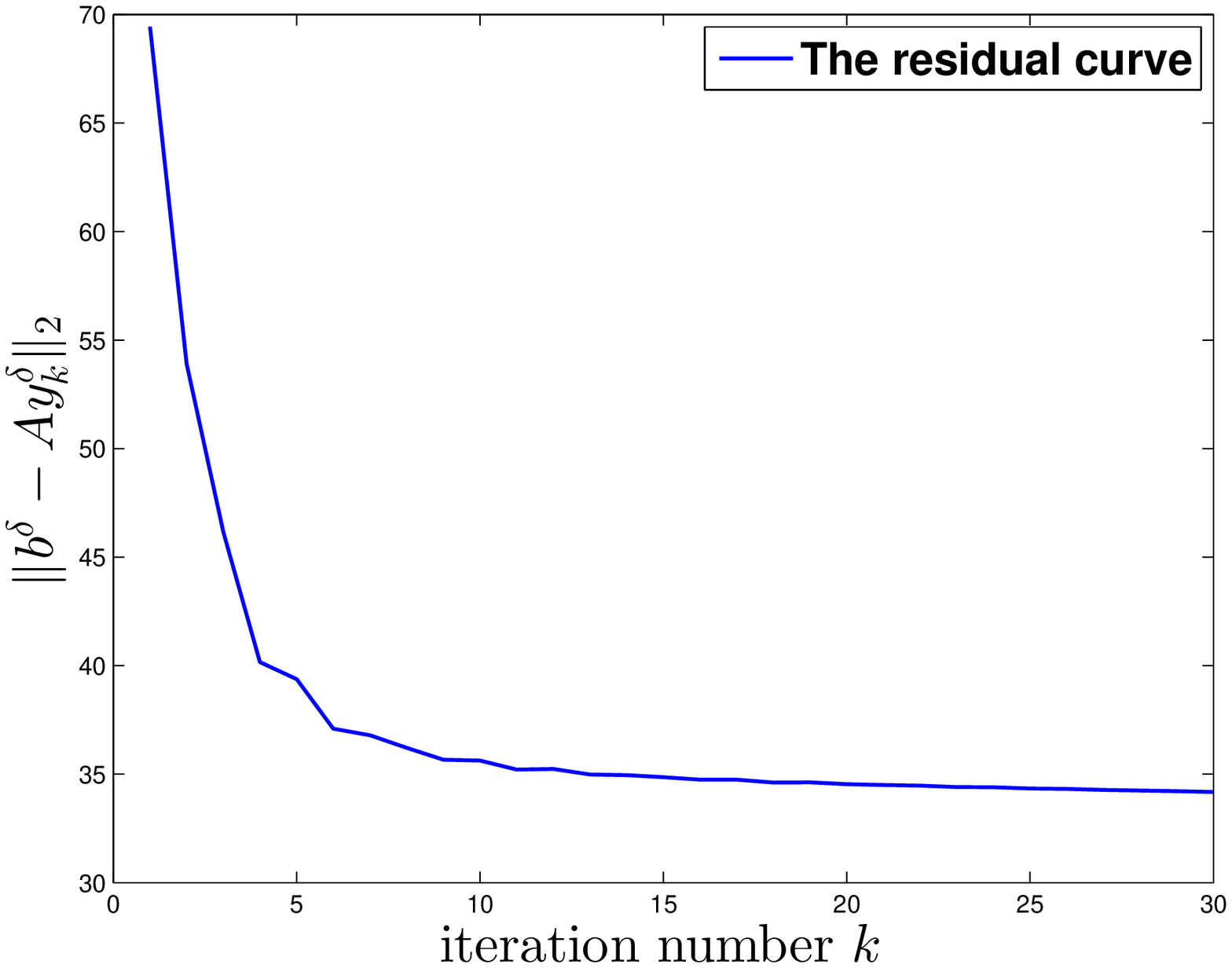}
  \end{minipage}}
  \caption{\footnotesize The residual curves of the Kaczmarz-Tanabe method for Model Problem 3 with Gaussian noise}\label{the_Kaczmarz_Tanabe_method_residue_Gaussian}
\end{figure}
%

\begin{figure}[!hbt]
  \centering
  \subfigure[$\eta=0$]{
  \begin{minipage}[ht]{.22\linewidth}
      \includegraphics[width=\linewidth]{figures/Phantom/Kaczmarz/Phantom_residue_0.eps}
  \end{minipage}}
  \subfigure[$\eta=0.023$]{
  \begin{minipage}[ht]{.22\linewidth}
      \includegraphics[width=\linewidth]{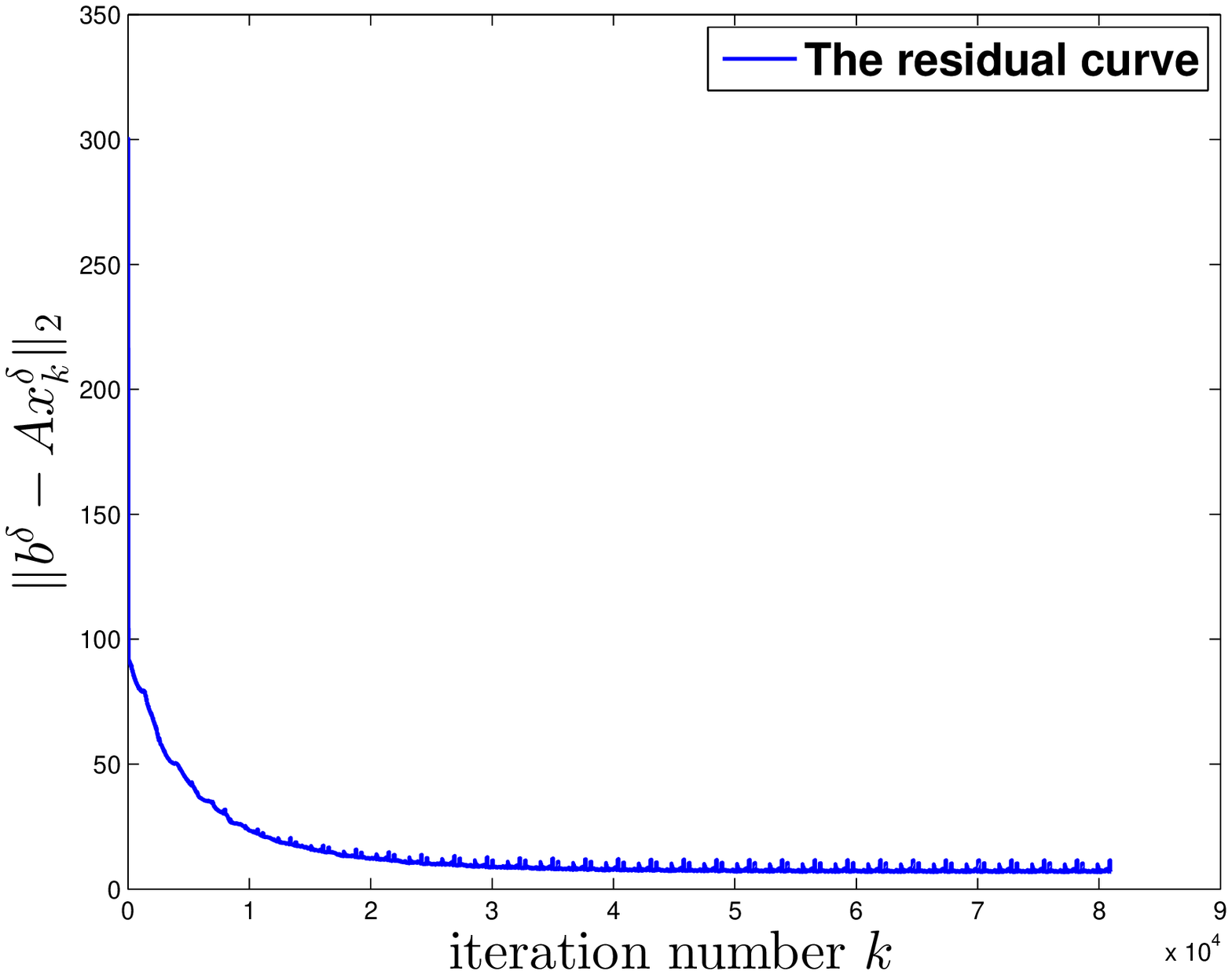}
  \end{minipage}}
    \subfigure[$\eta=0.046$]{
  \begin{minipage}[ht]{.22\linewidth}
      \includegraphics[width=\linewidth]{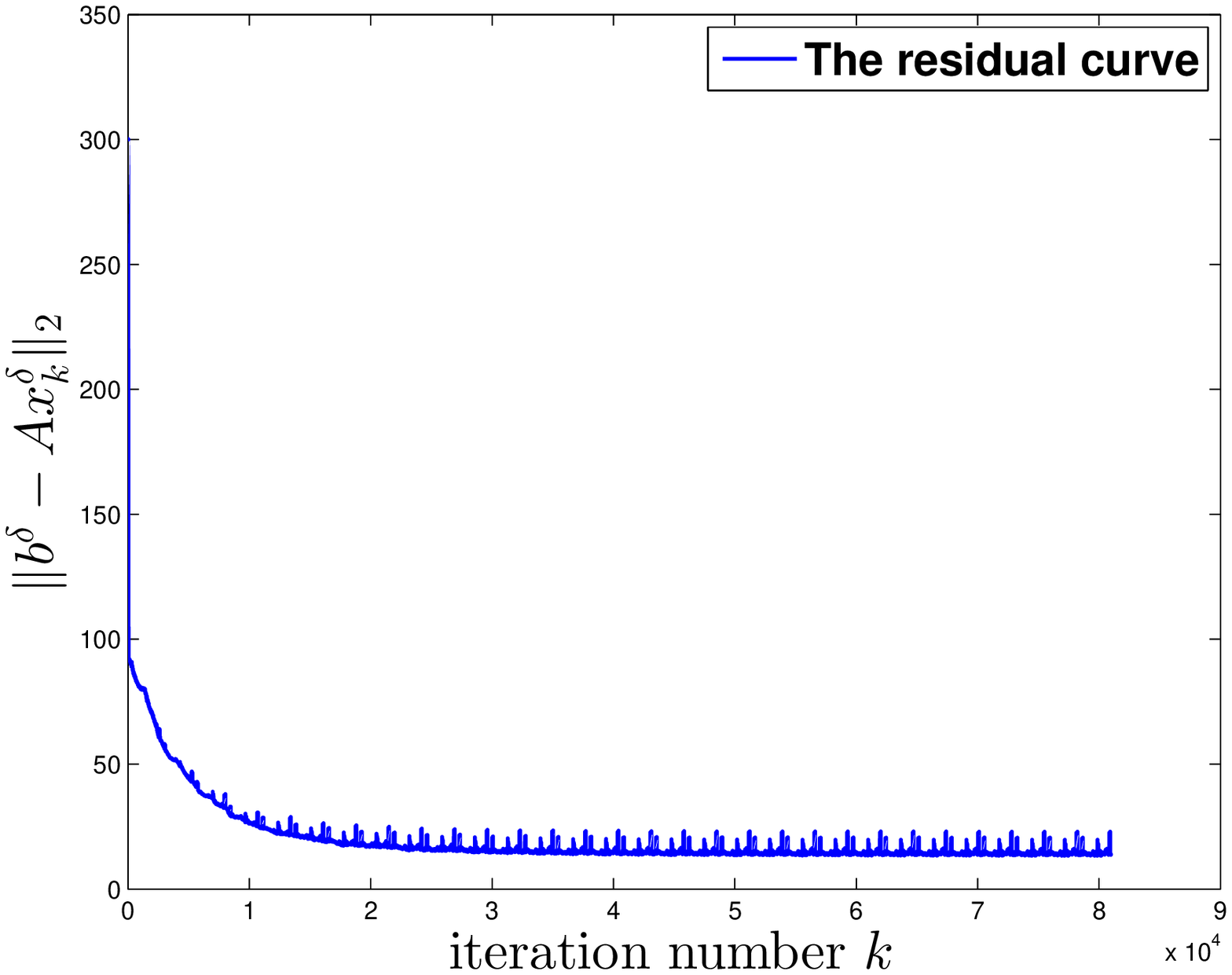}
  \end{minipage}}
  \subfigure[$\eta=0.115$]{
  \begin{minipage}[ht]{.22\linewidth}
      \includegraphics[width=\linewidth]{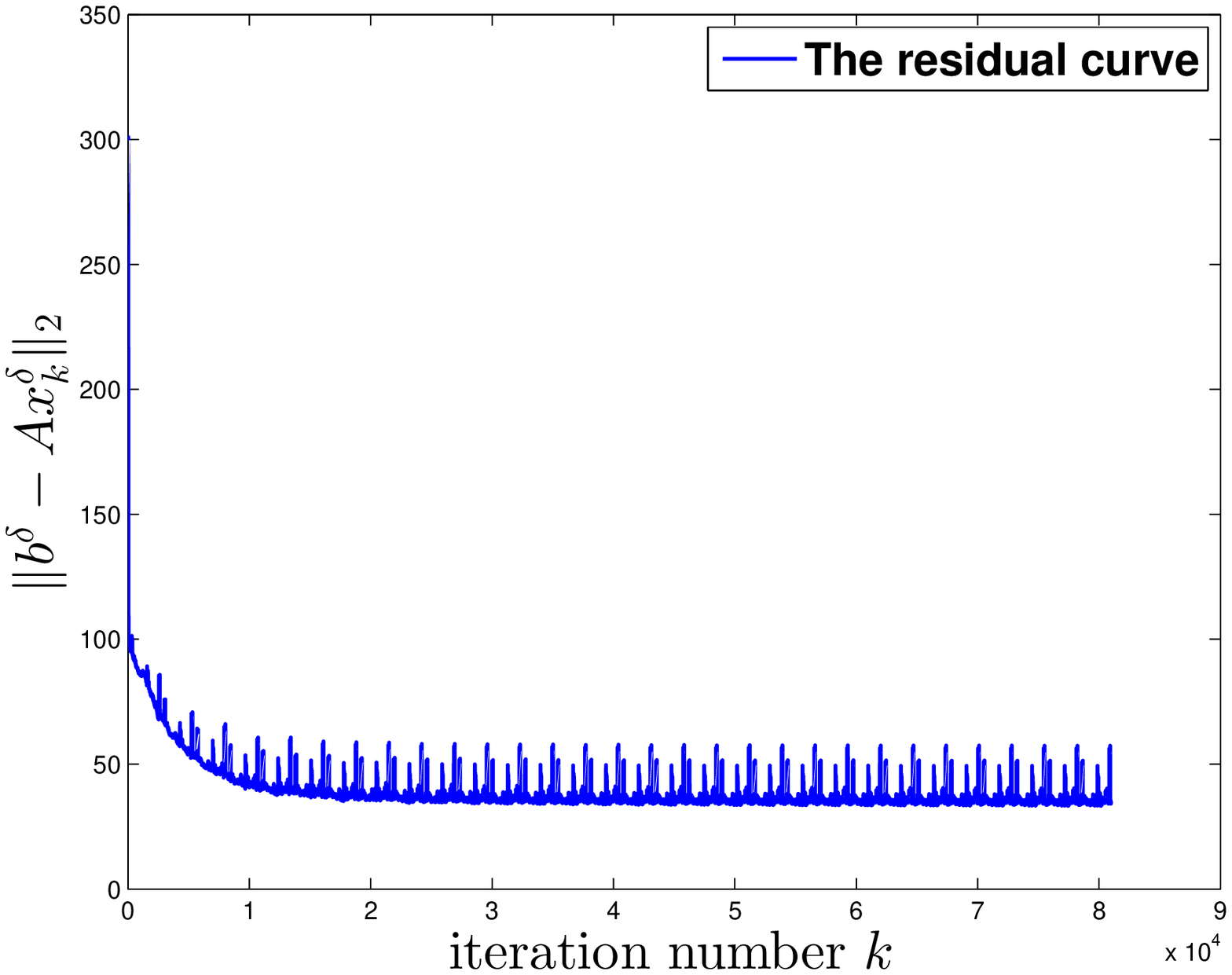}
  \end{minipage}}
  \caption{\footnotesize The residual curves of the Kaczmarz method for Model Problem 3 with Gaussian noise}\label{the_Kaczmarz_method_residue_Gaussian}
\end{figure}
%

\section{Conclusion}\label{section6}

The Kaczmarz-Tanabe method is essentially the Kaczmarz subsequence method. We extract periodically the last iteration of each epoch for Kaczmarz's iteration and construct a new sequence of vectors, which can actually be generated by a matrix multiply a vector. These operations are beneficial to study the Kaczmarz method for linear system as a whole rather than as a standalone equation. Consequently, we can obtain more general results of the Kaczmarz method.

Numerical results in this paper show that the Kaczmarz-Tanabe method (even the Kaczmarz method) is very suitable for a exact linear system, which can also be verified by Theorem \ref{convergence.error.free}. When there exists perturbation in linear system, the amplitude of fluctuation of the iterative error is determined by the second maximum singular value of $Q$ and the minimum non-zero singular value of $A$. We list the maximum and second maximum singular value of $Q$, and the minimum non-zero singular value of $A$ for our numerical problems in Table \ref{singularvalue_of_Q_and_A}.

\begin{table}[!htb]
  \centering
 \caption{The information of the singular values}
  \begin{tabular}{cccc}
    \hline
    \text{Menu}                                      &\text{Model Problem 1} &\text{Model Problem 2} & \text{Model Problem 3} \\
    \hline
    \text{Maximum singular value of $Q$}             & 1.0000                & 0.9913                 & 0.9967 \\
    \text{Second maximum singular value of $Q$}      & 0.7773                & 0.9912                 & 0.9959 \\
    \text{Minimum non-zero singular value of $A$}    & 1.6855                & 8978                   & 0.5455 \\
    \hline
  \end{tabular}
\label{singularvalue_of_Q_and_A}
\end{table}

Tanabe has pointed out that $\|Q\|\le 1$ in \cite{Tanabe1971}, which is in accordance with the maximum singular values in our numerical tests.
But from Lemma \ref{lemma1.1}, there also holds $\|Q\|<1$ in $N(A)^\bot$, which is very important for the convergence rate of the Kaczmarz-Tanabe method and to establish the equivalence between the convergence and the convergence rate for the Kaczmarz-Tanabe method (or the Kaczmarz method). Fortunately, Lemma \ref{Range_result.3} and Corollary \ref{range.result4} help us to restrain $e_k\in N(A)^\bot$ and to guarantee the convergence derived from the result of convergence rate (see Theorem \ref{convergence.error.free}).

\vskip 12pt
\noindent\textbf{Funding information} This work was partially supported by the Natural Science Foundation of Tianjin No. 18JCYBJC88000.


\end{document}